\documentclass[opre,nonblindrev]{informs3Itai}

\usepackage{hyperref} 

\usepackage[sort]{natbib}
 \def\BIBand{and}%
\usepackage{amssymb}
\usepackage{bbm}
\newcommand{\Ex}{\mathbb{E}}
\newcommand{\wtilde}{\widetilde}
\newcommand{\tpi}{\widetilde{\pi}}
\newcommand{\tDelta}{\widetilde{\Delta}}
\newcommand{\tP}{\widetilde{P}}
\newcommand{\tV}{\widetilde{V}}
\newcommand{\tW}{\widetilde{W}}
\newcommand{\rcol}{\textcolor{red}}
\newcommand{\trans}{^{\intercal}}
\newcommand{\1}{\mathbbm{1}}
\newcommand{\bcol}{\textcolor{blue}}
\newcommand{\fraks}{\frak{s}}
\newcommand{\lsb}{\left[}
\newcommand{\rsb}{\right]}
\newcommand{\calB}{\mathcal{B}}
\newcommand{\br}{\mathcal{B}_r}
\newcommand{\brh}{\mathcal{B}_{\varrho}(x)}
\newcommand{\brhh}{\mathcal{B}_{\frac{\varrho}{2}}(x)}
\newcommand{\lpb}{\left(}
\newcommand{\rpb}{\right)}

\newcommand{\tEx}{\wtilde{\Ex}}
\newtheorem{theorem}{Theorem}
\newcommand{\calS}{\mathcal{S}}
\newcommand{\calC}{\mathcal{C}}
\newcommand{\calN}{\mathcal{N}}

\newcommand{\calO}{\mathcal{O}}
\newcommand{\calM}{\mathcal{M}}
\newcommand{\hV}{\widehat{V}}
\newcommand{\hX}{\widehat{X}}

\newtheorem{lemma}{Lemma}[section]
\newtheorem{assumption}{Assumption}[section]

\newtheorem{proposition}{Proposition}[section]
\newtheorem{remark}{Remark}
\newtheorem{example}{Example}

\newcommand{\bsq}{\vrule height .9ex width .8ex depth -.1ex}
\newcommand{\eProof}{\hfill  \Halmos\vspace*{0.4cm} }
\newcommand{\bProof}{\noindent\textbf{Proof: }}

\newcommand{\bg}{G} 
\newcommand{\sg}{g}
\newcommand{\bu}{U}
\newcommand{\su}{u}

\newcommand{\gridfn}{q} 

\newcommand{\moma}{\textsc{MoMa}}

\usepackage{algorithm}
\usepackage[noend]{algpseudocode}

\usepackage[section]{placeins}
\usepackage{mathtools}

\TITLE{A Low-rank Approximation for MDPs\\ via Moment Coupling}

\ARTICLEAUTHORS{%
\AUTHOR{Amy B.Z. Zhang}
\AFF{School of Operations Research and Information Engineering, Cornell University}
\AUTHOR{Itai Gurvich} 
\AFF{Kellogg School of Management, Northwestern University\footnote{This work was supported by NSF grant CMMI-1662294}}}
\RUNTITLE{ADP via Moment Coupling}

\RUNAUTHOR{Amy B.Z. Zhang, Itai Gurvich}

\ABSTRACT{\textbf{Abstract.} We introduce a framework to approximate a Markov Decision Process
that stands on two pillars: state aggregation---as the algorithmic infrastructure; and central-limit-theorem-type approximations---as the mathematical underpinning of optimality guarantees. The theory is grounded in recent work \cite{braverman2018taylor} that relates the solution of the Bellman equation to that of a PDE where, in the spirit of the central limit theorem, the transition matrix is reduced to its local first and second moments. Solving the PDE is {\em not} required by our method. Instead, we construct a ``sister'' (controlled) Markov chain whose two local transition moments are approximately identical with those of the focal chain. Because of this {\em moment matching}, the original chain and its ``sister'' are coupled through the PDE, a coupling that facilitates optimality guarantees. Embedded into standard soft aggregation algorithms, moment matching provided a disciplined mechanism to tune the aggregation and disaggregation probabilities. The computational gains arise from the reduction of the effective state space from $N$ to $N^{\frac{1}{2}+\epsilon}$ is as one might intuitively expect from approximations grounded in the central limit theorem.

} 

\begin{document}\maketitle
\vspace*{-1cm} 

% \noindent \textbf{Keywords}: \\
% dynamic programming -- optimal control: discrete discounted infinite horizon MDP \\
% probability -- Markov processes: moments of Markov chain in PDE \\
% algorithms -- approximation: approximate policy iteration with state aggregation \\
% \newpage
\section{Introduction}

Dynamic programming is the fundamental technique for solving sequential decision problems. The key object of analysis is the Bellman optimality equation. As the dimension of the state space increases, the computational burden of solving the Bellman equation becomes prohibitive. Approximate dynamic programming (ADP) is a family of algorithms developed to address this challenge by reducing---through various mechanisms---the computational complexity of the problem. 

A central ADP theme is that of value-function approximation. Here, one a priori imposes a lower dimensional structure on the value function---assuming, for example, that it is an affine combination of pre-specified basis functions---and optimizes the combination parameters. One expects computational gains if the number of basis functions is small relative to the size of the state space. Even when these {\em architecture-based} algorithms converge, the performance---in terms of optimality gaps---depends on the choice of the basis functions; these are often chosen based on ad-hoc knowledge of the problem's structure. 

The approach we propose here does not impose an architecture on the value function directly. Instead it treats two Markov chains on the same state space as ``sisters'' in terms of their value if, in the spirit of the central limit theorem, their (controlled) transition matrices share the same first and second moments. The value closeness is based on the results---in \cite{braverman2018taylor}---that these two ``close'' chains have a common approximation via (and hence can be related by) a suitable PDE. The optimality gap of using the optimal control derived for one chain for the other, is no larger than the approximation error introduced by the PDE. 

We bring algorithmic relevance to this idea by using state aggregation as the scaffolding on which to build a low-rank ``sister'' chain. Viewed through our lens, aggregation in effect produces an alternative chain on the full state space, but with a lower-rank transition law governed by the aggregation design matrices. Moment matching provides a principled way to choose these design matrices and we introduce a simple and intuitive algorithm to do so: one that is grounded in the mathematical analysis of the approximation gap and, in turn, produces optimality-gap guarantees. 

A central contribution of this work lies in drawing algorithmic relevance from the powerful toolbox of central limit theorem (CLT) approximations for control problems. We produce an appealingly simple algorithm that, under ``CLT-type'' assumptions, results in an optimality gap that is an order-of-magnitude smaller than the value function itself. In our use of state aggregation as the infrastructure, this is also one of the first performance guarantees known for soft aggregation. 

Our CLT-based algorithm does not generally break the curse of dimensionality; one can hardly expect so without the use of ad-hoc properties of specific problems such as linearity of the value or low-rank characteristics of the transition matrix. It does, however, achieve the reduction that is consistent with a CLT approximation --- it reduces the effective number of states from $N$ to $N^{\frac{1}{2}+\epsilon}$ where $\epsilon>0$ can be chosen by the designer, even for problems with full rank. 

This paper provides the conceptual and mathematical foundation for CLT-based aggregation tuning. It is yet another tool in the ADP toolbox and, for specific problems, could be deployed in combination with other tools that exploit their idiosyncratic characteristics.

\subsection{Overview of results}

The lower-rank {\em sister chain} is a ``non-identical twin'' of the original chain coupled to it through the local-transition first and second moments: 
$$\mu(x)=\Ex_x[X_1-x],~~ \sigma^2(x)=\Ex_x [(X_1-x)(X_1-x)^\intercal].$$  These are collapsed statistics of the full transition matrix. Implicit in these definitions is our focus on chains where there is a natural notion of physical distance; we fix attention to state spaces of the form $\mathbb{Z}^d \cap \times_{i=1}^d [\ell_i,u_i]$. 

The premise that coupling two chains via their moments should produce small approximation gaps is grounded in recent work \cite{braverman2018taylor} that connects the Taylor expansion of value function to nearly optimal policies. While the math that supports this statement is non-trivial, the intuition is rather simple. Let $(X_t,~t=1,2,...)$ be a chain $\mathbb{Z}^d$ with transition probability matrix $P$, and consider the infinite horizon $\alpha$-discounted reward $$V(x)=\Ex_x\lsb \sum_{t=0}^{\infty}\alpha^t c(X_t)\rsb.$$ The value $V$ solves the fixed point equation $V(x)= c(x)+ \alpha PV(x),$ which we re-write as $$0=c(x)+\alpha (PV(x)-V(x)) -(1-\alpha)V(x).$$ If we pretend that $V$ has a continuously thrice differentiable extension to the reals $\mathbb{R}^d$, then $$PV(x)-V(x)=\mu(x)'DV(x) + \frac{1}{2}trace(\sigma^2(x)'D^2V(x))
+ \mbox{ Remainder}, $$
where the remainder depends on the third derivative of the continuous extension.  At least intuitively, the fixed-point equation translates to the solution of the partial differential equation
$$c(x)+ \mu(x)'DV(x)+\frac{1}{2}trace(\sigma^2(x)'D^2V(x))=0.$$

We {\em do not} advocate using this PDE as a computational alternative but, rather, as a link (a ``coupling'') between the chain $P$ and a more tractable one. Put simply, any chain $\wtilde{P}$ on $\mathbb{Z}^d$ with the same local moment functions $\mu(\cdot)$ and $\sigma^2(\cdot)$, too, would induce the same PDE. To the extent that the quality of the PDE as an approximation depends only on those moments (as functions over the state space), we have a mechanism to bound the gap between the value of the two chains. Among all sister chains, we want one that is tractable in terms of value-function computation. 

At this point, we plug moment matching as a module into the known {\em aggregation} method in ADP. Aggregation reduces the dimensionality of the Bellman equation by solving it for a small (in relative terms) number of ``meta-states'', denoted by $L$. The extent of the reduction in computational effort depends on how $L$ compares to the number of ``detailed'' states, $N$; the fewer the meta states, the less demanding the computation of the value function. The design parameters of aggregation---the so-called aggregation and disaggregation matrices---are typically chosen in an ad-hoc manner.  {\em Moment matching offers a principled way to choose these that is grounded in approximation/optimality gap bounds. } 

Matching both moments, simultaneously for all states, is generally impossible if one wishes to have $L<N$; see \S \ref{sec:aggregation}; we prioritize the first moment over the second. Our meta states are a grid of spaced out states. These ``representative states'' are the effective support of the chain $\tilde{P}$. We match perfectly the first local moment $\mu(\cdot)$ while maintaining, through {\em non-constant} spacing, a handle on the second-moment mismatch. The coarseness of the reduced state space is directly informed by the mathematical analysis. Interestingly, our moment matching gives mathematical justification to what is a rather intuitive choice of an aggregation matrix; see Theorem \ref{thm:backtodelta}. 

We prove with $L=\calO(N^{\frac{1-\varepsilon}{2}})$ meta states, the value of the sister chain $\tV$, is a good approximation to the true value $V$ in the following sense: 
\be |V(x)-\tV(x)| =\calO\left(\Ex_x\left[\sum_{t=0}^{\infty}\alpha^t\frac{c(X_t)}{(1+\|X_t\|)^{\varepsilon}}\right]\right)=o\left(V(x)\right),\label{eq:guaranteeintro}\ee where $\varepsilon\in (0,1)$ is a design variable. The closer it is to $0$, the fewer meta states (so computation is easier) but the larger the gap bound. 
% The gap is proportional to the infinite-horizon discounted value with a scaled down cost function $c(x)/(1+\|x\|)^{\varepsilon}$; with $\varepsilon=0$ computation is easier but the gap is of the order of the value itself.

When we embed moment matching in an approximate policy iteration algorithm the computational gains are further magnified with savings in not only policy evaluation but also the update steps. Importantly, our moment-based design of the aggregation and disaggregation matrices is {\em  policy independent}. In turn, they are computed once and do not have to be updated on each iteration. 

These accuracy guarantees are not fully general; to use PDE theory, we require that the first and second transition moments to satisfy certain smoothness properties (see Assumption \ref{asum:primitives}). These reveal the connection to the central limit theorem. To shed some light on these assumptions, while avoiding for now the formality it requires, consider the simple(st) case where the Markov chain is just a random walk: $X_t = \sum_{s=1}^t Z_s$ where $Z_s^t,~s=1,\ldots,t,$ are i.i.d random variables. Then, $\widehat{X}_t= \frac{X_t}{\sqrt{t}} = \sum_{s=1}^t \frac{Z_s}{\sqrt{t}}$ is well approximated by a normal for large $t$; $\widehat{X}_t$ is, as well, a random walk with increments of size $=\calO(1/\sqrt{t})$. What we require is a CLT-type relationship between the effective horizon $1/(1-\alpha)$ and the chain's speed---namely that the increments are inversely proportional to the square root of the horizon:  $$\mu(x)\sqrt{\frac{1}{1-\alpha}} \sim 1.$$ Put otherwise, we require that the effective horizon is not too long relative to the chain's ``speed''. Arguably, if this is violated ---if the chain moves at fast speeds, the horizon is enough to reach stationarity and an ergodic formulation may be more suitable.

\subsection{Literature\label{sec:lit}}

%ADP is concerned with approximating solutions to complex control problems where the size of the state space prohibits exact computation of the value function and/or the optimal control policy. The literature on ADP is vast. 

%In ADP, gains in computation are often achieved by restricting the search for value functions to an {\em architecture}---a pre-specified family of functions. In linear architectures, for example, value functions are restricted to linear combinations of pre-specified features. More recent methods use neural networks as the underlying architecture; see e.g. \cite{bertsekas2018feature} and \cite{vanvuchelen2020use,gijsbrechts2019can} for recent applications to operations management problems. 

Architecture-based ADP has a long and successful history. With more recent methods using neural networks as the underlying architecture; see e.g. \cite{bertsekas2018feature} and \cite{vanvuchelen2020use,gijsbrechts2019can} for recent applications to operations management problems. Two questions must be posed to any architecture-based ADP algorithm: (1) Does the algorithm converge to the best choice of parameters {\em within the given architecture}. In the case of a linear architecture, for example, does the algorithm produce the best feature coefficients; (2) Such convergence may not mean much if the architecture is inadequate for the problem at hand, so we must also ask how well the ``best'' choice within the given architecture approximates the original problem of interest. 

The first question was answered affirmatively for linear architectures in \cite{tsitsiklis1996feature,tsitsiklis1997analysis}, and followed by improvements to convergence rates; see e.g. \cite{devraj2017fastest}. There are, however, few approximation algorithms with general theoretical guarantees on the optimality gaps--- namely, on how well the prescribed (approximate) control performs in the original system; see e.g. \cite{pires2016policy}. Furthermore, while approximate dynamic programming has known significant practical success, the choice of the architecture often builds on ad-hoc intuition about the problem at hand, rather than on a principled approach to its construction. 

Our focus is not on convergence rates for a given approximation architecture but, rather, on a new architecture with optimality-gap guarantees. Instead of using value function approximation, we approximate the controlled Markov chain by matching its local moments. We produce an algorithm by piggy backing on the state-aggregation approach: given the original (controlled) chain we build a new chain that with the same {\em transition} moments, through the choice of aggregation parameters; it is a principled approach to tune the {\em aggregation design variables}.

Conceptually speaking, our sister chain---constructed, as it is, through aggregation---is a ``linear model approximation'' and as such, has an intimate connection to linear architectures; see \cite{parr2008analysis,jin2020provably,mahadevan2009learning, yao2014pseudo}. We here start from the model, rather than the value function. We create a linear approximation for {\em for the model}; moment matching guides this approximation's design. We subsequently draw implications for the value function. 

State aggregation has a long history; see \cite{whitt1978approximations, bean1987aggregation, tsitsiklis1996feature} to name a few. We primarily follow the exposition in \cite{bertsekas2017dynamic}. Various forms of aggregation and/or abstraction have been considered in the literature, see e.g. \cite{li2006towards,abel2016near,deng2011optimal}; aggregation in the context of discretizing a continuous state space is considered in e.g. \cite{chow1991optimal,kushner2001numerical}. Some of the literature considers aggregation that is learned from samples and adapted during the algorithm's progression, see e.g. \cite{bertsekas1989adaptive, baras2000learning,sinclair2019adaptive} and the literature review contained.

\iffalse 
Much of the literature on aggregation and abstraction focuses on hard aggregation, where each state has membership in only one cluster. The clusters are usually determined based on some notion of lumpability; see \cite{li2006towards} for a unified survey on exact abstraction and \cite{abel2016near} for approximate abstraction. Spectral analysis is another popular approach to partitioning, see e.g. \cite{deng2011optimal} and the literature referenced. Many works also focus on adaptive aggregation with unknown transitions, updating memberships based on algorithm progression, e.g. \cite{bertsekas1989adaptive, baras2000learning}. Of related flavor is the active area of adaptive discretization, see e.g. \cite{sinclair2019adaptive} and the literature review within. Many demonstrate impressive empirical success, but either provide little theoretical guarantees, or focus mainly on the efficiency of learning from samples, especially in the model-free setting. Continuous space discretization can also be interpreted through a similar lens, but available guarantees are often regarding convergence in the limit of the number of partitions; for classical analysis see e.g. \cite{kushner2001numerical,chow1991optimal}
\fi

We rely on soft aggregation and determine the parameters (the aggregation and dis-aggregation matrices) in advance. The ability for a state to have membership in more than one cluster in soft aggregation provides flexibility that is generally useful (see \cite{singh1995reinforcement, sarich2010optimal}) and turns out particularly suitable for our needs. Appealingly our choice of aggregation and disaggregation parameters does not depend on the transition matrix $P$ but only on the state space and desired level of accuracy.

\iffalse 
The power of soft aggregation, where a state can have membership in more than one cluster, has been recognized in works such as \cite{singh1995reinforcement}, but largely with only convergence results. Of similar idea is fuzzy partitioning, which considers probabilistic memberships in partitions; see \cite{sarich2010optimal} for an example and more references. However fundamentally, these work are concerned with exploring different metrics and orient their results as such. In contrast, what matters to us is the effect on the value function from the approximate transition matrix, which leads to the moment matching criteria and the guarantees directly on value approximation quality.
\fi 

Our implicit assumption is that the transition matrix $P$ is full rank. The computational gain then arises from constructing the lower-rank ``sister'' chain. In an application context, it makes sense to test that implicit assumption and, if the original chain itself is low-rank, to first use matrix factorization techniques to identify the aggregation parameters; see e.g.  \cite{duan2019state, ghasemi2020identifying,zhang2019spectral}. 
% pires2016policy
More generally, it makes sense to take advantage of latent structure such as low-rank $Q$ or $P$, particularly in learning; \cite{agarwal2020flambe, arumugam2020randomized}. %\bcol{This awareness of value function is seen in recent works such as \cite{farahmand2017value}.}

% Also related to our work is \cite{parr2008analysis,jin2020provably} that relate linear-model approximations to linear value-function approximations. 

Our work---in so far as it shows that two chains with similar transition moments produce similar values---can inform this literature in suggesting the use of these moments as a value-function-relevant metric to guide the choice of factorization.

\iffalse 
In so far as model compression can be approached from linear approximation, there is significant work on basis construction and representation learning; see e.g. \cite{mahadevan2009learning}. On the other hand, if a kernel is (approximately) linearly decomposible, the action-value function is also (approximately) linearly decomposible under some boundedness conditions, as shown in Proposition 2.3 in \cite{jin2020provably}. Many sophisticated methods have been carefully designed to address the numerical delicacy associated with stochastic factorization (or for pseudo-MDPs, nonnegative matrix factorization); see e.g. \cite{zhang2019spectral} and other references within. In terms of the effect on value function accuracy by using a low-rank approximation of the transition kernel, recent papers such as \cite{pires2016policy} are making progress towards quantifying it. We have seen limited existing work on taking this a step further by presenting procedure to minimize this error, and provide a more informative bound that results from such; indeed our grounding in moment coupling allows us to go that extra step, which we see also as opportunity for future fusion with more numerically efficient  methods.

\fi 

Moment-based approximations---inspired by the central limit theorem and functional version thereof---have been extremely successful in queueing theory facilitating the  analysis and optimization of highly complex queueing networks. Some of the ``import'' of the mathematical theory from the control of queues to general dynamic programs has been achieved in \cite{braverman2018taylor} where the connections to queueing theory are thoroughly discussed. 

We use the mathematical constructs in \cite{braverman2018taylor} as a starting point for an algorithmic framework. What we adopt is the view that matching local moments---a collapsed ``statistic'' of the full transition matrix---has the potential to produce small optimality gaps. How to do so algorithmically --- how to construct the sister chain $\tP$ for computational gains --- is the question we address in the current paper. In the process of developing our algorithm, we expand on \cite{braverman2018taylor} to allow for some mismatch in the second moment between the focal chain and its sister in our bounds.

The classical moment problem in probability has a long history; see \cite{prekopa1990discrete} and the references therein. Our challenges here deviate from the classical moment problem. Most fundamentally, we are facing a {\em simultaneous} problem as we are trying to match the first two moment of all of $N=|\mathcal{S}|$ random variables --- one for each state, where the random variable for state $x$ has the distribution $P_{x\cdot}$---via convex combination of the {\em same} (small) set of random variables. This, as will be made evident through simple examples, is generally impossible. 

Finally, the value-function effect of replacing $P$ with $\tP$, is an instance of sensitivity analysis for MDPs; see e.g. \cite{ross2009sensitivity, mastin2012loss} and, in model-based reinforcement learning, \cite{sun2018dual,janner2019trust}. \vspace*{0.2cm} 

\noindent {\bf Organization.} 
This paper is organized as follows. In \S \ref{sec:themodel} and \S \ref{sec:tayloring}, we study moment matching in the setting of Markov reward process (control is fixed). Being simpler than the control problem, this serves expositional clarity. In \S \ref{sec:aggregation}, we connect aggregation to a sister chain -- a low rank model on the same state space. \S \ref{sec:grid}, contains our explicit design of the aggregation and dis-aggregation matrices. The formal optimality guarantees appear in \S \ref{sec:guarantees}. 

We bring all of this to bear on optimization (optimal control) in \S \ref{sec:optimization}. Empirical performance for two numerical examples are presented in \S \ref{sec:numerical}. All proofs for lemmas appear in the appendix.

\section{The model\label{sec:themodel}}

We consider the infinite-horizon discounted reward for a discrete-time Markov chain on a finite state space $\calS\subseteq  \mathbb{Z}^d \cap \times_{i=1}^d[\ell_i, u_i]$. 
Let $N=|\calS|$ be the size of the state space. $P$ is the transition matrix with $p_{xy}$ equal to the probability of transitioning from $x$ to $y$ in one step; $c:\mathcal{S}\to \mathbb{R}^+$ is the cost function. We assume that the function $c$ is norm like; that there is a $k\in \mathbb{Z}_+$ and a point $x_0\in\calS$ such that \be \frac{1}{\Gamma} \|x-x_0\|^k\leq |c(x)|\leq \Gamma \left(1+\|x-x_0\|\right)^k.\label{eq:normlike} \ee Since one can shift the state space, we assume w.l.o.g. that $x_0=0$. Finally, $\alpha \in (0,1)$ is the discount factor. This so-called ``Markov reward'' process is characterized by the tuple $\mathcal{C}=<\mathcal{S},P,c,\alpha>$. 

The value function is then given by 
$$V(x) =\Ex_x\left[ \sum_{t=0} ^{\infty} \alpha ^ t c(X_t)\right],~x\in \mathcal{S},$$ where $\Ex_x[\cdot]$ is the expectation with respect to the law $P_{x\cdot}$. 

For a function $f:\mathcal{S}\to \mathbb{R}$ we use the operator notation $Pf(x):=(Pf)(x)=\sum_{y}p_{xy}f(y)=\Ex_x[f(X_1)]$. As is standard, the function $V:\calS\to \mathbb{R}$ is the unique solution to the equation $T V = V$, where
\[ T V(x)  = c(x) +\alpha P V(x). \] We refer to this as the Bellman equation despite the absence of a control decision here. This allows for continuity of language with optimization in \S \ref{sec:optimization}. Since the state space is finite, $V$ can be computed via the matrix inversion formula $V=(I-\alpha P)^{-1}c.$

In our analysis we sometimes refer to the maximal jump size of $P$ from $x$  \be \Delta_x:=\sup_{y: p_{xy}>0}\|y-x\|. \tag{maximal jump}\label{eq:maxjump}\ee  
Because our state space is finite, this quantity is bounded. 

\vspace*{0.2cm} 
\noindent {\bf Notation.} Unless stated otherwise, $\|\cdot\|$ corresponds to the Euclidean norm on $\mathbb{R}^d$ ($d$ will be clear from the context). We  write $y=x\pm \epsilon$ to denote $\|y-x\|\leq \epsilon$. We use $\mathbb{R}_+^d$ and $\mathbb{Z}_+^d$ to denote the non-negative reals in $\mathbb{R}^d$ and integers in $\mathbb{Z}^d$, and use $\mathbb{R}_{++}^d$ and $\mathbb{Z}_{++}^d$ when they are strictly positive. For a function $f:\mathcal{A}\to \mathbb{R}^d$ and a set $\mathcal{B}\subseteq \mathcal{A}$,  $|f|_{\mathcal{B}}^*=\sup_{x \in\calB}\|f(x)\|$. We use $\Gamma$ to denote a universal constant whose value might change from one line to the next but that does not depend on the state $x$ or the discount factor $\alpha$. Where useful we point out its dependencies. We write $f(x)\lesssim \gridfn(x)$ to mean $f(x)\leq \Gamma \gridfn(x)$ and $f(x)\cong \gridfn(x)$ if both $f(x)\lesssim \gridfn(x)$ and $\gridfn(x)\lesssim f(x)$.  \vspace*{0.2cm} 

\section{Tayloring reconsidered\label{sec:tayloring}} 

Consider two Markov Reward Processes. The first, $\mathcal{C} = <\mathcal{S}, P, c, \alpha>$,  is driven by the {\em focal} chain $P$. The other, $\wtilde{\mathcal{C}} = <\mathcal{S}, \wtilde{P}, c, \alpha>$, is driven by the {\em sister} chain $\wtilde{P}$; $\wtilde{\mathcal{C}}$ differs from $\mathcal{C}$ only in terms of the transition probability matrix. 

%If $\|P-\wtilde{P}\|_1\leq \epsilon$ then it is easy to show that $|\wtilde{V}-V|_{\calS}^*\leq \frac{\alpha\epsilon |c|_{\calS}^*}{(1-\alpha)^2};$ see e.g. \cite{ross2009sensitivity}. The value function itself has $|V|_{\calS}^*\leq \frac{|c|_{\calS}^*}{(1-\alpha)}$, so this bound is valuable only if $\epsilon=o(1-\alpha)$ (as $\alpha$ approaches one); also $|c|_{\calS}^*$ itself might be large. Moreover, 
A ``replacement'' of a chain with a proxy is useful only insofar as it yields computational benefits by, say, being of lower rank. It seems ambitious to require $P$ and a lower rank $\wtilde{P}$ to be close in some reasonable matrix norm unless $P$ is itself low rank. Instead, it makes sense to measure the distance between transition matrices in terms of their impact on the value function.

\begin{equation}
    |V- \wtilde{V}|_{\calS}^*\leq \frac{\alpha}{1-\alpha}\left(\delta_{V}[P,\wtilde{P}]+\delta_{\widetilde{V}}[P,\wtilde{P}]\right)
    \label{eq:tilde_operator_gap1}
\end{equation} where, for a function $f:\calS\to \mathbb{R}$, \begin{align*} \delta_f[P,\wtilde{P}]:=\sup_{x\in\calS}\lvert   \wtilde{P}{f}(x) - Pf(x)\rvert &=\sup_{x\in\calS}\lvert \tEx_x[f(X_1)]-\Ex_x[f(X_1)]\rvert \end{align*} 
\iffalse  

For transition probability matrices $P,\widetilde{P}$ and function $f:\mathcal{S}\to \mathbb{R}$ we have that 

we define the metric:
$$\delta_f[P,\wtilde{P}](x):=   \lvert \wtilde \Ex_x[f(X_1)]-\Ex_x[f(X_1)]\rvert = \lvert   \wtilde{P}{f}(x) - Pf(x)\rvert,$$
and let $ \delta_f[P,\wtilde{P}] :=  |\delta_f[P,\wtilde{P}] (\cdot)|_{\calS}^*$. Then
\fi 
This bound is reminiscent of classical results. It seems, however, of limited value as it requires information about $V$, the very construct whose computation we seek to avoid. It is nevertheless useful in that it identifies the one-step-ahead expectation $|\Ex_x[V(X_1)]-\tEx_x[V(X_1)]|$ as a central object of study. It makes clear that, in comparing two chains, what matters is the {\em local behavior}: how the single-step-change in value under $P$, $\Ex_x[V(X_1)]-V(x)$, compares to that under $\tP$, $\wtilde \Ex_x[V(X_1)]-V(x)$.

This localization makes Taylor-expansion a natural lens through which to consider approximation gaps. We make the following observation. If $V$ has a thrice continuously differentiable extension to $\mathbb{R}^d$, then \be \Ex_x[V(X_1)]= V(x) + \mu(x)'DV(x)+\frac{1}{2}trace(\sigma^2(x)'D^2V(x))\pm \frac{1}{6}\|D^3V\|\Delta_x^3,\label{eq:tayloring}\ee 
where, recall, $\Delta_x:=\sup_{y: p_{xy}>0}\|y-x\|$ is the maximal jump of the chain from state $x$, and $\mu,\sigma^2$ are the local moments 
$$\mu(x)=\Ex_x[X_1-x],~~ \sigma^2(x)=\Ex_x [(X_1-x)(X_1-x)^\intercal].$$ 

The expectation $\tEx_x[V(X_1)]$ for the sister chain is expanded analogously. If $\tP$ and $P$ share these moments, i.e., for all $x\in\calS$ 
$$\wtilde{\mu}(x):=\tEx_x[X_1-x]\approx \mu(x) \mbox{  and } \wtilde{\sigma}^2(x) :=\tEx[(X_1-x)(X_1-x)^{\intercal}]\approx
\sigma^2(x),\mbox{ for all }x\in \calS,$$ then 
\begin{align*} 
\tEx_x[V(X_1)]& \approx
V(x) + \wtilde{\mu}(x)'DV(x)+\frac{1}{2}trace(\wtilde{\sigma}^2(x)'D^2V(x))\\
&\approx V(x) + \mu(x)'DV(x)+\frac{1}{2}trace(\sigma^2(x)'D^2V(x)) \approx \Ex_x[V(X_1)],
\end{align*} 
so that 
$$\delta_V[P,\tP]= \sup_{x\in\calS}\mid \Ex_x[V(X_1)]- \tEx_x[V(X_1)]\mid \approx0.$$ Here $\approx 0$ should be interpreted as ``$\delta_V$ being small relative to the value function $V$''; the precise mathematical meaning of $\approx 0$ is exposed in \S \ref{sec:guarantees}. This informal derivation makes clear that (1) if a low-rank sister chain has the same moments as the focal chain, its value may provide a good approximation to that of the focal chain. In an effort to be low rank, the sister chain might have larger jumps, so that (2) in designing this sister chain we must keep its jumps small, at least in regions of the state space where the third derivative is substantial.

\begin{example}[The simple random walk] {\em Consider the simple absorbing random walk on the integers:  $P_{x,x+1}=P_{x,x-1}=1/2$ for all $x=1,...,n-1$ and  $P_{00}=P_{nn}=1$. Then $\Ex_x[X_t]=x$ for all $t\geq 0$ so that $\Ex_x[\sum_{t=0}^{\infty} \alpha^t X_t]=\frac{x}{1-\alpha}$. The same conclusion holds for the  ``simpler'' chain that jumps in one step to one of the end points: $\wtilde{P}_{xn}=1-\wtilde{P}_{x0}=x/n$. Because $\mu(x)=\wtilde{\mu}(x)=0$ for all $x$ $P$ shares the local first moment with a sister chain $\tP$ that has only two states. } \hfill \bsq \label{example:simple} \end{example} 

Example \ref{example:simple} is rather unique. One should not expect a perfect value-function match in general, certainly not with such a coarse state-space. Our bounds in \S \ref{sec:guarantees} will capture the dependence of the approximation's accuracy on the ``density'' of the meta-states. 

The informal derivation through Taylor expansion is useful for developing intuition but does not provide a basis for algorithm design. The value $V$ is not apriori known so it is impossible to ``refer'' to its continuous extension. To circumvent this, \cite{braverman2018taylor} develops a framework for obtaining {\em indirectly} an approximate continuous solution. A short summary of this earlier work is useful. Consider a chain on $\mathbb{Z}^d$. The value $V$ solves the Bellman equation $V(x)=c(x)+\alpha PV(x)$ which we find useful to re-write as
\be \label{eq:PDE} 0 = c(x) + \alpha (PV(x)-V(x))-(1-\alpha)V(x).\ee Pretending that the function $V$ is twice continuously differentiable, 2nd-order Taylor expansion yields the partial differential equation (PDE) 
$$ 0= c(x)+\alpha\left[\mu(x)'DV(x)+\frac{1}{2}trace(\sigma^2(x)'D^2V(x))\right]-(1-\alpha)V(x),$$ defined now over $\mathbb{R}^d$.   
While this equation has been arrived-to purely formally, the following is a valid mathematical question: what is the relationship between a solution $\hV$ (if it exists) to this equation on $\mathbb{R}^d$, and $V$ that solves the original discrete-state-space Bellman equation.

Two chains $<\calS,P,c,\alpha>$ and $<\calS,\tP,c,\alpha>$ with the same local moment functions $\mu(\cdot)$ and $\sigma^2(\cdot)$ induce the same reduction to a continuous-state space PDE so that bounds  $|V-\widehat{V}|$ and $|\wtilde{V}-\widehat{V}|$ produce, as a corollary, a bound on $|V-\wtilde{V}|$. This is the path we take. 

\section{Sister-chain construction via aggregation \label{sec:aggregation}}

Aggregation effectively creates a new Markov chain on a smaller state space. The tuning of the aggregation parameters is tantamount to selecting for this chain a transition matrix from a restricted family of such. The flexibility this offers makes it an ideal vehicle for our moment-matching algorithm. 
\subsection{Aggregation Preliminaries}

Recall that $N = \mid \calS \mid$ denotes the number of states in the original MDP, and let $\calM=\{1,\ldots,L\}$ be a set of meta states; obviously $L\leq N$. We refer to \cite[Chapter 6]{bertsekas2012approximate} for a thorough introduction to  aggregation and include below the minimal ingredients for a self-contained exposition. Two weight matrices govern the mapping between $\mathcal{S}$ and $\calM$;  
\begin{itemize} 
\item {\em Aggregation probabilities $(\bf G)$: } For each detailed state $x\in \mathcal{S}$, the probability that $x$ aggregates (or ``groups'') into $k$, $\sg_{xk}\geq 0$, represents the degree of membership of detailed state $x$ in meta-state $k\in\calM$. The $N \times L$ matrix $G=\{g_{xk}\}$ is non-negative and row-stochastic. {\em Hard aggregation} is the special case where the meta-states form a partition of the state-space, and each state $x\in\calS$ ``belongs'' to a single partition: $\sg_{xk}=1$ for one and only one $k\in\calM$. The more general case is referred to as {\em soft aggregation}.

\item {\em Disaggregation probabilities $(\bf U)$:} For each meta-state $l\in\calM$, the probability that $l$ disaggregates (or ``un-groups'') into $x$, $\su_{lx}\geq 0$, is the degree to which meta-state $l$ is represented by detailed state $x \in \calS$; the $L\times N$ matrix $\bu=\{\su_{lx}\}$ is non-negative and row-stochastic. If a $l\in\calM$ is represented by a single state $x_l$, i.e. $\su_{lx_l} = 1$, we refer to this $x_l$ as the {\em representative state} of meta-state $l$. It is then convenient to think of $\calM$ as the \textit{set of representative states} $\calS^0:=\{x\in\calS: x_l=x \mbox{ for some } l\in\calM\}$. \end{itemize} 

Having fixed matrices $G$ and $U$, one solves an {\em aggregated} Bellman equation on the meta-states: 
\begin{align}\label{eq:aggregateBellman}
R(l)& =\sum_{x\in\mathcal{S}}\su_{lx} (c(x)+ \alpha \sum_{y\in \mathcal{S}}p_{xy}  \sum_{k\in \calM}\sg_{yk}R(k)), ~l\in \calM, 
\end{align}  whose matrix form $R = \bu c+\alpha \bu P\bg R$  reduces to \be  R = (I-\alpha \bu P\bg)^{-1}\bu c.\label{eq:agginverse}\ee
The function $R$ is the value function of the aggregate problem. 

In hard aggregation, the true value function is assumed to be constant over each subset in the partition. We say that $x\in S_k$ (or in ``cluster'' $k$) if $\sg_{xk}=1$, and approximate its value with the aggregate value $R(k)$. Most works in the literature study this setting.

In the case of soft aggregation (stochastic G), but combined with a representative state for each meta-state, detailed states are interpolated from the representative ones. This is sometimes referred to as {\em coarse grid} scheme, which we opt for in our algorithm as explained below.

% Another special case is when each meta-state $l$ has a representative state $x_l$ (i.e. $\bu$ is binary); there  When combined with soft aggregation (non-binary $\bg$), this is sometimes referred to as {\em coarse grid} scheme, and interpolates detailed states from the representative ones. This is what we use in our algorithm, for reasons explained later.

\subsection{A low rank chain on $\calS$\label{sec:aggregation_connection}}

Aggregation produces a Markov chain on the meta states $\calM$ with transition law $\bu P\bg$: one first transitions (via $\bu$) to $\calS$, then within $\calS$ (via $P$) and finally back to $\calM$ (via $G$). Theorem \ref{thm:lifted} makes a simple but powerful observation that permuting these steps produces a chain on $\calS$ with closely related value. 
\begin{theorem} 
[aggregation as sister chain]
     Consider the value $\tV$ of a Markov chain on the original detailed state space $\mathbb{S}$ with the transition matrix 
$$\wtilde{P}=P\bg \bu~~~(\widetilde{P}_{xy} = \sum_{z\in \calS,l\in\calM}p_{xz}\sg_{zl}\su_{ly}).$$
The aggregate value $R$ in \eqref{eq:aggregateBellman} equals $\bu\tV$, and $\tV=c+\alpha P\bg R$. 
\label{thm:lifted} 
\end{theorem} 

\bProof 
From the Bellman equation for this chain, we have that the value $\widetilde{V}$ satisfies 
\begin{align*} \widetilde{V}(x)&=c(x) + \alpha\widetilde{P}\widetilde{V}(x)= c(x) +\alpha P\bg \bu \widetilde{V}(x).
\end{align*} 
Defining $\widetilde R:=\bu \widetilde V$ we have  $\tV=c+\alpha P\bg \widetilde R$. Multiplying both sides by $\bu$ gives 
$\widetilde R = \bu c + \alpha \bu P\bg\widetilde R. $. This $\widetilde R$ is in fact the unique solution to the aggregate Bellman equation.  \eProof

In this way, aggregation gives rise to a family of lower rank chains with law $\wtilde P[\bg,\bu] := P \bg \bu$, which we call the \textit{$(\bg,\bu)$-\emph{lifted chain}}, on the detailed space. Note that these values are still obtained from the lower dimensional $R$ and thus enjoy the reduction  in computational complexity. 

Theorem \ref{thm:lifted} provides a concrete mechanism to construct, from $P$, a lower-rank chain on the same state space $\calS$. One can then tune the parameters $\bg,\bu$ to meet certain objectives. We tune these for moment matching, i.e., so that for all $x\in \calS$, the first and second moments under this lifted chain given by
\begin{align*}  \wtilde \mu[\bg,\bu](x)&=\tEx_x[X_1]=\sum_{y} \wtilde{P}_{xy}[\bg,\bu] (y-x)\mbox{, }\\ 
\wtilde \sigma^2[\bg,\bu](x)&=\tEx_x[(X_1-x)(X_1-x)^{\intercal}]=\sum_{y}\wtilde{P}_{xy}[\bg,\bu](y-x)(y-x)^\intercal  \end{align*}
are matched to $\mu(x),\sigma^2(x)$ of the original transition law $P$. If there is first moments match, i.e., 
\begin{align*}
   \mu(x) - \wtilde{\mu}[\bg,\bu](x) =& \Ex_x[X_1]- \tEx_x[X_1]=0,  \end{align*} then the second moment gap reduces to 
\begin{align*}
     \sigma^2(x) - \wtilde{\sigma}^2[\bg,\bu](x)  = \Ex_x[X_1X_1^{\intercal}]- \tEx_x[X_1X_1^{\intercal}].
\end{align*} 

\subsection{Design considerations}
One can, at this point, set up a direct optimization to choose $G,U$ that minimize a norm of the gaps $\mu-\wtilde{\mu}$ and $\sigma^2-\wtilde{\sigma}^2$. 
% This is related to long-studied  moment problem in probability; see \cite{prekopa1990discrete} and the references therein. \iffalse and the references therein. There is a fundamental deviation from the classical moment problem. Most fundamentally, we are facing a {\em simultaneous} problem as we are trying to match the first two moment of all of $N=|\mathcal{S}|$ random variables --- one for each state, where the random variable for state $x$ has the distribution $P_{x\cdot}$---via convex combination of the {\em same} (small) set of random variables. \fi 
We avoid the subtleties of choosing the correct norm as well the computational cost of solving such an optimization problem. We take, instead, an approach that gives rise to a design of $G,U$ that is, appealingly enough, independent of the specific matrix $P$.  

As noted above, one simplifying structure uses hard aggregation---where $\bg$ is taken to be binary---but may be too restrictive, related as it is to assuming that the values are similar across all states in a partition. We use a non-binary (stochastic) $\bg$, but a binary $\bu$. \iffalse As can be seen from \S \ref{sec:lit}, $\bg$ is taken to be binary (i.e. hard aggregation / partitioning) in much of literature, which provides a venue for significantly reducing the complexity. However, this scheme performs well only if the value function is close to piece-wise constant, as each partition is approximated with a single value. It is also an intuitively restrictive assumption, and ill-suited when one might wish to represent states as having membership in more than one cluster, which is frequently the case in high-dimensional problems. \fi This turns out to provide sufficient flexibility while having better interpretability and allowing us to avoid full policy optimization on $\calS$ during policy iteration; see \S \ref{sec:optimization}. 

Specifically, for each meta state $l$, the disaggregation distribution $\bu_{l,\cdot}$ has 0's everywhere except for the state $x_l$ that ``represents'' the meta state $l$. \iffalse For this reason this scheme is often referred to as aggregation ``with representative states'' or ``a coarse grid scheme''.\fi The set \iffalse $$\calS^0=\{x\in\calS: \bu_{l,x}=1 \mbox{ for some } l\in\calM\}$$\fi $\calS^0$ of representative states has a one-to-one correspondence with $\calM$. Moment matching then means that, for all $x\in\calS$,
\begin{align}
    & \Ex_x[X_1] = \tEx_x[X_1] = \sum_{l \in \calM} [P\bg]_{xl} x_l \tag{1st-moment match} \\
    & \Ex_x[X_1X_1^{\intercal}] = \tEx_x[X_1X_1^{\intercal}]=\sum_{y,z \in \calS, l \in \calM} [P\bg]_{xl} x_l x_l^{\intercal} \tag{2nd-moment match}
\end{align}

Perfect matching \iffalse with representatives states \fi thus necessitates that any point in the $n+n^2$ dimensional scatter $\{\Ex_x[X_1],\Ex_x[X_1X_1^{\intercal}]),x\in\calS\}$ can be written as a convex combination of $\{ x_l, x_l x_l^{\intercal}, l \in \calM \}$; achieving this with $L<N$ is generally impossible; see Example \ref{example:impossibility} below. 

 \begin{example}[The (im)possibility of 2nd moment matching] {\em
A simple example makes it abundantly clear that $L$ might have to be no smaller than $|\mathcal{S}|$ for perfect matching of both moments. Consider the simple absorbing random walk on $\{0,1,\ldots,n\}$, with $P_{x,x+1}=P_{x,x-1}=1/2$ for all $x=\{1,\ldots,n-1\}$ and $\{ 0, n \}$ are absorbing states. Here we have $\Ex_x[X_1]=x$ for all $x$ and $\Ex_x[X_1^2]=x^2+\1\{x\notin \{0,n\}\}$. Given the scatter $\{(\Ex_x[X_1],\Ex_x[X_1^2]),x\in\mathcal{S}\}$, one cannot express all points as a convex combination of a small number of (common) points; see the round markers in Figure \ref{fig:simpleRW}. A piecewise linear approximation allows for matching the first moment while controlling, through the number of breakpoints, the quality of second moment match. }

{\em For contrast consider a chain that has $P_{x0}=1-x/n$ and $P_{xn}=x/n$ (absorbing in one step at the boundary). Here $\Ex_x[X_1]=x$ for all $x$ and $\Ex_x[X_1^2]=nx$ and both moments can be matched using only two representative states corresponding to the end/corner points of the state-space; see the square markers in Figure \ref{fig:simpleRW}.   } 

\begin{figure}[h]\centering
\includegraphics[scale=0.3]{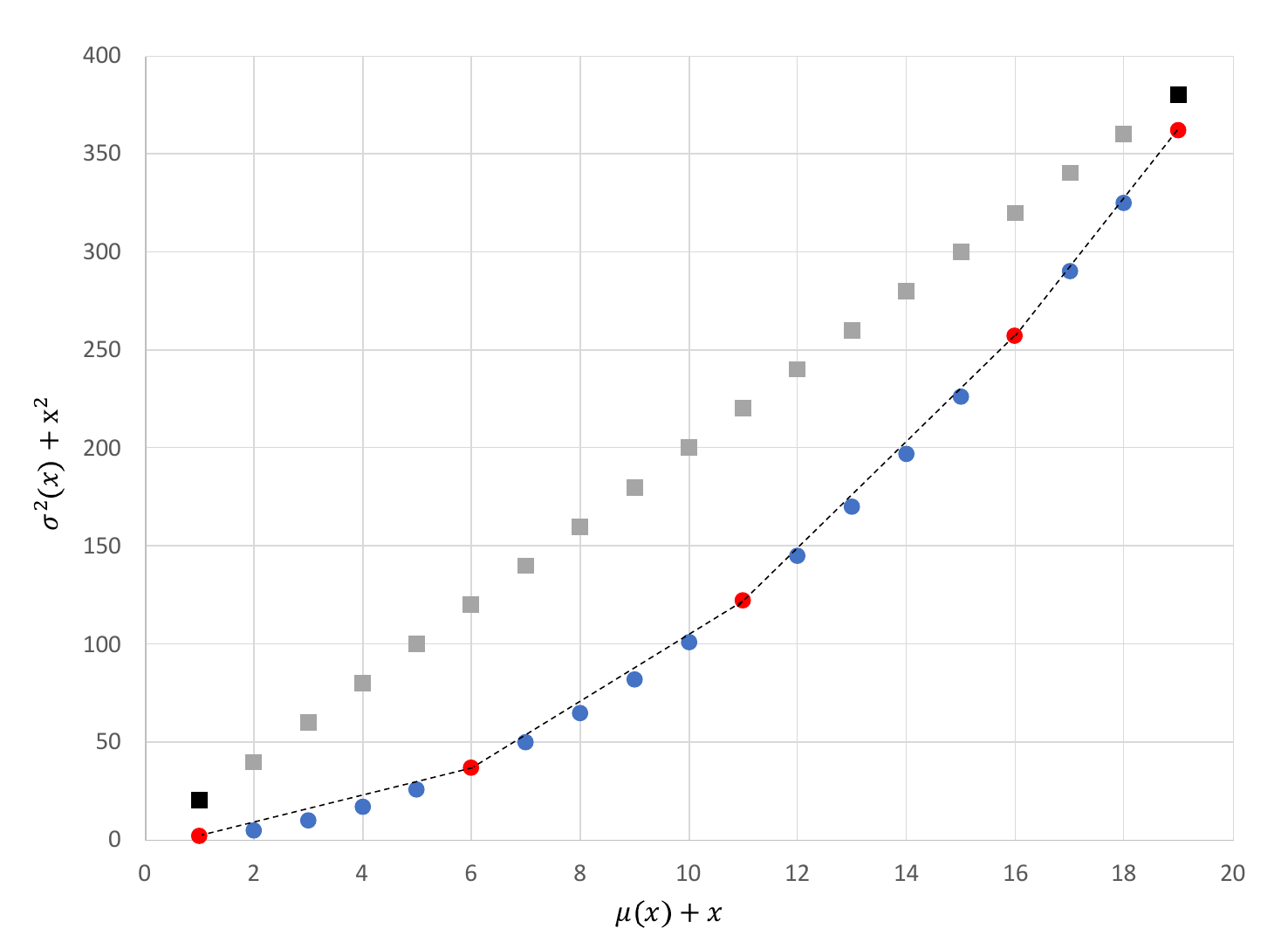}
     \caption{The moment scatter plot for two different random walks on $\{0,1, \ldots,20\}$. Circles: Simple random walk with absorbing end points. Simultaneous matching of both moment is impossible here---one cannot express a point as a convex combination of other points. Squares: a random walk where each point in the moment scatter can be written as a convex combination of the two end points.\label{fig:simpleRW}}\end{figure}
\label{example:impossibility}\hfill \bsq \end{example}

Given the general impossibility of matching both moments, we prioritize the matching of the first moment. The intuition of Tayloring, as reflected in \eqref{eq:tayloring}, is informative: errors in matching $\mu$ should translate into approximation error that are proportional to the first derivative of $V$, whereas errors in the matching of $\sigma^2$ would only be multiplied by the second derivative. We will insist then on matching the first moment exactly while controlling the second moment mismatch.

\section{Moment-Matching (MoMa) Aggregation Design\label{sec:grid}}

\subsection{Designing $\bg, \bu$}

Recall that we take a binary $\bu$. We design $\bg$ for first moment matching and doing so will be exceedingly simple. 

\begin{lemma}
For perfect first moment matching between $\wtilde P[\bg,\bu]$ and $P$, it suffices that for each $y$, $(\bg \bu)_{y,\cdot}$ is the distribution of $y+Z$ where $\Ex[Z]=0$.
\label{lem:convolution}
\end{lemma}

We will refer to a design of $G,U$ that has these properties as a {\em zero-mean} construction. A matrix $\bg$ with rows $\{\sg_{y,\cdot},y\in\calS\}$ satisfies the condition of Lemma \ref{lem:convolution} if for each state $y\in\calS$, $\sg_{y,\cdot}$ is a distribution over $\calM$ 
such that \be \sum_{l}\sg_{yl}x_l=y.\label{eq:convcomb}\ee Indeed, then
$$\wtilde{\Ex}_x[X_1]=\sum_{y}p_{xy}\sum_{x_l}\sg_{yx_l}x_l = \sum_{y}p_{xy} y=\Ex_x[X_1], \mbox{ for all } x \in \calS. $$ 

It remains to choose the set of representative states $\calS^0$ (and then $\bg$) so that \eqref{eq:convcomb} holds. \vspace*{0.2cm} 

\noindent {\bf Constructing $\bu$.} 
Given a grid over the state space, each point on the grid is taken to be meta-state $l \in \calM$, and the original state at that point is the representative state $x_l$. Define $\mathcal{S}^0 = \{x_l: l \in \mathcal{M} \}$. Construct $\bu$ to be a $L \times N$ matrix. In the $l^{th}$ row, assign 1 to $\su_{lx_l}$, and 0 otherwise.

\begin{figure}[h!]\centering
\includegraphics[width=0.4\textwidth]{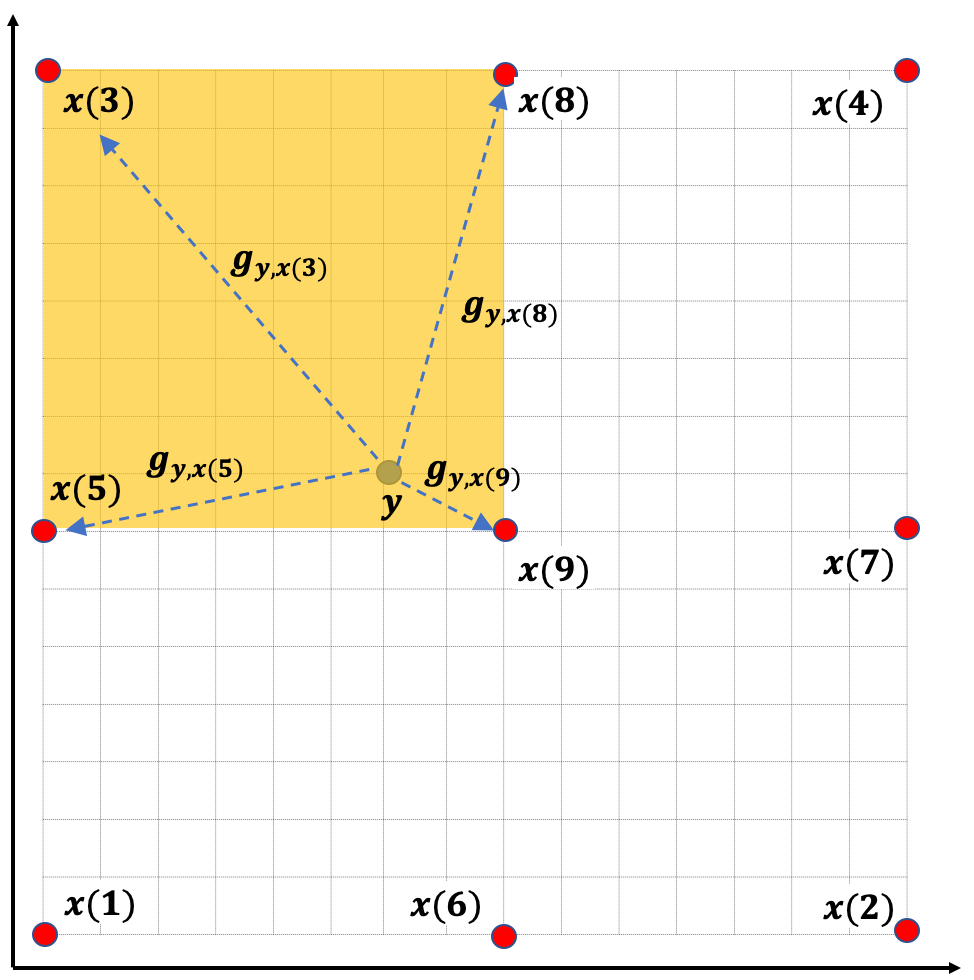} \hspace*{1cm} 
\includegraphics[width=0.4\textwidth]{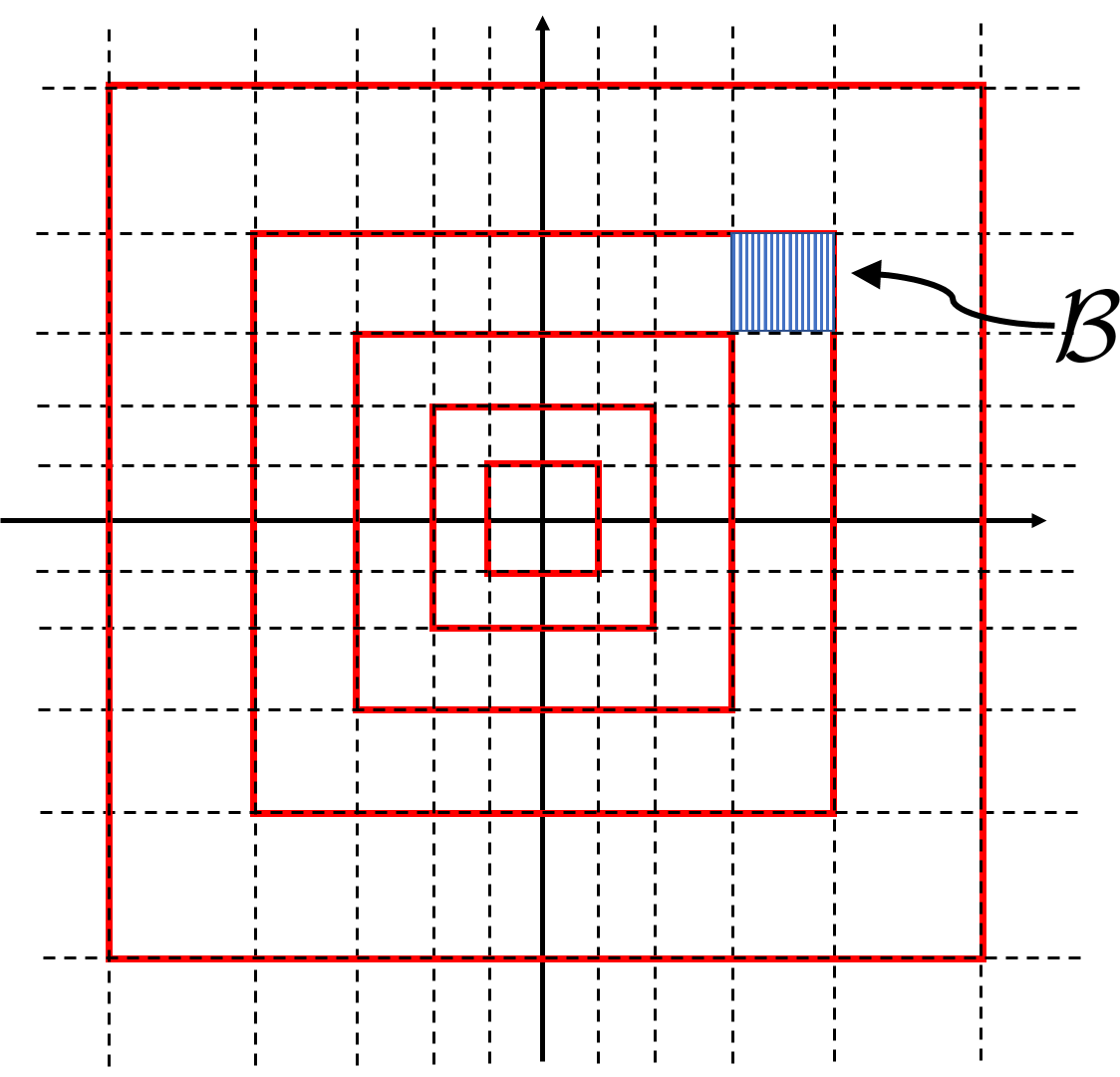}
\caption{(LEFT) Aggregation matrix $\bg$ is used to express each state $y$ as a convex combination of the closest points on the grid. (RIGHT) The grid with spacing exponent $\fraks=0.5$ and one box $\calB$ within it. \label{fig:grid}}
\end{figure}

The spacing in this grid will be judiciously chosen in \S \ref{sec:grid} below but, regardless of this choice, it is obvious that each $y \in \calS$ can be represented as a convex combination of the vertices of its encasing hypercube $\calB(y)$ simply by taking a linear interpolation, as illustrated in Figure \ref{fig:grid} (LEFT). \vspace*{0.2cm} 

\noindent {\bf Constructing $\bg$.} 
For a given $y\in\calS$, let $\calB = \{ \vec{k}_1, ..., \vec{k}_{2^d} \} \in \calM$ be the set of $2^d$ meta (i.e. representative) states that form its encasing box, restrict $\sg_{yl'} = 0$ for $l' \notin \calB$, and solve for 
$$
    \sum_{l \in \calB}\sg_{yl}x_l=y, \mbox{ where } \sum_{l \in \calB}\sg_{yl}=1 \mbox{ and } \sg_{yl} \geq 0 \mbox{ for } l \in \calB.
$$

This set of linear constraints has an explicit solution. Given a box $\calB$, for $i \in [d]$, let  $\bar{s}_i=\max_{x\in\calB}x_i$ and $\underline{s}_i=\min_{x\in\calB}x_i$. Then, given $y$ and its enclosing box $\calB$, we write 
\be \sg_{yl} = \Pi_{i=1}^d \left[ \1\{(x_l)_i = \bar{s}_i \} * \frac{y_i -  \underline{s}_i}{\bar{s}_i-\underline{s}_i} + \1\{ (x_l)_i = \underline{s}_i\} * \frac{\bar{s}_i - y_i}{\bar{s}_i-\underline{s}_i}  \right] \label{eq:phiconstruction}.\ee

Intuitively, $\sg_{yl}$ weighs nearby representative states $x_l \in \calB$ proportional to their relative distance to state $y$. The following summarizes the properties of $\bg$, verifying that it is a valid aggregating matrix, and induces perfect moment matching with the grid-based $\bu$.

\begin{lemma} \label{lem:phi_construct}
The construction of $\bg$ in \eqref{eq:phiconstruction} satisfies $\sum_{l \in \calB} \sg_{yl} = 1$ and $\sum_{l \in \calB}\sg_{yl}x_l=y$. Also, $\sg_{yl} = 1$ when $y = x_l$  for $l \in \calM$. 
\end{lemma}

\subsection{Explicit construction of the grid}

Denote the value of $V$ (/ $\tV$) restricted to $\calS^0$ with $W$ (/ $\tW$), so that
$$\bg W = \sum_{l}\sg_{yl}V(x_l) = \bg \bu V, ~~ \bg \tW = \sum_{l}\sg_{yl}\tV(x_l) = \bg \bu \tV.$$

\begin{theorem}
 Under a coarse grid scheme , the approximation error between the exact value $V$ and the approximate value $\tV$ from the $(\bg,\bu)$-lifted chain satisfies
 \begin{align} |V-\tV|_{\calS}^*&\leq \frac{1}{1-\alpha}\left(|V-\bg W|_{\calS}^* +|\tV-\bg \tW|_{\calS}^*\right).\label{eq:theorem2} \end{align}  
  \label{thm:backtodelta}
\end{theorem}

This generalizes a known result for hard aggregration. Indeed, when further restricting to a binary $G$, we obtain \cite[Proposition 4.2]{bertsekas2018feature}: 
\be \label{eq:recover_bound}   |V(x)-R(k)| \leq \frac{|\epsilon|_{\calM}^*}{1- \alpha}, ~ k\in \calM,~ x\in S_k,\ee where 
$\epsilon(k) = \max_{x,y\in S_k}  | V(x) - V(y) |.$ 
  This simple reduction is formally proved in the appendix.
 
 Theorem \ref{thm:backtodelta} provides a useful vehicle to motivate our detailed design choices. Pretending that a smooth extension of $V$ exists,
$$V(y)- [\bg W] (y)\approx -\sum_{l}\sg_{yl}DV(y)'(x_l-y) + \calO(\Delta_y^2 D^2V(y))= \calO(\Delta_y^2 D^2V(y)),$$ where the equality follows from $\sum_{l}g_{yl}x_l=y$; similarly, $$\tV(y)- [\bg \tW] (y)\approx \calO(\widetilde\Delta_y^2 D^2\tV(y)).$$

We make this more formal in \S \ref{sec:guarantees}, where $D^2V, D^2\tV$ will both be replaced by $D^2\hV$, the second derivative of the PDE solution $\hV$ to \eqref{eq:PDE}. 

It is clear then that our choice of grid should be such that the product of $\widetilde \Delta_y$ and the derivative $D^2\hV$ is, in a suitable sense, small. Theorem \ref{thm:guaranteemain} says that choosing the spacing so that 
\be \tDelta_y\leq (1-\alpha)^{\frac{1}{4}}\|y\|^{\frac{1-\varepsilon}{2}},\label{eq:jumpreq}\ee
guarantees the optimality gap states which in informal terms is $$ |V(x)-\tV(x)| =\calO\left(\Ex_x\left[\sum_{t=0}^{\infty}\alpha^t\frac{c(X_t)}{(1+\|X_t\|)^{\varepsilon}}\right]\right),$$
for design variable $\varepsilon \in (0,1)$.
\iffalse 

This is stated more formally in the following theorem. Define  $$V_{-\varepsilon}[x] = \Ex_x\lsb \sum_{t=0}^\infty \alpha^t \frac{|c(X_t)|}{(1+\|X_t\|)^{\varepsilon}}\rsb.$$

\begin{theorem}[approximation gap] \label{thm:guaranteemain} Suppose assumptions stated in \ref{asum:primitives} holds, and that $\widetilde{\Delta}_x\lesssim (1+(1-\alpha)^{\frac{1}{4}}\|x\|^{\frac{1-\varepsilon}{2}})$. Then, given $\varepsilon \in (0,1)$, for any $\kappa\geq 2 +\varepsilon - \frac{k}{2}$ and all $x:\|x\|\geq (1-\alpha)^{-(1+\kappa)}$, \be |\tV(x)-V(x)|\lesssim  
\frac{1}{\sqrt{1-\alpha}}(V_{-1}(x)+\tV_{-1}(x)+
V_{-\varepsilon}(x)+\tV_{-\varepsilon}(x)).\ee The $\lesssim$ here does not depend on $x,\alpha$. 
\end{theorem}

Then we need to design the grid so that $\widetilde{\Delta}_x$ is bounded relative to $\lVert x \rVert$. Then it is natural to consider spacing that is dependent on the state, and in fact on that increases with the magnitude. Define an increasing \textit{spacing function} $\gridfn(z): \mathbb{R} \mapsto \mathbb{R}_+$, such that on  axis $i \in [d]$, the space between the $k^{th}$ point on the grid with coordinate value $f_i(k)$ and the $(k+1)^{th}$ point is at least $\gridfn(|f_i(k)|)$.
\fi

To save on computation it makes sense to use the spacing in the grid as large as possible to minimize computation costs while maintaining the guarantee. 

To that end, recall that $\calS=\mathbb{Z}^d\cap \times_{i=1}^d [l_i,u_i]$ and consider first the positive portion of an axis $i$. Let the coordinate of the first gridpoint be $f_i(0) = \max \{ 0, \ell_i \}$. Let the axis value at index $k$ be given recursively by 
\be f_i(k+1) = \lceil f_i(k)+ (f_i(k))^{\fraks} \rceil + 1,\label{eq:f_recursion} \ee
where the spacing between grid point $f_i(k)$ and $f_i(k+1)$ is determined by a \textit{spacing function} $\gridfn(f_i(k))$, chosen here to be $(f_i(k))^{\fraks}$ with $0 \leq \fraks < \frac{1}{2}$.
The largest index on axis $i$ is $\bar{n}_i:=\min\{k:f_i(k) \geq u_i\}$, where we set $f_i(\bar{n}_i)=u_i$. Construct symmetrically for the negative axis the set of gridpoints $\{f_i(-k)\}$. 

This construction guarantees that one-step transitions in the sister chain---which jumps from a state $x$ to the corners of the box that contains some $y$ that is reachable from $x$---satisfy the requirement \eqref{eq:jumpreq}. \footnote{In the mathematical guarantees we use the spacing function $q^{\alpha}(z)=(1-\alpha)^{\frac{1}{4}}q(z)$. For $\alpha=0.99$ for example, $(1-\alpha)^{\frac{1}{4}}\geq 0.3$. We simplify the algorithm exposition by dropping this multiplicative constant.}

\begin{lemma}
The choice of $\gridfn(z) = z^{\fraks}$, $0 \leq \fraks < \frac{1}{2}$, for the grid construction and $\bg, \bu$ design described above, produces a sister chain $\tP[\bg, \bu]$ that satisfies $\tDelta_x\leq (1-\alpha)^{\frac{1}{4}}\|x\|^{\frac{1-\varepsilon}{2}}$ for $\varepsilon = 1 - 2 \fraks$.
\label{lem:gridfn}
\end{lemma}

Each point on the grid is thus characterized by an index set $\vec{k} = [k_1, ..., k_d]$, where $k_i\in [-\underline{n}_i,\bar{n}_i]$, $i=1, ..., d$. {\em These grid-points are the representative states}: 
$$ \calS^0 = \{ x({\vec{k}}): x({\vec{k}})_i = f(k_i) \}.$$
For notational simplicity, when the explicit value of $\vec{k}$ is immaterial we revert to using $l$ and $x_l$ for a meta state and its representative state.  

% In terms of spacing, recall that the goal is to use spacing in the grid that is as large as possible, while ensuring the maximal one-step jump in the sister chain is bounded. The lemma below details why a spacing of $(f_i(k))^{\fraks}$ between grid point $f_i(k)$ and $f_i(k+1)$ is suitable.

\iffalse 
As discussed above, we need $\fraks < \frac{1}{2}$ to satisfy Theorem \ref{thm:guaranteemain}, which informally, guarantees that
$$ |V(x)-\tV(x)| =\calO\left(\Ex_x\left[\sum_{t=0}^{\infty}\alpha^t\frac{c(X_t)}{(1+\|X_t\|)^{\varepsilon}}\right]\right)=o\left(V(x)\right),$$ 
where $\varepsilon = 1 - 2 \fraks$.\fi

Complexity analysis will reveal that it is ideal to take $\fraks \geq \frac{1}{3}$ for efficiency. So the desirable range for the spacing exponent is $[\frac{1}{3}, \frac{1}{2})$. Figure \ref{fig:grid} (RIGHT) illustrates the general pattern over $\mathbb{Z}^2$, with the red lines highlighting how the spacing along each axis scales with the distance to the origin on that axis.  The following lemma shows that \iffalse with a binary $\bu$ corresponding to a grid constructed using spacing function $\gridfn(z) = z^{\fraks}$, and $\bg$ that achieves perfect first moment matching,\fi with the above construction, the second moment $\wtilde \sigma^2(x)$ indeed has bounded difference from $\sigma^2(x)$.
\begin{lemma}[second moment mismatch]\label{lem:secondapprox} Consider a Markov chain on $\calS=[\ell_i,u_i]^d\cap \mathbb{Z}^d$ with second transition moment of $\sigma^2(\cdot)$. For a $(\bg,\bu)$-lifted chain with $\bu$ corresponding to the grid constructed above and $\bg$ obtained using zero-mean construction that achieves first moment matching, its second transition moment $\wtilde \sigma^2(\cdot)$ satisfies $$\|\wtilde \sigma^2(x) - \sigma^2(x) \|\leq  \Gamma \|x\|^{2\fraks}$$ for some constant $\Gamma$ that depends only on $d$.
\end{lemma}

\subsection{Moment-Matching (\moma) Value Approximation}

Putting the three steps together, we have the \emph{Moment-Matching (\moma) Pre-processing} procedure with coarse grid scheme, which guarantees exact first moment matching and bounded second moment gap; this is summarized in Algorithm \ref{alg:moma}. 
\begin{algorithm}
	\floatname{algorithm}{Algorithm}\caption{Moment-Matching (\moma) Pre-processing}
	\renewcommand{\thealgorithm}{} \label{alg:moma}
	\begin{algorithmic}[1]
		\Require State space $\calS=\mathbb{Z}^d\cap \times_{i=1}^d [l_i,u_i]$, spacing exponent $\fraks \in [\frac{1}{3}, \frac{1}{2})$.
		\Ensure Aggregation structure with parameters $\bu, \bg$. 
		\State {\em Construct $\fraks$-spaced grid}: Create grid-points $\{f_i(k)\}, \{f_i(-k)\}$ as appropriate for each axis $i$.
		\State  {\em Construct \bu}: For each meta-state $\vec{k}= [k_1, ..., k_d]$ on the grid, assign $\su_{\vec{k}x_{\vec{k}}}=1$ for $[x_{\vec{k}}]_i = f_i(k_i)$.
		\State {\em Construct $\bg$}: For each state $y$, compute distribution $\sg_{y\cdot}$ using (\ref{eq:phiconstruction}).
	\end{algorithmic}
\end{algorithm}

The lemma below characterizes the number of meta-states produced.

\begin{lemma}[number of meta-states]
With spacing exponent $\fraks$, the number $L=|\calS^0|=|\calM|$ of representative (and hence meta-) states satisfies  
\[L\leq  \left( \frac{\sqrt{2}}{1-\fraks} \right) ^d |\calS|^{1- \fraks}.\] 
\label{lem:grid_growth}
\end{lemma}

Then the number of meta-states is bounded by
$(2\sqrt{2})^d N^{1-\fraks}$ for $\fraks\in [0,1/2]$. In the special case where  $\calS=[0,r]^d$, 
$$\frac{L}{N} = \frac{|\calM|}{|\calS|} \leq  \left(\frac{2\sqrt{2}}{r^{\fraks}}\right)^d,$$
implying that the bigger $r- (2\sqrt{2})^{\frac{1}{\fraks}}$ is, the more substantial the dimensionality reduction.

This construction mechanism can then be plugged into a value approximation algorithm using aggregation, which approximates the value of a Markov rewards process, or equivalently, a Markov decision process with a given policy; see Algorithm \ref{alg:eval}.

\begin{algorithm}
	\floatname{algorithm}{Algorithm}\caption{Policy evaluation with \moma~aggregation}
	\renewcommand{\thealgorithm}{} \label{alg:eval}
	\begin{algorithmic}[1]
		\Require Markov reward process $\mathcal{C}= <\mathcal{S}, P, c, \alpha>$, spacing exponent $\fraks \in [\frac{1}{3}, \frac{1}{2})$.
		\Ensure Approximate value $\tV$. 
		\State {\em \moma~Pre-processing}: Obtain $\bu, \bg$ using Algorithm \ref{alg:moma}.
		\State  Solve $R=(I-\alpha \bu P\bg)^{-1}\bu c$.
		\State Compute approximation $\tV=c+\alpha P\bg R$.
	\end{algorithmic}
\end{algorithm}

\subsection{Computational complexity \label{sec:complexity}}

The computational complexity of Markov Decision Problems (MDP) is well studied; see \cite{littman2013complexity}, \cite{blondel2000survey} for a detailed exposition. The discussion here focuses on the evaluation step. We embed it in the context of policy optimization in \S \ref{sec:optimization}. 

Recall that $N=|\calS|$ is the number of states and $L=|\calM|$ is the number of meta states. Per Lemma \ref{lem:grid_growth}, we know $L=\calO(N^{1- \fraks})$. Value of $\fraks$ closer to $1$ are less expensive but more inaccurate. Let range $r$ be smallest integer such that $u_i - \ell_i \leq r,\mbox{ for all } i\in [d]$, we have also that $N \leq (r+1)^d$ and, in turn, that $L=\calO(r^{d(1-\fraks)})$.

 Two ingredients determine the computational value of our approach. The first is the {\em gain} from matrix inversion. This is a gain that is embedded in aggregation and is independent of moment matching. The second is the {\em loss} inherent to our moment-based computation of the aggregation matrices $\bg, \bu$. We treat these two ingredients separately.\vspace*{0.1cm} 
 
 \noindent \paragraph{\bf Matrix inversion.} Computationally speaking, the key step in solving for $V=(I-\alpha P)^{-1}c$ is the inversion of the $N\times N$ matrix $(I-\alpha P)$. The complexity of matrix inversion is  $\Omega(N^2\log(N))$ (see \cite{tveit2003complexity}) but $\mathcal{O}(N^3)$ is achieved by the standard Gauss-Seidel inversion. Solving for the aggregate value $R \in \mathbb{R}^L$, on the other hand, requires the inversion of the smaller $L\times L$ matrix $(I-\alpha \bu P\bg)$ so that
$$\mbox{Gain} = \Omega(N^2\log N)- \calO(L^3).$$ 

When $\fraks \geq \frac{1}{3}$, $L^3=\calO(N^{3(1- \fraks)})= \calO(N^2)$, so we have a gain of  
$$\Omega(N^2 \log N-N^{2})=\Omega(N^2\log N).$$ 

% \textcolor{red}{for gain we do not do $\calO$, only lower bound $\Omega$ so remove that. Then where and why does the $N^2\log N$ disappear and become linear. What is the value of bringing in $r$ here. As is this is confusing. }

This is a conservative estimate of the gain. On one hand, no known algorithm achieves the $N^2\log N$ lower bound and, on the other, various algorithms are faster than Gauss-Seidel and require less than $\calO(L^3)$ for the aggregate problem. If we fix the inversion algorithm (to, say, Gauss-Seidel inversion) the aggregate matrix inversion takes $\calO (N^{3- 3\fraks})$ compared to $\calO (N^3)$ for the full one, approximately square root the time complexity when $\fraks$ is close to $\frac{1}{2}$.\vspace*{0.1cm} 

\noindent \paragraph{\bf Moment matching.} The matrix $\bg$ can be constructed as a linear program\footnote{It would have, per state $y$, $2^d$ variables (as the number of box corners) and $d+1$ constraints (one constraint for each dimension $i\in[d]$ and an additional stochasticity constraint)}. Leveraging our coarse grid scheme, we construct $\bg$ explicitly in \eqref{eq:phiconstruction}. These operations 
take $\calO(Nd2^d)$ time. This is compared against the gain of at least $\calO(N^2 log N)$ with $\fraks \geq \frac{1}{3}$. The total gain is then 
$$\Omega(N^2\log N-Nd2^d) .$$ 

\section{Policy optimization\label{sec:optimization}  }

Some control notation is needed first. We let $\mathrm{A}(x)$ be the set of feasible controls in state $x\in\calS$. We use the notation $\pi$ for a stationary policy; it is a function from $\calS$ to $\mathrm{A}:= \cup_{x \in \calS}\mathrm{A}(x)$ such that $\pi(x)$ is the action the policy takes in state $x$. Let $p_{xy}^a$ denote the probability of transitioning from $x$ to $y$ under the action $a \in \mathrm{A}(x)$, and $P^{\pi}$ for the transition matrix under policy $\pi$; $\Ex_x^a$ (or respectively $\Ex^{\pi}$) is the corresponding expectation.

The Bellman operator for a fixed policy $\pi$ is given by 
$$T^{\pi} V(x)=c(x,\pi(x))+\alpha [P^{\pi}V](x),$$  so that the value under $\pi$ is the solution to the fixed point equation $V^{\pi}=T^{\pi}V^{\pi}$ which is solved by matrix inversion; recall \S \ref{sec:themodel}. The optimization Bellman operator $T$ is given by
$$T V(x)=\max_{u \in \mathrm{A}(x)}\{c(x,a)+\alpha [P^aV](x)\},$$
and the optimal value $V^*$ is the unique solution of the Bellman optimality equation $V^* = T V^*$. Denote the minimizing policy by $\pi^*$; if multiple policies are optimal, we arbitrarily pick one. 

The first and second local moments depend on the state and the action taken in that state:
$$\mu_a(x)=\Ex_x^a[X_1-x],\mbox{ and } \sigma_a^2(x) = \Ex_x^a[(X_1-x)(X_1-x)\trans],~x\in\calS,$$ and we denote with $\wtilde{\cdot}$ all analogous definitions for a sister chain.

The optimal aggregate value function is the fixed point of 
$$R(k)=\sum_{x\in\mathcal{S}}\su_{kx} {\min_{a \in\mathrm{A}(x)}} \sum_{y\in \mathcal{S}}p_{xy}^a [c(x,a)+\alpha \sum_{l\in \calM}\sg_{yl}R(l)].$$
Although the value is defined only for $k\in \calM$, the minimizing policy, note, is defined on the full state space $\calS$. We measure the performance of the approximate policy $\pi'$ by comparing its value ($V^{\pi'}$) to the optimal value  ($V^*$): $\mid V^*(x) - V^{\pi'}(x) \mid $ is the \emph{optimality gap}.

\subsection{Approximate PI with \moma~aggregation}

The bound for a fixed policy in equation \eqref{eq:tilde_operator_gap1} extends to the setting of optimal control. Consider a focal chain $\mathcal{C}$ with optimal value and policy $V^*,\pi^*$, and a sister chain $\wtilde{\mathcal{C}}$ with optimal $\tV^*,\tpi^*$, we can conclude the following about the optimality gap
\begin{equation}
    \mid V^*(x) - \tV^*(x) \mid  \leq  \frac{ \alpha}{1-\alpha} (\delta_{V^*}[P^{\pi^*}, \wtilde{P}^{\pi^*}] + \delta_{\tV^*}[P^{\tpi^*}, \wtilde{P}^{\tpi^*}]),
    \label{eq:opt_delta}
\end{equation} where $\delta$ is the one-step-ahead gap in \eqref{eq:tilde_operator_gap1}. This, in combination with our coarse grid scheme, leads to the following optimization analog of Theorem \ref{thm:backtodelta}.

\begin{theorem}  
  $$\mid V^* - \tV^* \mid_{\calS}^*  \leq \frac{ 1}{1-\alpha}\left( |V^*-GW^*|_{\calS}^*  +|\tV^*-G\tW^*|_{\calS}^*\right).$$
  \label{thm:backtodelta2}
\end{theorem}

The gap depends then on how well the optimal value under focal chain $P$ and the sister chain $\wtilde{P}$ are approximated by interpolating values at their nearest grid-points. 

With fixed policies $\pi^*$ and $\wtilde{\pi}^*$ are fixed, moment matching supports---as seen in earlier sections---a small value-approximation gap with non-negligible computational gains. In embedding this within policy iteration it is important (indeed central) that our construction of the aggregation matrices $\bg,\bu$ does not depend on the policy and, hence, does not have be updated in each iteration. 

The base algorithm is the aggregate analogue of  standard policy iteration (PI) and alternates between evaluation and updating steps. Exact evaluation and/or update are replaced by approximate computations in approximate policy iteration (API); using aggregation, the $k^{th}$ iteration proceeds as follows:

\begin{itemize}
    \item[(i)] Evaluation: for current policy $\pi^k$ and induced $P^{\pi ^k}$, compute
    $R^k=(I-\alpha \bu P^{\pi ^k}\bg)^{-1}\bu c$ 
    \\(see \eqref{eq:aggregateBellman},\eqref{eq:agginverse}). 
    \item[(ii)] Update: find policy $\pi ^{k+1}$ that satisfies 
    $$\pi^{k+1}(x)\in \argmin_{a \in \mathcal{A}(x)}\left\{c(x,a)+\alpha \sum_{y\in\mathcal{S}}p_{xy}^a\sum_{l\in\calM}\sg_{yl}R^k(l)\right\}$$
    \end{itemize} 
This is nothing but policy iteration for a chain on $\calM$ with transition matrix $\bu P\bg$ and, as such, it is guaranteed to converge; see \cite[Proposition 6.4.2]{puterman1994markov}.

\begin{algorithm}
\floatname{algorithm}{Algorithm}\caption{(\moma~API)\label{alg:onestep}}
\renewcommand{\thealgorithm}{}
\begin{algorithmic}[1]
	\Require Spacing exponent $\fraks \in [\frac{1}{3}, \frac{1}{2}).$ %Fix $\texttt{Grid}$ and the matrices $\bu=\bu(\texttt{Grid})$ and $\bg=\bg(\texttt{Grid})$.
	\Ensure Policy $\pi^*$.
	\State {\em \moma~Pre-processing}: Obtain $\bu, \bg$ using Algorithm \ref{alg:moma}. 
	\State Set initial control on representative states $\bar \pi^0: \calS^0 \rightarrow \mathrm{A}$. 
	\State Compute the induced transition $\bar P^{\bar{\pi}^0}: \calS^0 \rightarrow \calS$. Set $\bar{P}^0 \leftarrow \bar P^{\bar{\pi}^0}$.
	\While{convergence criterion is not met}
	\State {\em Policy Evaluation}: Compute $R^k=(I-\alpha  \bar P^k \bg )^{-1}\bu c$.
	\State {\em Policy Update}: \begin{align*} &\bar \pi^{k+1}(x_l)\leftarrow \argmin_{a \in \mathrm{A}(x_l)}\{c(x_l,a) + \alpha [\bar P^a \bg R^k](x_l)\},\mbox{ for } x_l\in\calS^0,\\&
	\bar P^{k+1}\leftarrow \bar P^{\bar{\pi}^{k+1}}.\end{align*} 
   	\EndWhile
   	\State {\em Full Update}: %For resulting $R$, compute For each $x \in \calS$, compute 
   	$$\wtilde \pi(x)\leftarrow \argmin_{a}\{c(x,a) + \alpha [P^a \bg R](x)\}, ~ \forall x \in \calS.$$
\end{algorithmic}
\end{algorithm}
\noindent Into this general schema we add two ingredients: \vspace*{0.1cm} 
\begin{itemize} 
\item[1.] {\bf \moma~preprocess.} As detailed in Algorithm \ref{alg:moma}: we build the $\fraks$-spaced \texttt{Grid}, and create the binary disaggregation matrix $\bu(\texttt{Grid})$ that has $\su_{lx_l}=1$ for all grid points $x_l$. Next we compute the non-negative row-stochastic matrix $\bg(\texttt{Grid})$, so that the $y^{th}$ row is a $y$-mean distribution over the representative states. This construction does not depend on the transition matrix and, in turn, neither on the control. It is computed once and requires no update during the PI iterations. \vspace*{0.1cm} 

\item[2.] {\bf Reduction to PI on representative states.} To further reduce computational burden---especially in the updating step---we leverage a useful implication of the coarse grid scheme. In (aggregate) evaluation we solve for $R=\bu c+\bu P^{\pi}\bg R$ by inversion $(I-\alpha \bu P^{\pi}\bg)^{-1}\bu c$. Because $\su_{lx_l}=1$ (and $\su_{ly}=0$ otherwise) we have $(\bu P^{\pi})_{ly}=p^{\pi(x_l)}_{x_ly}$, so that the only rows of $P^{\pi}$ used are those corresponding to the representative states $\calS^0=\{x_1,\ldots,x_L\}$. Defining $\bar{P}^{\pi}$ to be the $L\times N$ matrix with $\bar{P}^{\pi}_{x,\cdot}=P^{\pi}_{x,\cdot}$ for $x\in\calS^0$, we re-write $R=(I-\alpha \bar{P}^{\pi} \bg)^{-1}\bu c.$ 

Policy update can be similarly limited to $\calS^0$, as they are the only states for which we wish to compute the induced $\bar{P}$. Define $\bar{\pi}: \calS^0 \rightarrow \mathrm{A}$. We have at iteration $k$
\begin{align*} \bar \pi^{k+1}(x_l)\leftarrow \argmin_{a \in \mathrm{A}(x_l)}\{c(x_l,a) + \alpha [P^a \bg R^k](x_l)\}=
\argmin_{a \in \mathrm{A}(x_l)}\{c(x_l,a) + \alpha [\bar{P}^a\bg R^k](x_l) \},\end{align*} 
where the equality follows from the fact that the $1 \times N$ vector of probabilities $p_{x_l, \cdot}$ can be accessed from $\bar{P}$ instead of $P$.
% for each $x_k\in\calS^0$, $\sg_{x_kk}=1$ and $\sg_{x_kl}=0$ for $l \neq k$, so that $$ [P^a\bg R](x_k)= \sum_{z,l}p^a_{x_kz}\sg_{zl}R(l)= \sum_{l}p^a_{x_kx_l}\sg_{x_kx_l}R(l)= [\bar{P}^a\bg R](x_k).$$ 

Thus, we can first run a policy iteration on $\calS^0$, and only do a single full Bellman lookahead for the states $x\in \calS\backslash \calS^0$ after convergence; this is step 7 of the algorithm.  
\end{itemize} 

Convergence of  Algorithm \ref{alg:onestep} follows immediately from that of standard PI, applied here to the controlled chain $\bar{P}^{\bar{\pi}}$ on $\calS^0$. The value and policy to which this PI converges inherit an optimality-gap-bound from the approximation-gap-bound in Theorem \ref{thm:guaranteemain}.

\subsection{Complexity of optimization}

We expand the discussion of evaluation complexity in \S \ref{sec:complexity} to the entire policy iteration algorithm. 

A difficulty in the complexity calculations is that, while the approximation algorithm is more efficient {\em per iteration}, it might require more iterations to converge compared to the exact one, thus erasing any possible gains. Fortunately, the upper bound on the number of iterations is much smaller for  aggregation PI compared to the exact PI because, recall, we perform updates only for states $x\in \calS^0$; see  
\cite{ye2011simplex,singhmansour,hollanders2016improved,scherrer2013improved}. 

In \S \ref{sec:complexity} we showed that the time used for moment matching optimization is made up for by the time saved from evaluating policies for $L=|\calM|= \calO ( N^{1 - \fraks})$ instead of $N$ states.
The time savings are {\em not} however limited to the evaluation step. Computation is reduced also because we perform the policy iteration steps, up to convergence, only on the representative states $x\in\calS^0$.

Specifically, suppose the cost of policy update for a single state is $m$; in the worst case $m$ might correspond to comparing all feasible actions $a \in \mathrm{A}(x)$. In full PI, this is done for every state $x \in \calS$ so the complexity is $\calO(Nm)$. In Algorithm \ref{alg:onestep}, on the other hand, we  update only $x \in \calS^0$, giving $\calO(Lm) = \calO(N^{1-\fraks}m)$. Moreover, while implicit in the algorithm, obtaining the control-induced transition matrix $P^{\pi}$ at each iteration has non-negligible computation expense of $\calO(N^2)$ for full PI, and reduced to $\calO(LN) = \calO(N^{2-\fraks})$ each in Algorithm \ref{alg:onestep}.

Except for cases where the action space far exceeds the state space in magnitude, specifically $m > \calO(NlogN)$, the time complexity of matrix inversion in the evaluation step dominates, thus the gain in each iteration is still $\calO(N^2logN)$. These gains are multiplied by the number of iterations it takes for the aggregate values to converge; though one must also account for the time to perform one full policy update after convergence. With $T$ iterations we have
\[ Gain =   \Omega(T N^2\log N - Nd2^d - Nm) = \Omega(N^2\log N) . \]

\section{Numerical experiments\label{sec:numerical}}

We consider two operations-management problems that pose a computational challenge for exact methods. In both cases there is a natural way to scale up the complexity, starting from small instances where we can visualize the outcomes and proceeding to larger instances that test the computational benefits of \moma. Importantly, both were studied using alternative approximation methods, providing a benchmark for our approach.

\noindent {\em Notes:} (i) Where {\em an optimal policy $\pi$} is mentioned, the specific instance of the problem is exactly solvable and we compute the value and policy via standard policy iteration; (ii) All experiments reported in this section were run on a  machine with Intel(R) Core(TM) i7-6700 CPU @ 3.40GHz 3.41 GHz and 16.0GB of RAM, using 64-bit Python.

\subsection{\moma~pre-process}

\begin{figure} %[b!] 
\centering
\includegraphics[width=0.45\textwidth]{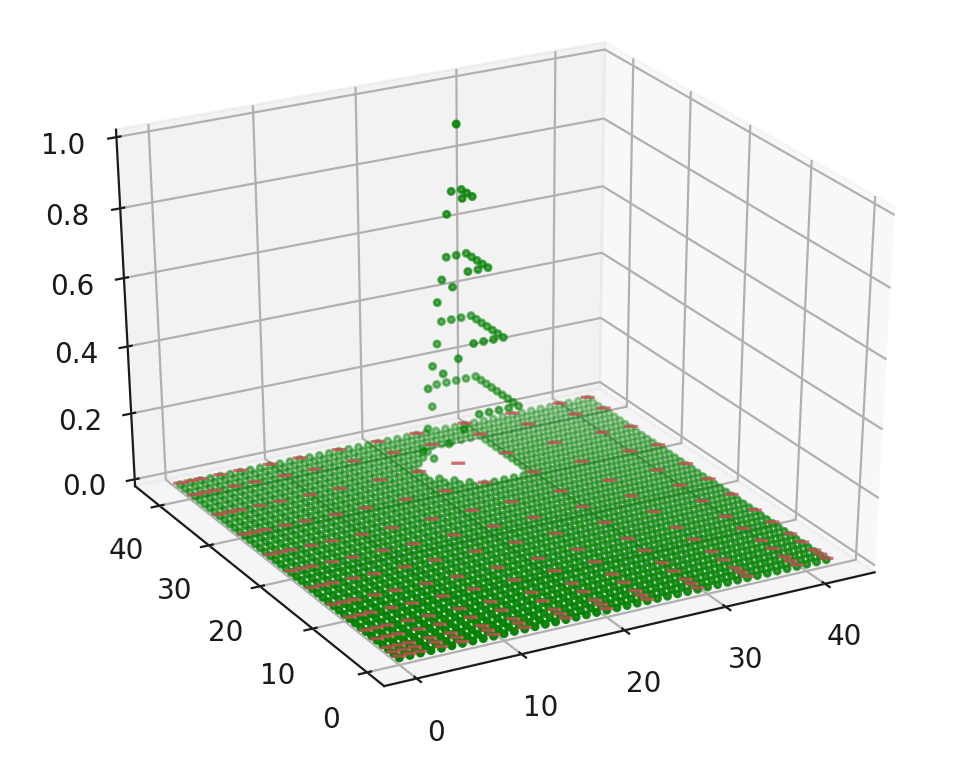}
\includegraphics[width=0.5\textwidth]{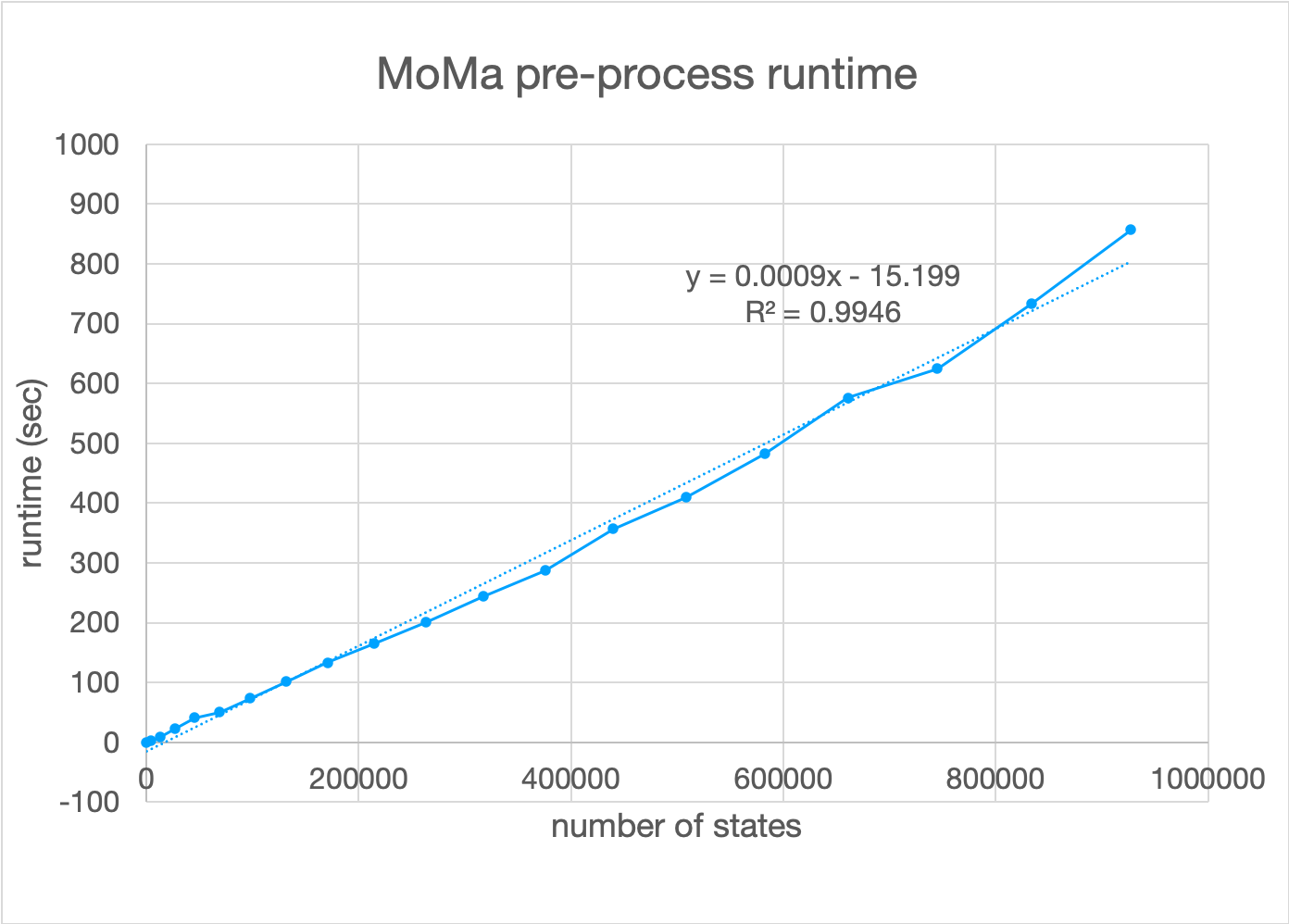}
\caption{(LEFT) Highlighted with red ticks the representative states that form the grid. In green we plot, for one fixed meta-state, the aggregating probabilities $\bg$ of each state into a fixed meta-state $l$; only states that are no further than the nearest neighboring meta-states aggregate into $l$ with positive probabilities; $\sg_{yl}$ is proportional to the distance from state $y$ to representative state $x_l$ (RIGHT) We scale up the (two-dimensional) state space $[0,u]^2$ by raising the value of $u$; pre-processing takes less than 25 minutes even for state space with size larger than a million.
  \label{fig:moma}}
\end{figure}

Common to both examples is the \moma~pre-processing step as detailed in  Algorithm \ref{alg:moma}. Figure \ref{fig:moma} visualizes the grid, representative states and the aggregating probabilities $\bg$ for the case of $d=2$ and state space $[0, 40]^2 \cap \mathbb{Z}^2$. The right-hand side in said figure confirms the linear scaling--in the number of states $N$---of the pre-processing. 

\subsection{Joint replenishment problem}

A retailer carries two types of products. Demand for each product is independent (across products and time periods). There are two types of {\em fixed} ordering costs: (i)  a \textit{minor} ordering cost for placing an order for product $i$; and (ii) orders of both products can arrive in the same truck and a \textit{major} ordering cost is incurred for each {\em truckload}. The number of truckloads then depends on the total amount ordered (of both products). 

We follow the standard setup as very clearly laid out in \cite{vanvuchelen2020use}. For simplicity, only full truckloads are considered. 

At time $t$, the order amount $q_{i,t}$ for each item type $i=1,2$ is determined based on the inventory level $I_{i,t}$ at the end of the previous period. Lead time is assume to be $0$ and orders arrive before the demand $d_{i,t}$ is realized. The system dynamics are given by $$I_{i,t} = I_{i, t-1}+q_{i,t} - d_{i,t}$$

Per-item holding cost $H_i$ is incurred for product-$i$ inventory per unit of time. Per-item backorder cost $B_i$ is incurred for unmet demand. The minor ordering cost for product $i$ is $k_i$, and $K$ is the cost per truckload. The immediate cost function at period $t$ is  then 
\[ c(I_t, q_t) = \sum_i ( H_i [I_{i, t}]^+ + B_i [I_{i,t}]^- + k_i \mathbbm{1}_{q_{i,t}>0}) + K \left\lceil \frac{\sum_i q_{i,t}}{TC} \right\rceil \]
where $TC$ is the truck capacity. 

\subsubsection{Small instance} 

Demand is $d_1 = \mathcal{U}\{0,5\}, d_2 = \mathcal{U}\{0,3\}$, where we use $\mathcal{U}\{a,b\}$ to denote the discrete uniform distribution between $a$ and $b$. Each truck can carry 6 items. The detailed parameters are listed in Table \ref{table:params_JRP_small}. The only difference between the two products is the minor ordering cost. The discount factor is $\alpha = 0.99$. 

\begin{table}[h!]
\centering
\begin{tabular}{| c |c| c| c| c| c|c| c| }
\hline
 \text{\textbf{item type}}  & d & H & B & k & K & $\ell_i$ & $u_i$  \\
 \hline
i=1 &	U\{0,5\} & 1 & 19 & 40 & 75 & -30 & 40  \\  
 \hline
i=2  &	U\{0,3\} & 1 & 19 &	10	& 75  & -30 & 40 \\
\hline
\end{tabular}
\caption{Demand parameters for the small instance of the joint replenishment problem.}
\label{table:params_JRP_small}
\end{table}

We truncate the inventory for each item at 40 and the backorder is capped at 30 units for each item type; this means the order quantity must satisfy $q_{i,t}\leq 40-I_{i,t}$. The total number of states is $N=5041$. We take the \moma~spacing exponent to be $\fraks = 0.45$, resulting in $L=400$ meta states. 

\begin{figure}
\centering
\includegraphics[width=0.43\textwidth]{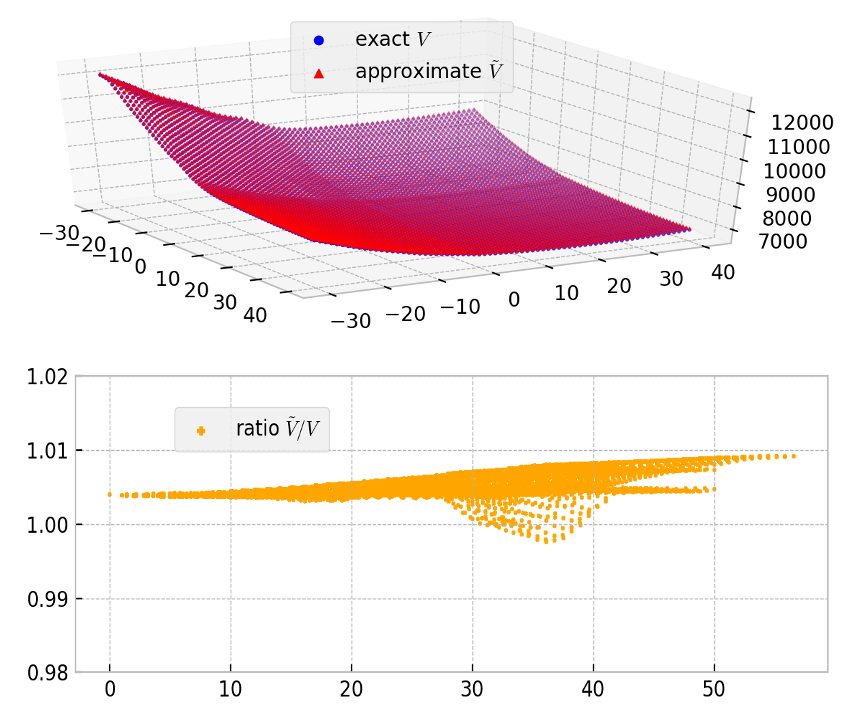}
\includegraphics[width=0.5\textwidth]{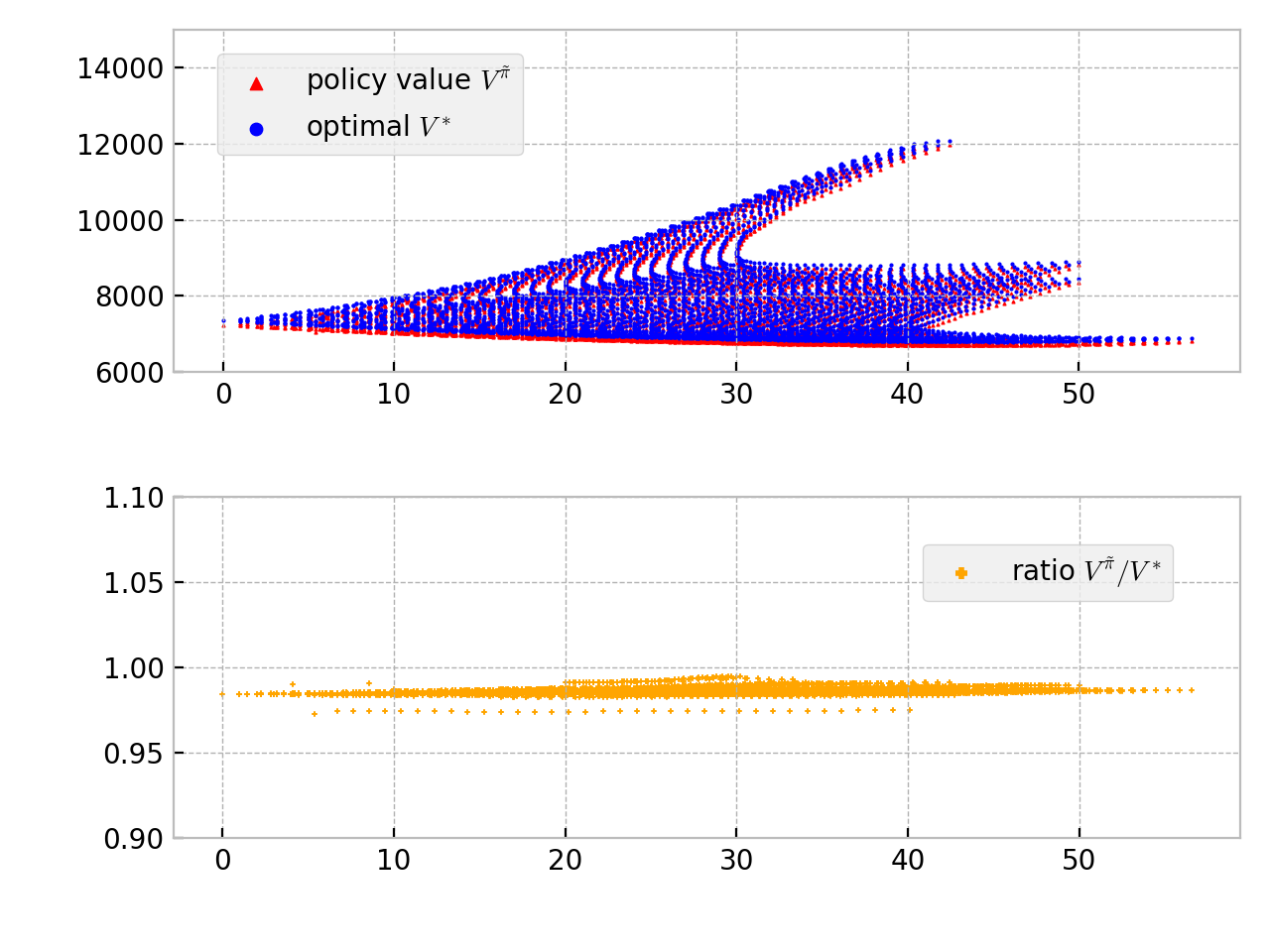}
\caption{(TOP LEFT) Comparison of the approximately evaluated value function vs the exact against the state space. (BOTTOM LEFT) Their ratio against the Euclidean distance to the origin. (RIGHT) Comparison of the performance of the approximate policy $V^{\tpi}$ vs the optimal value $V^*$ in (TOP), their ratio in (BOTTOM), both against the Euclidean distance to the origin.  \label{fig:JRP-small}}
\end{figure}

First, we test the evaluation performance of \moma; see Algorithm \ref{alg:eval}. The instance is small enough that we can compute the exactly optimal policy $\pi^*$. We take $P = P^{\pi^*}$ as the transition function for the focal chain and compute the approximate value $\tV(x) = c(x) +\alpha P\bg R(x)$, where $R=(I-\alpha \bu P\bg)^{-1}\bu c$ is the aggregate value. This is displayed against the exact value $V = (I-\alpha  P)^{-1}c$ in Figure \ref{fig:JRP-small}; the mean and max (over the state space) of the evaluation gap as percentages of the exact value are 0.51 \% and 0.92 \% respectively.

We consider optimization next. We obtain the candidate policy $\tpi$ using Algorithm \ref{alg:onestep}, and compare the value $V^{\tpi}$ of this policy against the optimal value $V^*$; see Figure \ref{fig:JRP-small} (RIGHT). The relative optimality gap  $\mid V^{\tpi} - V^* \mid / V^*$, has a mean of 1.38 \% and max of 2.73 \%. Based on simulation, \cite{vanvuchelen2020use} reports an optimality gap with mean 0.46 \% max 0.91 \% for their neural network method. Computation times are not reported in \cite{vanvuchelen2020use}.

\begin{figure}
    \centering
    \includegraphics[width=0.95\textwidth]{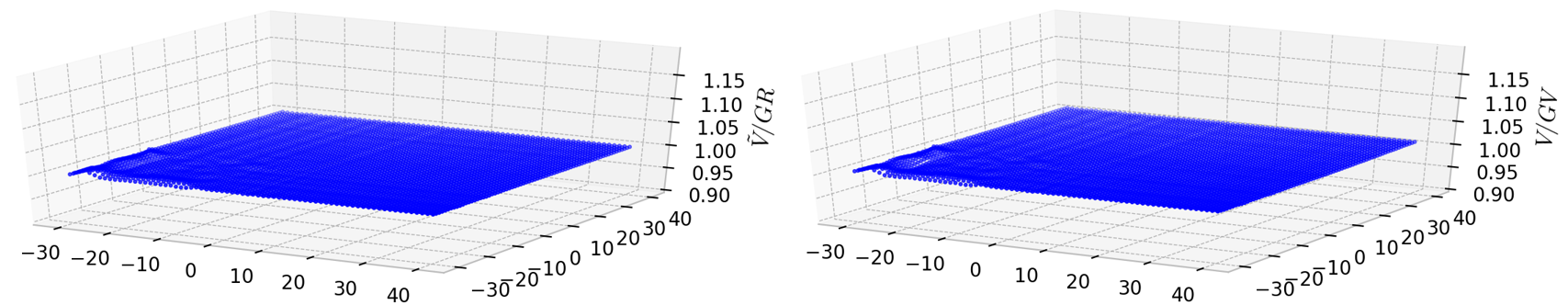}
    \caption{Ratio of $V, \tV$ against their interpolated values for the small instance of joint replenishment problem.}
    \label{fig:remark_JRP}
\end{figure}

Theorem \ref{thm:backtodelta} relates our approximation's quality  to the ``local linearity'' of the values $V$ and $\tV$, i.e, to how well the value at a state is a distance-proportional interpolation of the values at the gridpoints. Figure \ref{fig:remark_JRP} shows that such local linearity holds for the joint replenishment problem and explains why we are observing such impressive accuracy. 

\iffalse 
When we examine the ratio between the value functions $V, \tV$ and their interpolation $\bg \bar{V}, \bg R$, they are both close to 1 with little fluctuation, consistent with the small approximation gaps observed. It seems that an immediately available explanation for this good performance is in the second moments, which form roughly piece-wise linear manifolds, thus taking the first order Taylor approximation results in small remainders. Another factor here is the existence of patches of structure in a desirable policy. A common proxy for performance of approximate policy $\wtilde \pi$ is Bellman residual, $\wtilde V$, $\lVert T^{\wtilde \pi} \tV - \tV \rVert$, as examined in e.g. \cite{williams1993tight, antos2008learning,farahmand2010error}. We can see clear patch-patterns when Bellman error as a percent of the optimal value function is plotted against the state space. The policies themselves display certain consistency.
\fi 

\begin{figure}
    \centering
    \includegraphics[width=0.47\textwidth]{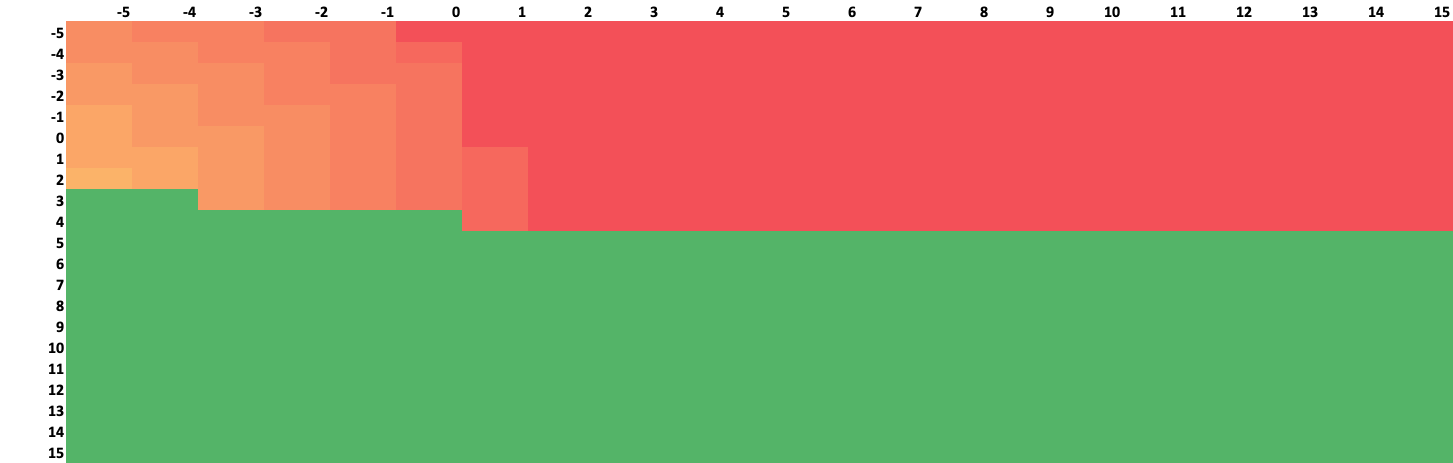}
    \includegraphics[width=0.47\textwidth]{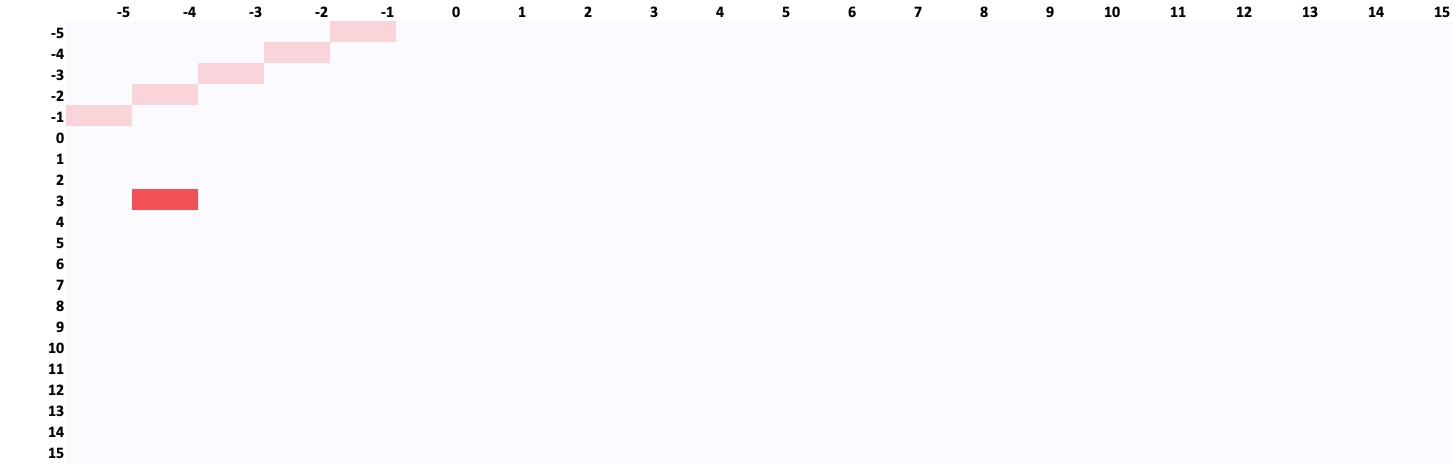}
    \includegraphics[width=0.47\textwidth]{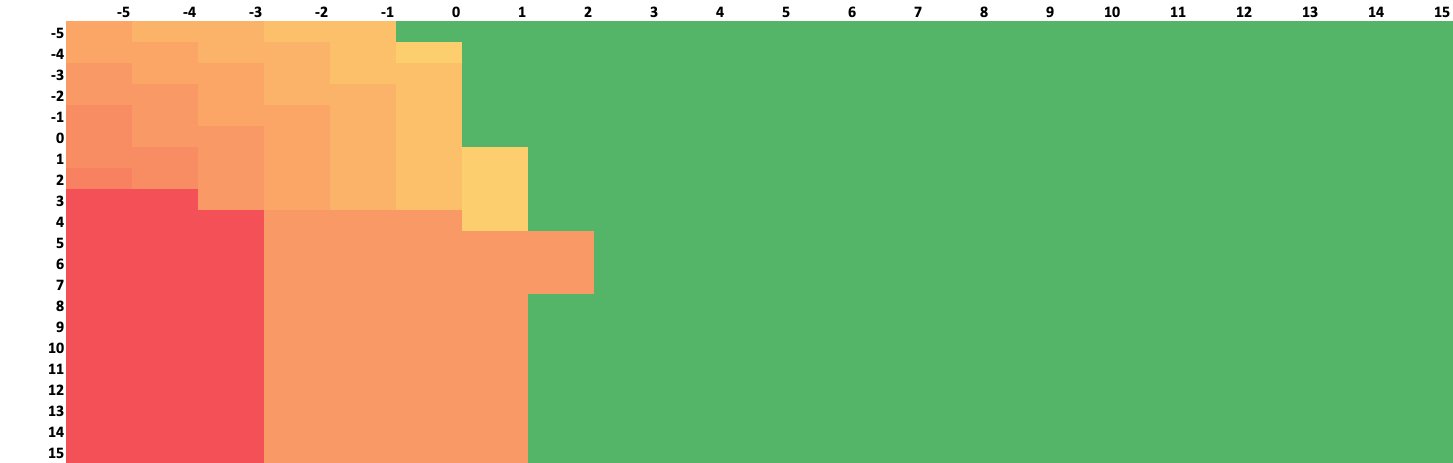}
    \includegraphics[width=0.47\textwidth]{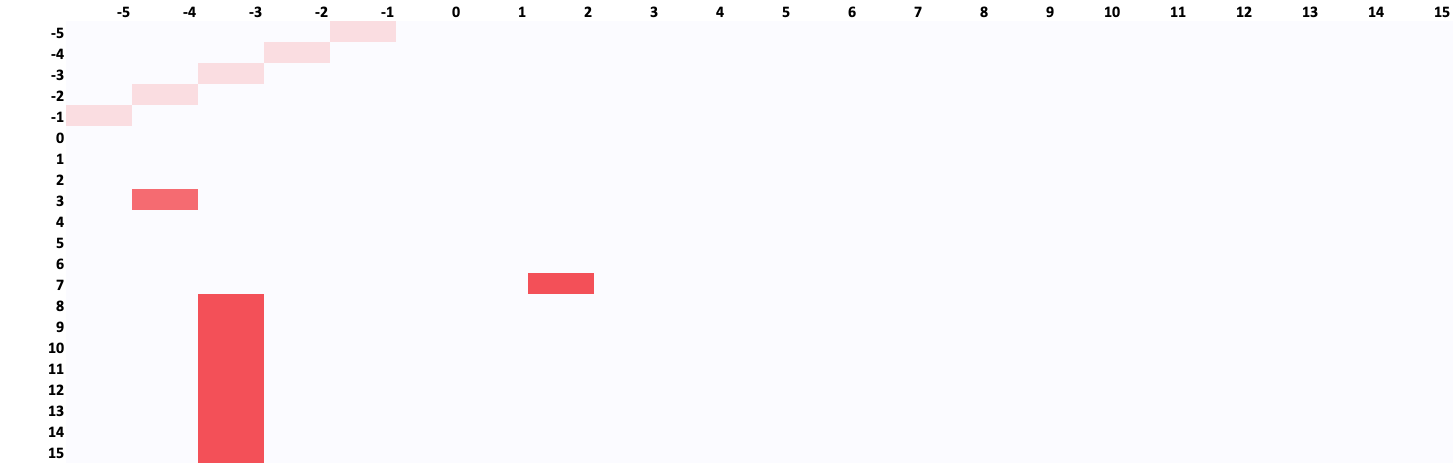}
    \caption{(LEFT) Optimal order policies for item type 1 in (TOP) and item type 2 in (BOTTOM), in color gradient. (RIGHT) Difference from approximate policies for item type 1 in (TOP) and item type 2 in (BOTTOM); differences are marked in shades of red depending on magnitude.}
    \label{fig:pol_JRP}
\end{figure}

\subsubsection{Large instance.} We consider a larger instance studied in  \cite{vanvuchelen2020use}. The parameters are as reported in Table \ref{table:params_JRP_large}. In addition each truck can carry 33 items and the discount factor is $\alpha = 0.99$.

We cap the capacity at 120. \footnote{That is we restrict order quantity to  satisfy $q_{i,t}\leq 120-I_{i,t} + \text{minimum demand}$.} For each item type, and backorder at 50 units for each item type. The state space is then $\calS = [-50,120]^2$ and the total number of states is $N=|\calS|=29241$; using $\fraks=0.45$ we have $L = 1089$ meta states.

\begin{table}[h!]
\centering
\begin{tabular}{| c |c| c| c| c| c|c| c| }
\hline
 \textbf{item type}  & d & H & B & k & K & $\ell_i$ & $u_i$  \\
 \hline
i=1 &	$\mathcal{U}$\{15,25\} & 7 & 19 & 40 & 400 & -50 & 120  \\  
 \hline
i=2  &	$\mathcal{U}$\{5,15\} & 1 & 19 &	10	& 400 & -50 & 120  \\
\hline
\end{tabular}
\caption{Demand parameters for the larger instance of the joint replenishment problem.}
\label{table:params_JRP_large}
\end{table}

\begin{figure}
\centering
\includegraphics[width=0.43\textwidth]{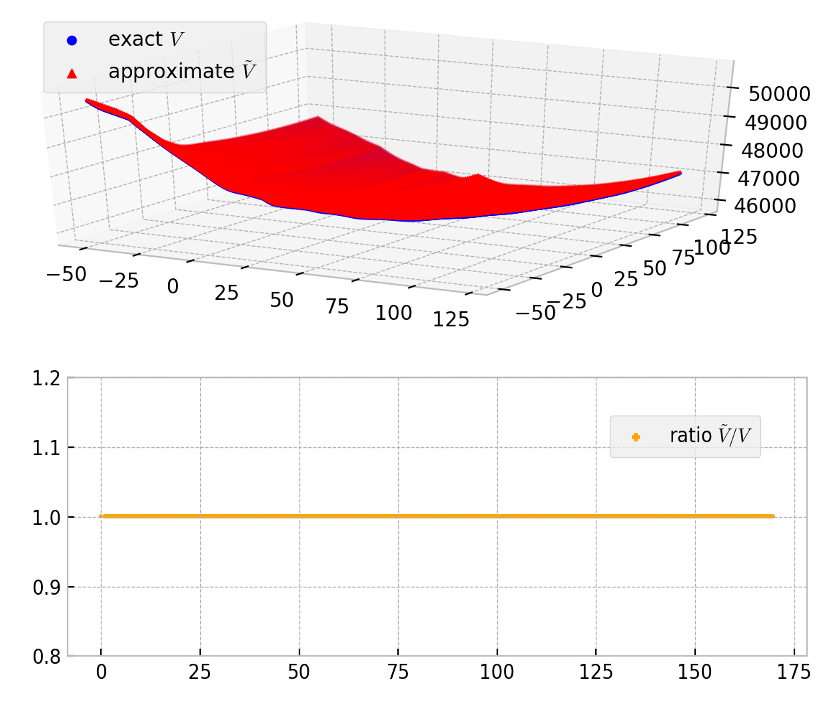}
\includegraphics[width=0.5\textwidth]{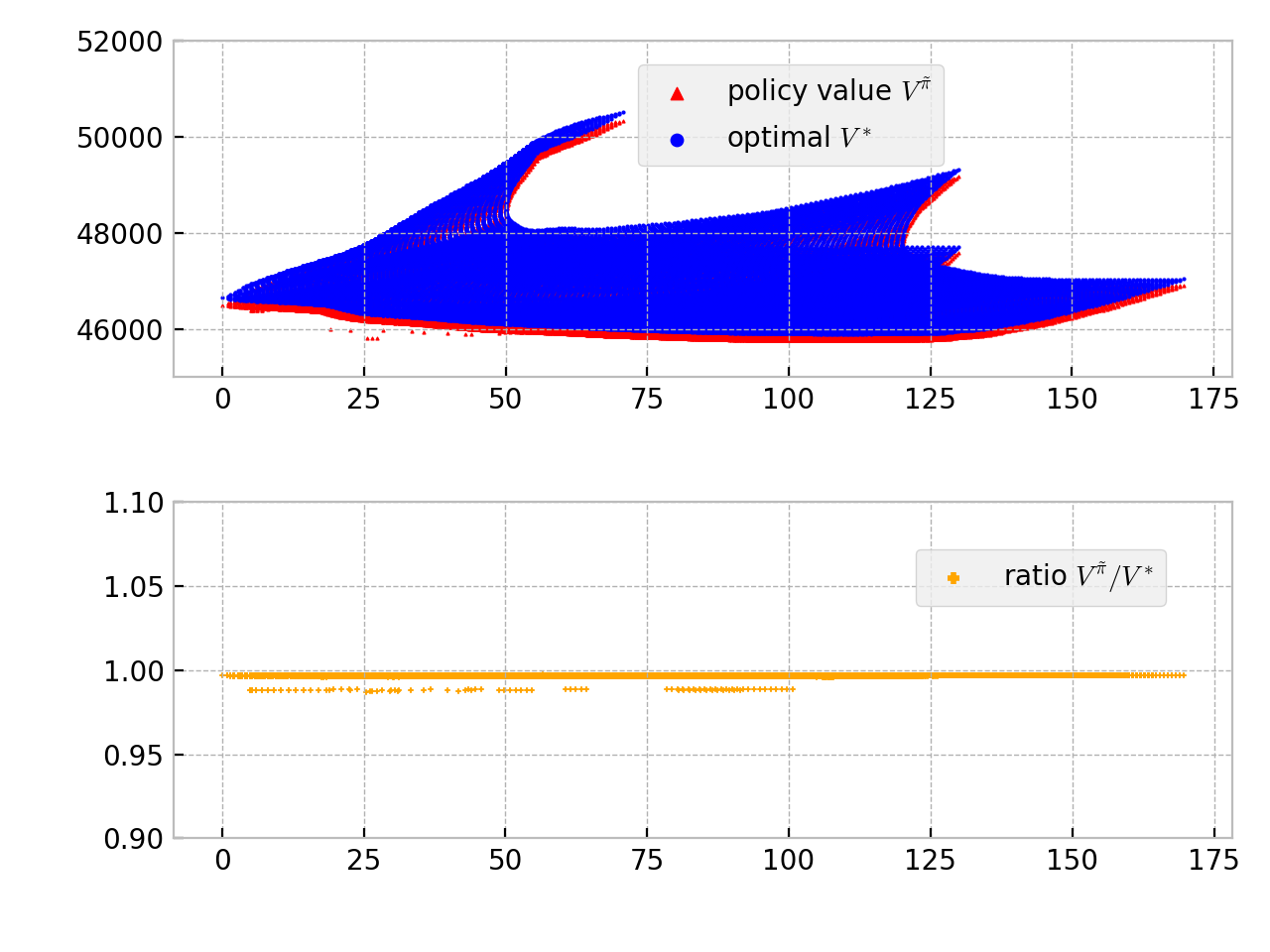}
\caption{(TOP LEFT) Comparison of the approximately evaluated value function vs the exact against the state space. (BOTTOM LEFT) Their ratio against the Euclidean distance to the origin. (RIGHT) Comparison of the performance of the approximate policy $V^{\tpi}$ vs the optimal value $V^*$ in (TOP), their ratio in (BOTTOM), both against the Euclidean distance to the origin.  \label{fig:JRP-large}}
\end{figure}

\begin{figure}
    \centering
    \includegraphics[width=0.4\textwidth]{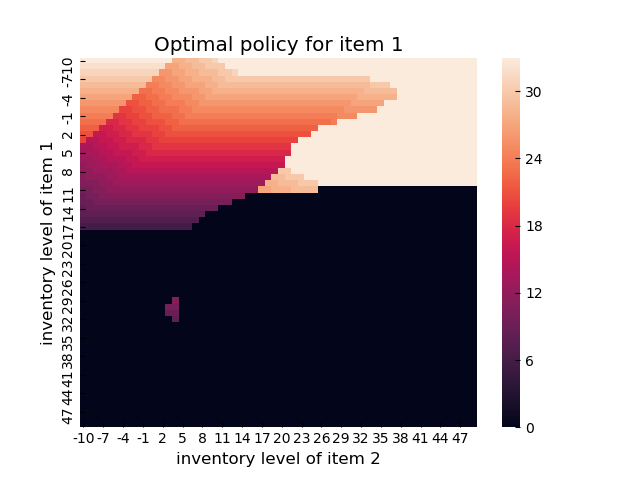}
\includegraphics[width=0.4\textwidth]{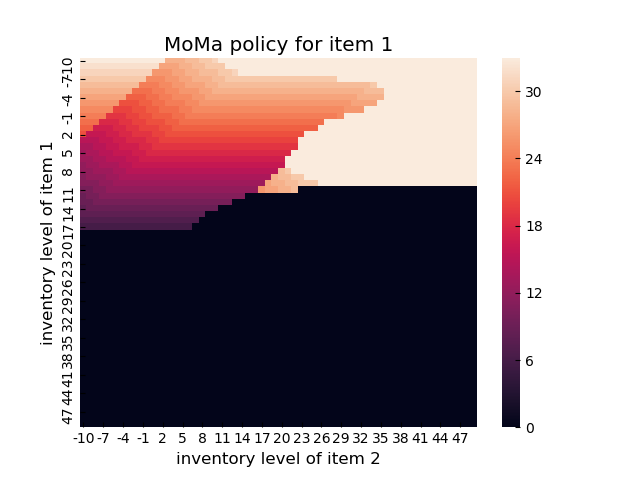}
\caption{Heat map of the prescribed ordering amount for item 1 against inventory levels, according to policies produced by exact PI (left) and \moma Aggregate PI (right).}
    \label{fig:pol_compare}
\end{figure}

The performance of \moma~is captured in Figure \ref{fig:JRP-large}. The evaluation gap as a percent of the exact value function has mean 0.11 \% and max 0.13 \%. The optimality gap as a percent of the optimal value function has mean 0.32 \% and max 1.29 \%. While the optimality gaps from the method in \cite{vanvuchelen2020use} are not reported for this larger instance, precluding direct comparison, the near optimality of the \moma~policy is illustrated in Figure \ref{fig:pol_compare}, where order quantities of item type 1 prescribed by the policies are displayed.

In terms of computation time, \moma~API took less than 6 minutes to converge, whereas exact PI took more than 2 hours; see detailed breakdown of runtime in seconds in Table \ref{table:JRP_time}. Note that the time used for each step is averaged across iterations (quoted with the unit of seconds per iteration), whereas the time cost of \moma~preprocess is incurred only once (quoted with the unit of seconds). Exact runtime is not reported in \cite{vanvuchelen2020use}.

The last step of the algorithm---to obtain the optimal actions for all states---always requires one full update. That cost is unavoidable unless one interpolates the control obtained for the representative states. 

\begin{table}[h!]
\centering
\begin{tabular}{| c |c| c| c| c| c|c| }
\hline
 \textbf{algorithm}  & \textbf{update} & \textbf{compute $P$} & \textbf{evaluation}  & \textbf{\# iter} & \textbf{\moma~preprocess} & \textbf{total}  \\
 \hline
Exact PI &	696.90 /iter & 12.76 /iter & 119.79 /iter & 9 & / & 7597.68  \\  
 \hline
\moma~API  &	35.64 /iter & 0.51 /iter & 0.10 /iter &	9	& 30.28  & 357.15  \\
% Exact PI &	524.79 sec/iter & 12.19 sec/iter & 75.76 sec/iter & 9 & / & 5602.70 sec \\  
%  \hline
% \moma~API  &	21.05 sec/iter & 0.42 sec/iter & 0.07 sec/iter &	9	& 23.78 sec & 218.10 sec \\
\hline
\end{tabular}
\caption{Runtime breakdown for the larger instance of the joint replenishment problem.}
\label{table:JRP_time}
\end{table}

\subsection{Inpatient-flow optimization}

This second example follows on \cite{dai2017two} which is already re-considered in \cite{braverman2018taylor}. It is a hospital routing problem with multiple patient types and dedicated hospital wards; patients waiting in queue for beds in their specialized wards could be routed (or ``overflown') to a different ward in the hospital at a cost. The $J$ internal wards are the server pools in this discrete-time queuing model, and the $N_j$ beds in each constitutes the servers. Arrivals of type-$j$ in a period $t$ follow a Poisson random variable with mean $\lambda_j$, and arrivals are independent across types and time periods. Once admitted to ward $j$, a patient's length of stay in pool $j$ (the time occupying a bed) is geometrically distributed with mean $1/p_j$. 

An arriving type-$j$ patient is immediately assigned a bed in pool $j$ when available, and otherwise waits for service in an (infinite) type-$j$ queue; the queue is truncated for numerical experiments. We let $X_j (t)$ be the number of patients either in service in pool $j$ or waiting in queue $j$; we use $X(t)$ for the vector process.

While waiting, a patient in queue $j$ incurs a holding cost $H_j$ per period of delay. A waiting type-$j$ patient can be re-routed to (an unsaturated) pool $j \neq i$ at a cost of $B_{ij}$ and served immediately. A re-routing from $i$ to $j$ can happen only if there are available server in pool $j$; we do not re-route a waiting customer to another queue. This \emph{overflow decision} is made at the start of the time period, before departures and arrivals are realized. Let $U_{ij}(t) = U_{ij}(X(t)) $ be the number of customers overflown from buffer $i$ to pool $j$ at time period $t$. The action space in state $x$ is then 
$$ \bigg\{ U \in \mathbb{Z}^{J  \times J}  \mid \sum_{i \neq j} U_{ij} \leq (N_j - x_j) ^+, \sum_{j \neq i } U_{ij} \leq (x_i - N_i) ^+ \bigg\}. $$ The number of type-$i$ patients routed to other pools cannot exceed those waiting in buffer $i$, $(x_i-N_i)^+$, and the number routed to pool $j$ cannot exceed the number of available servers there, $(N_j-x_j)^+$.

The discrete-time dynamics are given by 
$$ X_i(t) = X^P_i (t-1) + A_i (t-1) - D_i(X^P_i (t-1)),$$
where $A_i(t) \sim $ Poisson($\lambda_i$) is the number of type-$j$ arrivals, $D_i (x) \sim Binomial(x \wedge N_i, p_i)$ is the number of type-$i$ departures, and where $$ X^P_i (t-1) = X_i (t-1) + \sum_{j \neq i} U_{ji} (X(t-1) - \sum_{j\neq i} U_{ji}(X(t-1))$$ is the post-action state. The cost is incurred immediately after the re-routing but before arrivals and service completions (i.e., before the realization of randomness). It is given by 
$$\sum_{i} \sum_{j \neq i} B_{ij} U_{ij} + \sum_i H_i \times (x_i - \sum_{j \neq i} U_{ij} - N_i) ^+.$$

\subsubsection{Small instance} 
% \textcolor{red}{re-order: first the real parameters then the artificial ones -- like truncation and \moma~parameters}
Here we consider an instance with 2 wards so again $\calS \subseteq \mathbb{Z}_+^2$. Parameters are listed in Table \ref{table:params_hospital_small}. They are chosen so that $\lambda_1/p_1 = 14>12$ for the ward 1 and $\lambda_2/p_2= 8<12$ for ward 2, resulting in pressure to overflow patients from the overloaded ward 1 to ward 2; $B_{12}=5> B_{21} = 1$ so that such overflow is costly. 

Each of the queues is truncated at $30$ so the state space is $\calS=[0,42]^2\cap \mathbb{Z}_+^2$ and $N=|\calS|=1849$. We take the \moma~spacing exponent to be $\fraks = 0.45$ resulting in  $L=196$ meta states. 

\begin{table}[h!]
\centering
\begin{tabular}{| c |c| c| c| c| c| c| c| c| c|}
\hline
  & $\lambda_i$ & $p_i$ & $H_i$  & $B_{i1}$ & $B_{i2}$ & $N_i$ & $\ell_i$ & $u_i$\\
 \hline
i= 1 &	3.5 & 0.25 & 5 &  / & 5 & 12 & 0 & 42  \\  
 \hline
i= 2  &	2.8 & 0.35 & 5 & 1 & / & 12 & 0 & 42 \\
\hline
\end{tabular}
% \quad \quad \quad
% \begin{tabular}{| c |c| c| }
% \hline
% B_{ij}  & ward 1 & ward 2  \\
%  \hline
% ward 1 &  / & 5 \\  
%  \hline
% ward 2  & 1 & /	   \\
% \hline
% \end{tabular}
\caption{Parameter setting for the small instance of the hospital routing problem.
% with  separate table whose entry in row $i$ column $j$ shows the overflow cost $B_{ij}$.
}
\label{table:params_hospital_small}
\end{table}

\begin{figure} [t!] 
\centering
\includegraphics[width= 0.55 \textwidth]{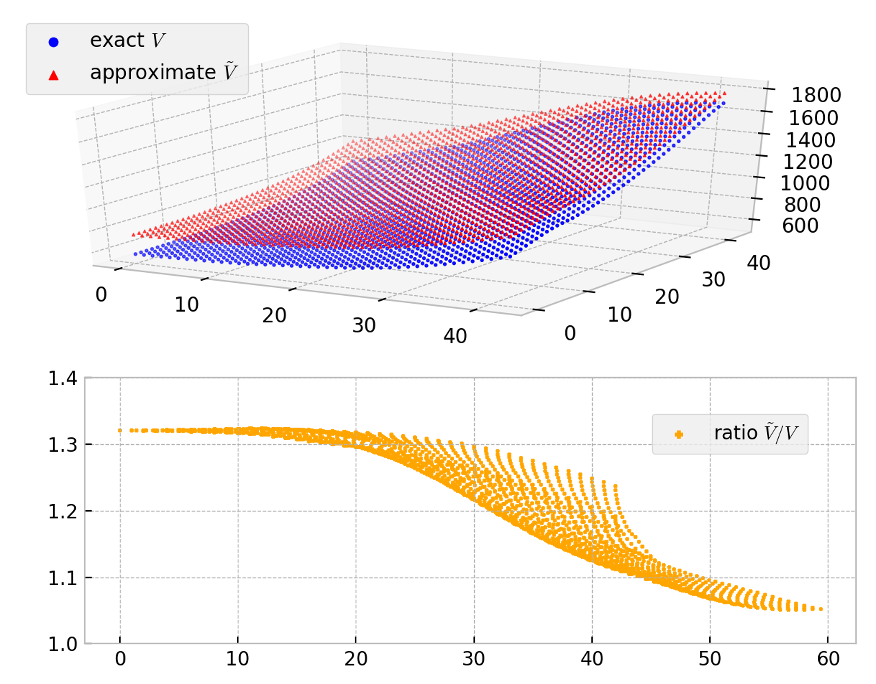}
\includegraphics[width= 0.42 \textwidth]{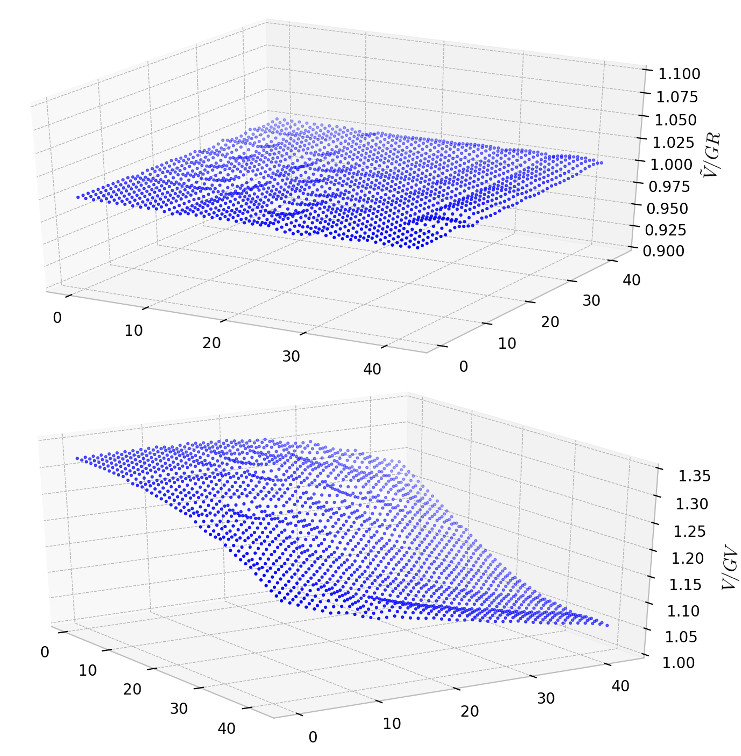}
\caption{(LEFT TOP) Exact vs approximate values: against the state space. (LEFT BOTTOM) Their ratio against the Euclidean distance to origin. (RIGHT) Ratio of $V, \tV$ against the interpolated value. 
 \label{fig:Vs-eval}}
\end{figure}

We first use \moma~for evaluation using the exact optimal policy and compare the value of the focal chain $V=V^*$ to $\tV$. Figure \ref{fig:Vs-eval} (LEFT) shows the gap. The approximation has a visible gap from the true value, and the largest deviation happens near the origin; this is consistent with our accuracy guarantees where the gap-bound is smaller the farther one is from the origin. Figure \ref{fig:Vs-eval} explains this in the terms of Theorem  \ref{thm:backtodelta}. We see that $V$ does not have the clear ``local linearity'' we observed in the replenishment problem; at least not near the origin.

\begin{figure}
\centering
\includegraphics[width=0.45\textwidth]{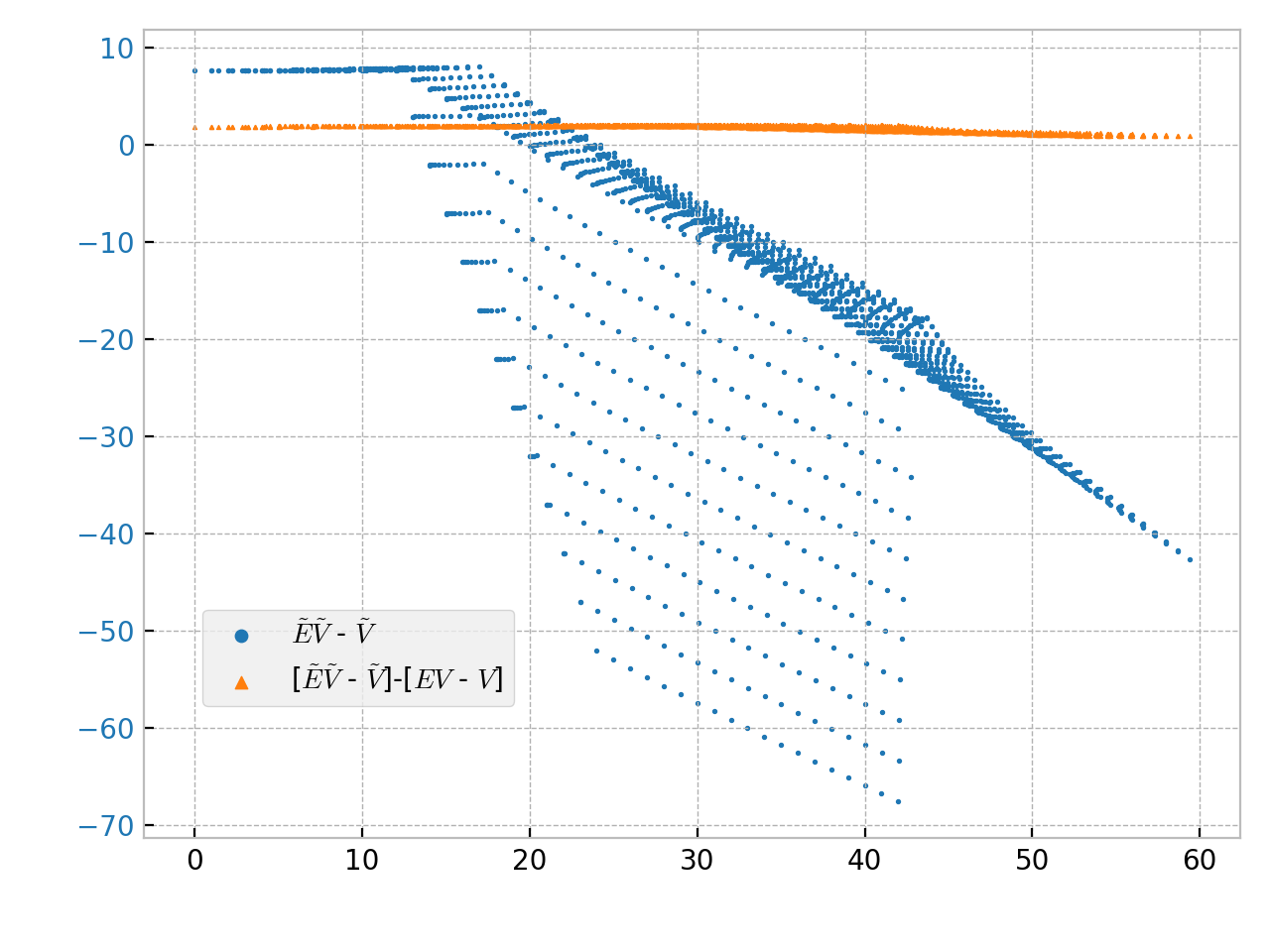}
\includegraphics[width=0.45\textwidth]{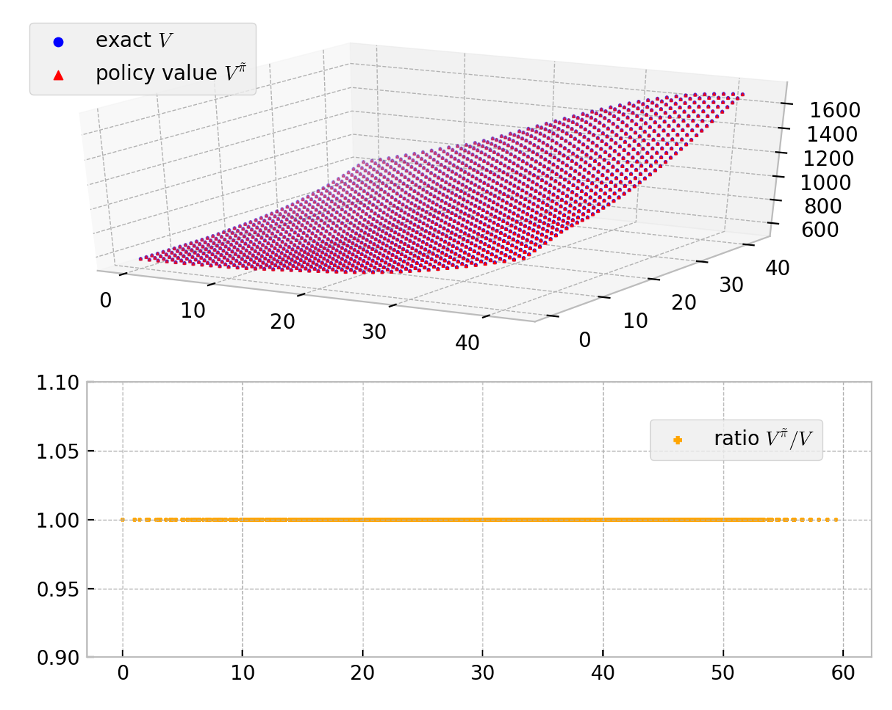}
\caption{(LEFT) The incremental changes against the Euclidean distance to the origin, (RIGHT) Their difference against the state space.  \label{fig:inc-c}}
\end{figure}

It is important to note however that the ``shapes'' of $V$ and $\tV$ are very similar. This matters for optimization: an optimal control $\pi^*$ satisfies  $$\pi^*(x) \in \argmax_{a\in\mathrm{A}(x)}\{c(x,a)+\alpha(\Ex_x^a[V(X_1)]-V(x))\},$$ so the main influence of value $V$ on the prescribed action is through the increment $\Ex_x^a[V(X_1)]-V(x)$. The control computed using $\tV$ will be influenced by the approximate increment
\begin{align*}
   \tEx_x[\tV(X_1)]-\tV(x)&= P\bg R(x)-\tV(x). %\tag{one-step}  
\end{align*}
Figure \ref{fig:inc-c}(LEFT) plots these increments and captures how close they are.

It is then less surprising that---where it matters most, i.e., in the context of optimization --- \moma~performs exceedingly well here. This is confirmed in Figure \ref{fig:inc-c}(RIGHT), where $\tpi$ is computed using \moma~API and its performance $V^{\tpi}$ compared to the optimal $V$.  

\subsubsection{Three and four wards}
We replicate an instance studied in \cite{braverman2018taylor} with 3 specialty wards with parameters as listed in \ref{table:params_hospital_3ward}. This instance has a load levels $0.7$, i.e., $\lambda_i = 0.7  N_i p_i$. We later consider also a higher $0.8$ load. The discount factor is set to $\alpha = 0.99$. 
 The total number of states is $N=15625$, and with $\fraks=0.45$, there are $L=1000$ meta states. 
\begin{table}[h!]
\centering
\begin{tabular}{| c |c| c| c| c| c| c| c| c|  }
\hline
  & $\lambda_i$ & $p_i$ & $H_i$  & $B_{i1}$ & $B_{i2}$ & $B_{i3}$ & $\ell_i $& $u_i$ \\
 \hline
i= 1 &	2.8	 & 0.4 & 10 & / & 5 & 2 & 0 & 24 \\  
 \hline
i= 2  &	4.2	 & 0.6 & 2 & 3 & / & 7  & 0 & 24 \\
\hline
i= 3  &	0.7 & 0.1 & 6 & 7 & 9 & /  & 0 & 24 \\
\hline
\end{tabular}
\caption{Parameter setting for the 3-ward instance of the hospital routing problem with load level 0.7.
% with  separate table whose entry in row $i$ column $j$ shows the overflow cost $B_{ij}$.
}
\label{table:params_hospital_3ward}
\end{table}

For this instance, exact PI converged in 8.4 hours. Meanwhile \moma~API converged in 46 minutes---a time saving of over 90 \%; see  
Table \ref{table:3ward_time} for a detailed runtime breakdown. Runtime of ``below 10 minutes" is reported in \cite{braverman2018taylor}, but without specifying whether it was for the setting of grid size 1728, 216 or 27.

\begin{table}[h!]
\centering
\begin{tabular}{| c |c| c| c| c| c|c| }
\hline
 \textbf{algorithm}  & \textbf{update} & \textbf{compute $P$} & \textbf{evaluation}  & \textbf{\# iter} & \textbf{\moma~preprocess} & \textbf{total}  \\
 \hline
Exact PI &	5530.81 /iter & 423.30 /iter & 22.23 /iter & 5 & / & 30327.26  \\  
 \hline
\moma~API  &	628.02 /iter & 46.58 /iter & 0.15 /iter &	4	& 23.20  & 2768.91  \\
\hline
\end{tabular}
\caption{Runtime breakdown for hospital routing instance with 3 wards.}
\label{table:3ward_time}
\end{table}

\begin{figure}
\centering
\includegraphics[width=0.45\textwidth]{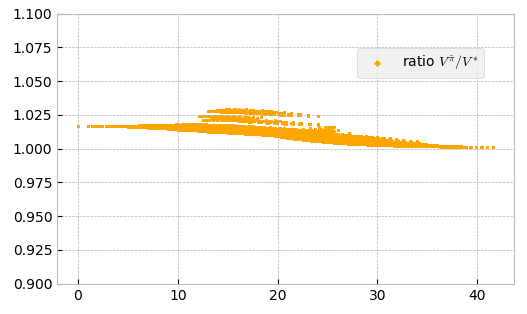}
\includegraphics[width=0.45\textwidth]{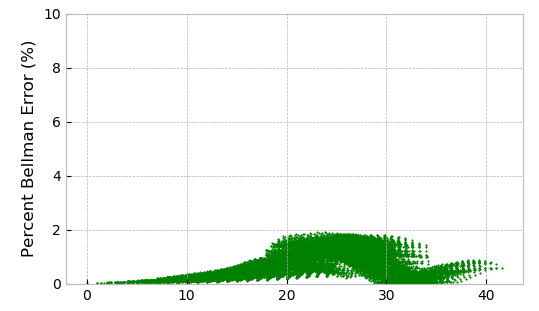}
\caption{(LEFT) Ratio of approximate values compared to convex combination of aggregate values for the constructed grid. (RIGHT) Bellman residual of the function approximation, as percentage of the approximate value function. \label{fig:agg-3ward}}
\end{figure} 

The optimality gap  $\frac{|V^{\tpi}-V^*|}{V^*}$ of \moma~is reported in Figure 
\ref{fig:agg-3ward}(LEFT). For context, the mean optimality gap of 0.92 \% from \moma~API is comparable to the mean of 1.1 \% in \cite{braverman2018taylor}, while the max of 2.97 \% compares quite favorably to the max of 20.6 \% quoted. We repeat this experiment with load of 0.8 and observe similar performance; see Table \ref{table:3ward_opt}.

Useful for larger instances, where exact PI is not computationally feasible, is to use the Bellman residual as a proxy for optimality. Given a candidate policy $\wtilde \pi$, the Bellman residual is  $\wtilde V$, $\lVert T^{\wtilde \pi} \tV - \tV \rVert$. Bounds---that relate the Bellman residual to the optimality gap---appear in \cite{antos2008learning,farahmand2010error}. In this instance the Bellman residuals, as a percentage of the approximate value, have a mean 0.80 \% and a max 1.91 \%; see Figure \ref{fig:agg-3ward} (RIGHT) and Table \ref{table:3ward_opt}.

\begin{table}[h!]
\centering
\begin{tabular}{ |c| c| c| c| c| c| }
\hline
 Load level & TAPI mean error &	TAPI max error & \moma~mean error & \moma~max error & max BR \\ 
 \hline
0.7 &	1.1 \% &	20.6\%	& 0.92\% &	2.97\%  &	1.91\%\\  
 \hline
0.8 &	0.5\% &	9.6\% &	0.91\%	& 3.58\% & 1.04 \% \\
\hline
\end{tabular}
\caption{TAPI and \moma~API error (percent optimality gap) comparison with two load levels, and maximum Bellman Residual (BR) as percent of approximate value function. }
\label{table:3ward_opt}
\end{table}

We use the Bellman residual to experiment with a 4-ward instance where our computer memory resources no longer allow for the computation of the exact policy. 

Here $N_1 = 2, N_2 = 3, N_3 = 1, N_4 = 2$ and the queue capacity is 12 for each. Note that the dimensions are no longer equal in size. Here the state space size is $N=50400$, which translates to the inversion of a matrix with 2.5 billion elements in full evaluation. Using $\fraks = 0.45$ we have $L=1512$. 

\begin{table}[h!]
\centering
\begin{tabular}{| c |c| c| c| c| c| c| c| c| c| }
\hline
  & $\lambda_i$ & $p_i$ & $H_i$  & $B_{i1}$ & $B_{i2}$ & $B_{i3}$ & $B_{i4}$ & $\ell_i$ & $u_i$ \\
 \hline
i= 1 &	0.32 & 0.2 & 10 & / & 5 & 2 & 1 & 0 & 14 \\  
 \hline
i= 2  &	1.68 & 0.7 & 2 & 7 & / & 1 & 2  & 0 & 15 \\
\hline
i= 3  &	0.4 & 0.5 & 6 & 7 & 9 & / & 3 & 0 & 13 \\
\hline
i= 4  &	0.48 & 0.3 & 6 & 1 & 2 & 3 & /  & 0 & 14 \\
\hline
\end{tabular}
\caption{Parameter setting for the 4-ward instance of the hospital routing problem.
% with  separate table whose entry in row $i$ column $j$ shows the overflow cost $B_{ij}$.
}
\label{table:params_hospital_4ward}
\end{table}

\moma~API converged in 5 iterations, taking a total of 5 hours.  Even with sparse representations our memory resources do not allow for policy evaluation. However, to get a sense of the time scale (and hence a benchmark for comparison), we point out that a single full update took 14.6 hours, 17 times as much as the aggregate policy update in \moma~API.

\begin{figure} 
\centering
\includegraphics[width=0.495\textwidth]{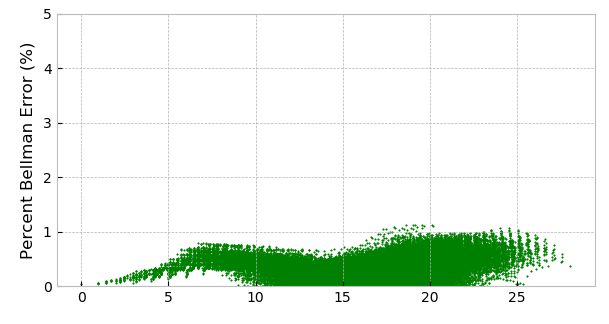}
\caption{ Bellman error as percentage of the function approximation in the 4-ward instance.
\label{fig:agg-4ward}}
\end{figure}

The Bellman residual of the resulting policy is shown in Figure \ref{fig:agg-4ward}, with mean of 0.32 \% and max of 1.13 \% relative to the approximate value function, smaller than even the 3-ward instance, which provides positive indication on the quality of the approximate policy.

\section{Accuracy guarantees\label{sec:guarantees} }

We build on, and expand upon, the results of \cite{braverman2018taylor}. There, the moments $\mu$ 
and $\sigma^2$ are derived directly from $P$ and the focus is on the accuracy and optimality gaps that arise from using the continuous state space PDE to approximate the Bellman equation. Inherent to this then is that the error terms depend on (informally speaking) the third derivative of the solution $\hV$ to the differential equation.  

Here we must correct the bounds for the case of mismatch 
of the second moment between $\tP$ and $P$. When comparing $\tV$ to $\hV$ the accuracy will depend, as well, on bounds on the second derivative which multiplies this mismatch. 

Because the guarantees depend on a PDE, some of the notation and language of that literature is unavoidable. The final result in Theorem \ref{thm:guaranteemain} can, however, be read without familiarity with PDE language and key results. 

To simplify the exposition we assume that the state space of $P$ is unbounded and equal to all of $\mathbb{Z}^d$. In various applications the states cannot go negative. Such ``reflection'' at the boundary causes issues that we ignore for the sake of exposition. These ``gaps'' however are easily completed by reference to \cite{braverman2018taylor}. Also, for computation the state space is often truncated, but this is of secondary importance. We assume that truncation is done at large enough values to have only minimal effect. This is formalized further below. 

The PDE ``induced'' by the equation $0=TV-V$ is given, recall, by 
\be 0 = c(x) + \alpha \mu(x)'DV(x)+\alpha\frac{1}{2}trace(\sigma^2(x)'D^2V(x))-(1-\alpha)V(x).\ee The state space for this PDE is all of $\mathbb{R}^d$. The moment functions $\mu(x)=\Ex_x[X_1-x]$ and $\sigma^2(x)=\Ex_x[(X_1-x)(X_1-x)^{\intercal}]$, and the cost $c(x)$, are defined only for $x\in\mathbb{Z}^d$. With some abuse of notation $c(x)$, $\mu(x)$ and $\sigma^2(x)$ in the PDE are the extensions of these to $\mathbb{R}^d$. Assumption \ref{asum:primitives} below imposes condition on these extensions. {\em Any chain} on $\mathbb{Z}^d$ whose local first and second moments are
given by the functions $\mu$ and $\sigma^2$ induces the same PDE. This PDE is defined relative to $\mu,\sigma^2$ and {\em not} with $\wtilde{\mu},\widetilde{\sigma}^2$. Under our construction in Algorithm \ref{alg:moma}, $\widetilde{\mu}=\mu$, but $\sigma^2\neq \widetilde{\sigma}^2$.  
We mark the PDE solution, when it exists, with $\hV$. The accuracy with which $\hV$ approximates the true value $V$ depends on smoothness properties of $\mu,\sigma^2$ as well on the maximal jumps $|\Delta|_{\calS}^*$. 

For a function $f:\mathbb{R}^d\to \mathbb{R}$, a constant $\vartheta \in (0,1]$ and a set $\calB\subseteq \mathbb{R}^d$, we write $$[f]_{\theta,\calB}=\max_{x,y\in \calB}\frac{|f(y)-f(x)|}{\|y-x\|^{\theta}}.$$ When $\theta=1$, this is the (local) Lipschitz constant over $\calB$ and we drop the subscript $\theta$. We also remind the reader that $x\pm z$ is the set of points $\{y\in\mathbb{R}^d: \|y-x\|\leq z\}$. 
Recall that $c$ is assumed to be norm-like in the sense of \eqref{eq:normlike} and that we take $x_0 = 0$ without loss of generality.

\begin{assumption}[primitives] The primitives $\mu,\sigma^2$ and $c$ satisfy the following assumptions
\begin{enumerate}  \item $\mu$ is globally bounded and Lipschitz and
$$(1-\alpha)^{-1/2}|\mu|_{\br}^*+(1-\alpha)^{-1}[\mu]_{\br}^*\leq \Gamma $$
\item $\sigma^2$ is globally bounded and Lipschitz with
$$|\sigma^2|_{\br}^*+(1-\alpha)^{-1/2}[\sigma^2]_{\br}^*\leq \Gamma, $$ and satisfies 
the ellipticity condition: there exists $\lambda>0$ such that 
$$\lambda^{-1}|\xi|^2\geq \sum_{i,j}\xi_i\xi_j\sigma_{ij}(x)\geq \lambda \|\xi\|^2, \mbox{ for all } \xi,x \in\mathbb{R}^d.$$
\item The cost function $c$ is norm-like and three times differentiable with 
$$|D^ic|_{\calB_r}^* \leq \Gamma(1+r^{k-i}),\mbox{ for } i=0,1,2.$$ 
\end{enumerate} 
\label{asum:primitives}
\end{assumption}

The requirement on $c$ (specifically on $[c]$) is satisfied, for example, if $c(x)=\sum_{i=1}^d c_i(x_i)$ where $c_i(\cdot):\mathbb{R}\to\mathbb{R}$ is a polynomial of degree less than $k$. More importantly, the requirements --- most importantly that on $\mu$ --- specifies a relationship between the drift and the discount factor. This is the relationship that introduces a ``central-limit-theorem-like'' behavior. In its most basic setting, we consider $n$ random variables (and ``horizon'' of length $n$) and scale space by $\sqrt{n}$. Interpreting discounting as a random exponentially distributed horizon---we observe on average $1/(1-\alpha)$ transitions. The requirement on $\mu$ means that the natural scale of the process fluctuation is $(1-\alpha)^{-1/2}$. 

\iffalse
\bcol{Recall the definition} $$V_{-\varepsilon}[x] = \Ex_x\lsb \sum_{t=0}^\infty \alpha^t \frac{|c(X_t)|}{(1+\|X_t\|)^{\varepsilon}}\rsb.$$ \bcol{We state below the main result that leads to Theorem \ref{thm:guaranteemain}}.
\fi

To present our main theorem, we introduce the following definition $$V_{-\varepsilon}[x] = \Ex_x\lsb \sum_{t=0}^\infty \alpha^t \frac{|c(X_t)|}{(1+\|X_t\|)^{\varepsilon}}\rsb.$$

\begin{theorem}[approximation gap] \label{thm:guaranteemain} Suppose the assumptions stated in \ref{asum:primitives} holds, and that $\widetilde{\Delta}_x\lesssim (1+(1-\alpha)^{\frac{1}{4}}\|x\|^{\frac{1-\varepsilon}{2}})$. Then given $\varepsilon \in (0,1)$, for any $\kappa\geq 2 +\varepsilon - \frac{k}{2}$, all $x:\|x\|\geq (1-\alpha)^{-(1+\kappa)}$

]\be |\tV(x)-V(x)|\lesssim  
\frac{1}{\sqrt{1-\alpha}}(V_{-1}(x)+\tV_{-1}(x)+
V_{-\varepsilon}(x)+\tV_{-\varepsilon}(x)).\ee The $\lesssim$ here does not depend on $x,\alpha$. 
\end{theorem}

The next proposition is a key step 

%TODO
\begin{proposition} \label{thm:guaranteemain1} 
Suppose that $ \widetilde{\Delta}_x\lesssim (1+(1-\alpha)^{\frac{1}{4}}\|x\|^{\frac{1-\varepsilon}{2}})$ and that Assumption \ref{asum:primitives} holds. Then, 
$$|V(x)-\hV(x)|\lesssim  \left(\frac{1}{1-\alpha}\right)^{\frac{k+3}{2}} +  \frac{1}{\sqrt{1-\alpha}}V_{-1}(x)+V_{-\varepsilon}(x),$$ and the same holds with $V$ replaced by $\tV$ everywhere.
\end{proposition}

To prove this result we must study the PDE \eqref{eq:PDE} and how well its solution approximates $V$ and $\tV$. The existence and uniqueness of the PDE solution is typically considered on a smooth bounded domain and one must specify values (or derivative conditions) on the boundaries. To this end, let $$\calB_r:=\{x\in\mathbb{R}^d: \|x\|<r\},~~  \tau(r) = \inf\{t\geq 0: X_t\notin \calB_r\}. $$

We will effectively consider a family of PDEs with growing $r$ and establish bounds that do not depend on $r$; taking $r\uparrow \infty$ will produce Theorem \ref{thm:guaranteemain}. Given a radius $r$, $\mathcal{C}^{2,\theta}(\overline{\calB_r})$ is the space of twice continuously differentiable functions $f:\overline{\calB_r}\to \mathbb{R}$ whose second derivative is H\"{o}lder continuous with parameter $\theta$, i.e., $[D^2f]_{\theta,\overline{\calB_r}<\infty}$. The next lemma follows from standard PDE results; see \cite[Theorem 6.14]{gilbarg2015elliptic}.

Define $$\varrho:= \frac{1}{\sqrt{1-\alpha}}, ~~~\brh:=x\pm \varrho.$$

\begin{lemma}[PDE derivative estimates] 
Fix a radius $r$ and suppose that Assumption \ref{asum:primitives}. Then, for any $\theta\in (0,1)$ the PDE \eqref{eq:PDE}, with the boundary condition $\hV(x)=0,x\in \partial \calB_{r}$, has a unique solution $u\in \mathcal{C}^{2,\theta}(\overline{\calB_{r}})$.
 Furthermore, for all $x: x\pm \varrho \in \calB_r$ 
 \be \tag{2nd derivative} |D^2u|_{\brhh}^* \leq \Gamma \lpb \frac{\|x\|^{k-1}}{\sqrt{1-\alpha}}+ \left(\frac{1}{1-\alpha}\right)^{\frac{k+1}{2}}\rpb.\label{eq:2ndDerBound}\ee
The constant $k$ is as in Assumption \ref{asum:primitives} and $\Gamma$ does not depend on $\alpha,r,x$ but may depend on $\vartheta$. \label{lem:interior} \end{lemma} \vspace*{0.2cm} 

The error in the approximation of the value is bounded by the ``integrated'' second derivative up to the stopping time plus the ``tail'' of the value. When one considers only initial states in $x\in\calB_r\subset \calB_{r^2}$ the latter is small.

\begin{lemma} Suppose that $\Delta_x,\widetilde{\Delta}_x\leq \Gamma(1+\sqrt{\|x\|})$ and that Assumption \ref{asum:primitives} holds. Then, fixing $r$ and letting  $u$ be the solution over $\calB_{r^2}$ and $\tau=\tau(r^2)$, we have 
$$|V(x)-\hV(x)|\lesssim 
\Ex_x\left[\sum_{t=0}^{\tau(r^2)}\alpha^t |D^2\hV|_{X_t\pm \Delta_{X_t}}^{\ast}\Delta_{X_t}^2\right]+\Ex_x\left[\sum_{t=\tau(r^2)+1}^{\infty}\alpha^t |c(X_t)|\right],~x\in\calB_r.$$ 

Further, given $\epsilon>0$, we can choose $r_0$ sufficiently large such that
$$\Ex_x\left[\sum_{t=\tau(r_0^2)+1}^{\infty}\alpha^t |c(X_t)|\right]\leq \epsilon,~ x\in \calB_{r_0} $$ and the same holds for $| \tV(x)-\hV(x)|$ and $\tV$ with $\Ex,\Delta$ replaced by $\tEx,\tDelta$. \label{lem:intergratedDer}
\end{lemma}

\noindent {\bf Proof of Proposition \ref{thm:guaranteemain1}.} Since we can make $\epsilon$ arbitrarily small in \ref{lem:intergratedDer}, we simplify exposition by dropping it from further calculations below. 

Plugging \eqref{eq:2ndDerBound} into Lemma \ref{lem:intergratedDer} 
\begin{align}
\nonumber |V(x)-\hV(x)| &\leq  \Ex_x\left[\sum_{t=0}^{\tau}\alpha^t \Gamma \lpb \frac{\|X_t\|^{k-1}}{\sqrt{1-\alpha}}+ \left(\frac{1}{1-\alpha}\right)^{\frac{k+1}{2}}\rpb \Delta_{X_t}^2\right]\\
&\lesssim  \frac{1}{\sqrt{1-\alpha}}\Ex_x\left[\sum_{t=0}^{\tau}\alpha^t\|X_t\|^{k-1}\Delta_{X_t}^2\right]+ 
\left(\frac{1}{1-\alpha}\right)^{\frac{k+1}{2}}\Ex_x\left[\sum_{t=0}^{\tau}\alpha^t\Delta_{X_t}^2\right]\label{eq:decomp}.\end{align}  

Because $\Delta_x\lesssim 1+ (1-\alpha)^{\frac{1}{4}}\|x\|^{\frac{1-\varepsilon}{2}}$,  $\Delta_x^2\lesssim 1+ \sqrt{1-\alpha}\|x\|^{1-\varepsilon}$ so that
\begin{align*} 
\frac{1}{\sqrt{1-\alpha}}\Ex_x\left[\sum_{t=0}^{\tau}\alpha^t\|X_t\|^{k-1}\Delta_{X_t}^2\right]&\leq  \frac{1}{\sqrt{1-\alpha}}\Ex_x\left[\sum_{t=0}^{\tau}\alpha^t\|X_t\|^{k-1}\rsb+ \Ex_x\left[\sum_{t=0}^{\tau}\alpha^t\|X_t\|^{k-\varepsilon}\rsb\\
&\lesssim \frac{1}{\sqrt{1-\alpha}}V_{-1}(x)+V_{-\varepsilon}(x),
\end{align*} 
\begin{align*}
    \left(\frac{1}{1-\alpha}\right)^{\frac{k+1}{2}}\Ex_x\left[\sum_{t=0}^{\tau}\alpha^t\Delta_{X_t}^2\right] & \lesssim \left(\frac{1}{1-\alpha}\right)^{\frac{k+3}{2}} + \left(\frac{1}{1-\alpha}\right)^{\frac{k}{2}}\Ex_x\left[\sum_{t=0}^{\tau}\alpha^t\|X_t\|^{1-\varepsilon}\rsb \\
    & \lesssim \left(\frac{1}{1-\alpha}\right)^{\frac{k+3}{2}}
\end{align*}
because $\|X_t\|$ bounded. The same derivations holds for the sister chain.\footnote{Notice that the assumptions of Lemma \ref{lem:intergratedDer} are satisfied because $\frac{1-\varepsilon}{2}< \frac{1}{2}$. For the focal chain, since jumps are bounded, we can take $\varepsilon=1$.} \eProof

Now we are ready to use Proposition \ref{thm:guaranteemain1} to establish the main approximation gap result.

\noindent {\bf Proof of Theorem \ref{thm:guaranteemain}.}
 We will show that for any $\kappa\geq 0:\kappa \geq 2 +\varepsilon - \frac{k}{2}$ and $x:\|x\|\geq (1-\alpha)^{-(1+\kappa)}$,
$$\left(\frac{1}{1-\alpha}\right)^{\frac{k+3}{2}}\lesssim V_{-\varepsilon}(x).$$ 

To see this, let $\bar{\Delta}=\sup_{x}\Delta_x$. Then, for all $t\leq t_0(\alpha)=\frac{1}{2\bar{\Delta} (1-\alpha)}$ and $x:\|x\|\geq (1-\alpha)^{-(1+\kappa)}$ $$\|X_t\|\geq \|x\|-\bar{\Delta}t_0\geq \frac{1/2}{(1-\alpha)^{1+\kappa}}.$$ 

Consider $\Ex_x \sum_{t=0}^{t_0}\alpha^t \frac{|c(X_t)|}{(1+\|X_t\|)^{\varepsilon}}\leq V_{-\varepsilon}(x)$. For any trajectory $X_t$, 
% true for $\lVert X_t \rVert \geq 1$; but the additive term should becomes 0 when c(x) also 0
$$\sum_{t=0}^{t_0}\alpha^t \frac{|c(X_t)|}{(1+\|X_t\|)^{\varepsilon}}\gtrsim \sum_{t=0}^{t_0}\alpha^t \left(\frac{1/2}{(1-\alpha)^{(1+\kappa)}}\right)^{k-\varepsilon}\gtrsim  \left(\frac{1}{1-\alpha}\right)^{k(1+\kappa)-\kappa \varepsilon-\varepsilon}.$$  
The last inequality follows from the fact that, as $\alpha\uparrow 1$,  $$\sum_{t=t_0(\alpha)}^{\infty}\alpha^t = \alpha^{t_0(\epsilon)}\sum_{t=0}^{\infty}\alpha^t = \frac{1}{1-\alpha}\alpha^{\frac{1}{2\bar{\Delta}(1-\alpha)}},$$ 
and noting that $\alpha^{\frac{1}{2\bar{\Delta}(1-\alpha)}}\rightarrow e^{-\frac{1}{2\Delta}}<1$, we have $$\sum_{t=t_0(\alpha)}^{\infty}\alpha^t \leq \gamma$$ for some $\gamma<1$ that does not depend on $\alpha,x$.

Finally, 
$$\left(\frac{1}{1-\alpha}\right)^{\frac{k+3}{2}}\lesssim \left(\frac{1}{1-\alpha}\right)^{k(1+\kappa)-\kappa \varepsilon-\varepsilon}\lesssim V_{-\varepsilon}(x),$$ 
if $\kappa \geq 2 +\varepsilon - \frac{k}{2}$. \eProof 

\iffalse 
\noindent {\bf The connection between the PDE-based bounds and $\delta$.} The PDE-based results do not translate easily to bounds on $\delta$ but they are qualitatively informative. Recall Theorem \ref{thm:backtodelta} where, starting from $\delta$, we developed an informal bound in terms of the derivative of extension of $V$ and $\tV$. Treating Replacing $\Delta^2V$ and $\Delta^2\tV$ by $\Delta^2\hV$: 
$$|V(x)-\tV(x)|\lesssim \frac{1}{1-\alpha}\sum_{y}p_{xy}(\Delta_y^2+\widetilde{\Delta}_y^2)\|D^2\hV(y)\|.$$

By Lemma \ref{lem:interior}, $$\widetilde{\Delta}_x^2\|D^2\hV\|_{x\pm \widetilde{\Delta}_x}^*\lesssim \left(\frac{c_{-(1+\varepsilon)}(x)}{1-\alpha}+c_{-\varepsilon}(x)\right).$$ Overall, we have 
$$|V(x)-\tV(x)|\lesssim \frac{1}{1-\alpha}\left(\frac{c_{-(1+\varepsilon)}(x)}{1-\alpha}+c_{-\varepsilon}(x)\right).$$

This is of the same order of magnitude as the bound we established rigorously. Indeed, for all $x:\|x\|\geq \frac{1}{1-\alpha}$, $c_{-(1+\varepsilon)}(x)\leq c_{-\varepsilon}(x)$ and because $P$ has bounded jumps $\frac{1}{1-\alpha}c_{-\varepsilon}(x)\lesssim V_{-\varepsilon},$ for all such $x$ and we arrive at the conclusion that \begin{equation}|V(x)-\tV(x)|\lesssim V_{-\varepsilon}(x),\mbox{ for all } x:\|x\|\geq \frac{1}{1-\alpha}.\label{eq:lastbound} \end{equation} 

\textcolor{red}{the back and forth with $\delta$ is still a bit weak.}\fi 

We conclude this section with a reference back to the \moma~algorithm.

\begin{theorem}[optimality gap of \moma~API]
    Consider a controlled chain on $\mathbb{Z}^d$. Let $\pi^*$ be the optimal policy and $\widetilde \pi$ be the \moma~policy (produced by Algorithm \ref{alg:onestep}). Suppose that Assumption \ref{asum:primitives} holds for  $c(x,\pi(x)),\mu_{\pi(x)}(x),\sigma^2_{\pi(x)}(x)$ for both $\pi \in \{ \tpi, \pi^*\}$.  Then, 
    $$V^{ \widetilde \pi}(x) - V^*(x) \lesssim  V^*_{-\varepsilon}(x)+ V^{\tpi}_{-\varepsilon}(x).$$
  
    \label{thm:aggOptGap}
\end{theorem} 

\bProof By Theorem \ref{thm:lifted}, the policy $\tpi$ produced by \moma~API is optimal for the lifted chain. In particular, $\tV^*=\tV^{\tpi}\leq \tV^{\pi^*}$. Then by Theorem \ref{thm:guaranteemain} 
\begin{align*}  V^{\tpi}(x)-V^{\pi^*}(x)&\leq V^{\tpi}(x)-\tV^{\tpi}(x)+\tV^{\tpi}(x)-\tV^{\pi^*}(x)+\tV^{\pi^*}(x)-V^{\pi^*}(x)\\ & \leq V^{\tpi}(x)-\tV^{\tpi}(x)+\tV^{\pi^*}(x)-V^{\pi^*}(x)\\& \lesssim V^{\pi^*}_{-\varepsilon}(x)+ V^{\tpi}_{-\varepsilon}(x).
\end{align*} 
In the second inequality we used $\tV^{\tpi}(x)\leq \tV^{\pi^*}(x)$ and in the third we applied (twice) the bounds from Theorem \ref{thm:guaranteemain}. \eProof

\section{Concluding remarks\label{sec:concluding}} 

What we offer here is an approach to ADP that achieves a synergy between a central-limit theorem view of control (the matching of moments) and a well-established algorithmic building block (aggregation). Our paper brings algorithmic relevance to some theoretical ideas introduced in \cite{braverman2018taylor}, which itself builds on a long list of papers on CLT-based approximations. 

The key here is the identification of aggregation as a stepping stone on which to build implementable algorithms that can be ``matched'' with the theory that approximates a discrete problem with a continuous one. The re-interpretation of aggregation as creating a new Markov chain on the original state space, gives rise to a flexible infrastructure on which to superimpose moment matching. The approximation is grounded in math that informs a choice of the aggregation parameters that is consistent with optimality guarantees. 

\noindent {\bf Next steps.} \begin{enumerate} 
\item {\bf m-step moment matching.}  Our choice of $\bg \bu$, recall, is such that
$$\sum_{z}(\bg \bu)_{yz} = y = \Ex_y[X_0].$$ The second equality is trivial---the chain at time $t=0$ is at its initial state $y$ and serves to re-interpret what we did here. By choosing $\bg \bu$ that matches the first moment at time ${\bf t=0}$, we guarantees that the $\tP = P\bg \bu$ matches the moment at time ${\bf t=1}$: 
$$\sum_{z} (P\bg \bu)_{yz} = \Ex_y[X_1] = y+\mu(y).$$ 

More generally, if one chooses $\bg,\bu$ to match the first moment at time $t=m-1$, i.e., $\sum_{z}(\bg \bu)_{yz} = \Ex_y[X_{m-1}]$, then $\tP$ matches the $m$-step moment:  
$$\sum_{z} (P\bg \bu)_{yz} = \Ex_y[X_m].$$ 

Why might this be valuable? In the coarse grid approximation it is inevitable that the chain $\tP$ has large jumps relative to $P$ itself. The coarser the grid the larger the jumps of $\tP$. This affects the second moment mismatch which, in turn, affects the quality of the approximation. Because it lumps multiple transitions together, the $m$-step chain $P^m$ has larger jumps than $P$ that are possibly better aligned with those of $\tP$. This, in some settings, may improve performance. To make this point more concrete, some initial developments are offered in \S \ref{sec:mstep} of the appendix. 
\vspace*{0.2cm} 

\item {\bf Approximation-Estimation tradeoff and a hierarchy of models:} Consider a setting where $P$ is not known in advance but rather estimated based on observations. As transitions are performed, the estimate of the matrix $P$ improves. It seems intuitively reasonable that the approximate model, based as it is on a Taylor approximations, will induce smaller variance over finite samples but, because of the approximation, will have a larger bias. Given the uncertainty early in the horizon it may, nevertheless, make sense to use a coarser and computationally cheaper model and gradually transition -- as more samples are collected --- to a more accurate model. Within our framework, such a gradual transition, is enabled naturally by the coarseness (spacing) exponent $\mathfrak{s}$. 
\iffalse 
That this bias/variance tradeoff is to be expected and is made sufficiently evident in a simple static setting. Take a random variable $Z$ and the expectation $\Ex[f(x+Z)]$ where $f(y)=y^4$. Consider the second-order Taylor approximation 
$\tilde{f} (x+y)  = f(x)+ f'(x)y + \frac{1}{2} f''(x)y^2,$ or, with expectations, $$\Ex[f(x+Z)]\approx  f(x) + f'(x)\Ex[Z] +\frac{1}{2}f''(x)\Ex[Z^2].$$ 
Obviously $\Ex[f(x+Z)]\neq \Ex[\tilde{f}(x+Z)],$ and the approximation error is of the order of $x^3\Ex[Z^3]$. Suppose that we have a finite number of samples $z_1,\ldots,z_n$ and we compute our estimates as $$\widehat{\Ex}[f(x+Z)]=\frac{1}{n}\sum_{k=1}^n f(x+z_k), ~ \widehat{\Ex}[\tilde{f}(x+Z)] = \frac{1}{n}\sum_{k=1}^n\tilde{f}(x+z_k).$$

If, for example, $x\geq 0$ and $Z$ is positive valued that it is easily verified here that $$\sqrt{n}\sigma_{\tilde{f}} = \sqrt{Var(\widehat{\Ex}[\tilde{f}(x+Z)]} < 
\sqrt{n}\sigma_{f} = \sqrt{Var(\widehat{\Ex}[f(x+Z)]}.$$ \fi 
The bias-variance view offers then a lens through which to explore the interaction of ADP approximation and parameter estimation, one that is natural within the Tayloring/moment matching framework we put forth here. 
\vspace*{0.2cm}
\item {\bf Aggregation based RL.}
Aggregation has already been espoused as an aid in simulation-based policy iteration in \cite{bertsekas2018feature}. Useful in our algorithm is the fact that the matrix $\bg$ does not depend on the transition probability matrix $P$ so that, in contrast to other architectures, it requires no updating within iterations. This may have useful implications for sample complexity. From an analysis perspective, the coarse grid approximation is appealing. Because the construction leads itself to a Markov chain (a lower rank one) on a subset of states, one can build on existing studies; see e.g. \cite{haskell2016empirical}. 
\end{enumerate} \vspace*{0.2cm} 

\noindent \paragraph{Acknowledgements.} This work is supported by NSF grant CMMI-1662294. We thank Anton Braverman and Nathalie Vanvuchelen for generously sharing their dynamic programming codes for the examples in our \S \ref{sec:numerical}

\iffalse 

\begin{figure}[h!] 
\begin{center} 
\includegraphics[scale=0.35]{BiasVariance.pdf}
\end{center} 
\label{fig:taylor}
\caption{\rcol{I don't like this figure standing here like a stand alone. At any rate, caption is missing} } 
\end{figure} 

Fixing $x=25, \mu=\sigma=5$ we compare the mean and the variance (computed via simulation with 10000 replications) of the two estimates. The result is plotted in Figure \ref{fig:taylor}. What is important is that the standard deviation of the estimate--- $\sigma_{\tilde{f}}$ --- is smaller that of $\sigma_f$ that correspond to the estimate for $f$. 
\fi 

% \bibliographystyle{informs2014}
% \bibliography{central1,central2}

\newpage
\renewcommand*{\thesection}{A.\arabic{section}}
\setcounter{section}{0}
\section*{Appendix}

% \newpage
% \section{Computation Configuration}
% All computations were done on machine with Intel(R) Core(TM) i7-6700 CPU @ 3.40GHz 3.41 GHz and 16.0GB of RAM, using 64-bit Python.
\section{Proofs of Auxiliary Lemmas} \label{sec:proofs}

\noindent {\bf Proof of Lemma \ref{lem:convolution}}
 In that case, $(P\bg \bu)_{x\cdot}$ is the convolution of $P_{x\cdot}$ and a zero-mean jump and, consequently, has the same mean as $P_{x\cdot}$: 
$\bg \bu M_1(y) = y \implies P\bg \bu M_1(x) = x+\mu(x).$
\eProof

\noindent {\bf Proof of Lemma \ref{lem:phi_construct}. } 
For fixed $y \in \calS$ and $x_l$ in its enclosing box $\calB^d$, we defined
\[ \sg_{yl} = \Pi_{i=1}^d \left( \mathrm{1}\{[x_l]_i = \bar{s}_i \}  \frac{y_i - \underline{s}_i}{\bar{s}_i-\underline{s}_i} + \mathrm{1}\{ [x_l]_i = \underline{s}_i\}  \frac{\bar{s}_i - y_i}{\bar{s}_i-\underline{s}_i}  \right) \] 
where $\bar{s}_i, \underline{s}_i$ are the upper and lower values along axis $i \in [d]$ for corners of the hyperbox. There are $2^d$ such corners.  To simplify notation, let 
$$h_i(l)=\left\{\begin{array}{ll} 1 &\mbox{ if } (x_l)_i = \bar{s}_i,\\ 0&\mbox{ if} (x_l)_i = \underline{s}_i.\end{array}\right. $$Also, set $w_i(y) = \frac{y_i - \underline{s}_i}{\bar{s}_i-\underline{s}_i}.$ so that $0 \leq w_i(y) \leq 1$ is such that $y_i=w_i (y) * \bar{s}_i + (1-w_i(y)) * \underline{s}_i = y_i $. Finally, define $F_i(l) = h_i(l) * w_i(y) + (1 - h_i(l))* (1-w_i(y)),$ so that $$\sg_{yl} = \Pi_{i=1}^d F_i(l).$$ Since $F_i(l) \in [0,1]$ also $\sg_{yl} \in [0,1]$. 

We argue by induction that $\sum_{l\in \calB^d}g_{yl}=1$. For $d=1$, $\sum_{l \in \calB^d}\sg_{yl} = w_1(y) + (1-w_1(y)) = 1$. Now suppose for $d'=d-1$, $\sum_{l \in \calB^d} \Pi_{i=1}^d F_i(l) = 1$. Suppose also, w.l.o.g, that the indices are ordered such that for $k=1, .., 2^{d'-1}$, $l_{2k-1}, l_{2k} \in \calB^{d'}$ differ only on axis $d$, specifically $h_{d'}(l_{2k-1})=1, h_{d'}(l_{2k})=0$.

Then we have 
\begin{align*}
    \sum_{l \in \calB^{d'}} \sg_{yl} 
    = & \sum_{l \in \calB^{d'}} \Pi_{i=1}^{d'} F_i(l)\\
    = & w_{d'}(y) \Pi_{i=1}^{d'-1} F_i(l_1) +  (1-w_{d'}(y)) \Pi_{i=1}^{d'-1} F_i(l_2) + ... +  w_{d'}(y) \Pi_{i=1}^{d'-1} F_i(l_{2^{d'-1}}) +  (1-w_{d'}(y)) \Pi_{i=1}^{d'-1} F_i(l_{2^{d'}})\\
    = & w_{d'}(y) \sum_{k =1}^{2^{d'-1}} \Pi_{i=1}^{d'-1} F_i(l_{2k-1}) + (1-w_{d'}(y)) \sum_{k =1}^{2^{d'-1}}  \Pi_{i=1}^{d'-1} F_i(l_{2k}) \\
    =& \left( w_{d'}(y)+1-w_{d'}(y) \right) \sum_{l' \in \calB^{d'-1}} \Pi_{i=1}^{d'-1} F_i(l') \\ 
    = & 1
\end{align*}
where the last equality is due to the inductive assumption, completing the induction. 

Now we want to show $\sum_{l \in \calB^d} \sg_{yl} x_l = y$, i.e. 
\[ \sum_{l \in \calB^d} [x_l]_j \Pi_{i=1}^d \left( \mathrm{1}\{[x_l]_i = \bar{s}_i \} * \frac{y_i - \underline{s}_i}{\bar{s}_i-\underline{s}_i} + \mathrm{1}\{ [x_l]_i = \underline{s}_i\} * \frac{\bar{s}_i - y_i}{\bar{s}_i-\underline{s}_i}  \right) = [y]_j\]
for $j = 1, ..., d$. Notice that $[x_l]_j = h_j(l)\bar{s}_j + (1- h_j(l)) \underline{s}_j$, so when $d=1$, the quantity on the left hand side simply evaluates to $LHS = \bar{s}_1 w_1(y) + \underline{s}_1 (1-w_1(y))$, satisfying the equality. 

For $d \geq 2$, 
\begin{align*}
LHS = & \bar{s}_j w_j(y) \sum_{\mathclap{\substack{l \in \calB^d\\ h_j(l) =1}}} \Pi_{i \neq j} F_i(l) + \underline{s}_j (1-w_j(y)) \sum_{\mathclap{\substack{l \in \calB^d\\ h_j(l) =0}}} \Pi_{i \neq j} F_i(l) \\
= & \left( \bar{s}_j w_j(y) + \underline{s}_j (1-w_j(y)) \right) \sum_{l \in \calB^{d-1} } \sg_{yl} \\
= & y_j
\end{align*}
using the first part. Since this is true for any $j = 1, ..., d$, this concludes the proof. 

\eProof

\noindent {\bf Proof of Theorem \ref{thm:backtodelta}.}
Recall from \S \ref{sec:tayloring} that $\delta_f[P,\wtilde{P}] :=  |\delta_f[P,\wtilde{P}] (\cdot)|_{\calS}^*$  
where 
$$\delta_f[P,\wtilde{P}](x):=   \lvert \wtilde \Ex_x[f(X_1)]-\Ex_x[f(X_1)]\rvert = \lvert   \wtilde{P}{f}(x) - Pf(x)\rvert.$$
With the coarse grid scheme, it takes on the form \begin{align*} \delta_V[P,\tP](x)&=|PV(x)-\tP V(x)|
=|\sum_{y}p_{xy}V(y)-\sum_{y}p_{xy}\sum_{l}\sg_{yl}V(x_l)|\\ & 
\leq \sum_{y}p_{xy}|V(y)-\sum_{l}\sg_{yl}V(x_l)|\end{align*} 
so that $\delta_V[P,\tP]\leq \max_{x \in \mathcal{S}} \sum_{y}p_{xy}|V(y)-\sum_{l}\sg_{yl}V(x_l)|.$ Then by equation \eqref{eq:tilde_operator_gap1}, we have 
\begin{align*}
    |V-\tV|_{\calS}^*\leq & \max_{x \in \mathcal{S}} \frac{1}{1-\alpha}\sum_{y}p_{xy}\Big(|V(y)-\sum_{l}\sg_{yl}V(x_l)|+|\tV(y)-\sum_{l}\sg_{yl}\tV(x_l)|\Big) \\
    \leq & \frac{1}{1-\alpha}\left(|V-GW|_{\calS}^* +|\tV-G\tW|_{\calS}^*\right)
\end{align*}
Equation \eqref{eq:tilde_operator_gap1} is a special case of equation \eqref{eq:opt_delta} obtained, trivially, by assuming that the actions space $\mathcal{A}(x)$ contains a single action (so that $\pi^*=\tpi^*$). We refer the reader to that proof later in this appendix.
\eProof

\noindent {\bf Proof of Equation (\ref{eq:recover_bound}).}
\begin{align*}
    \lvert V(x) - R(k(x)) \rvert = & \lvert V(x) - \bg \bu \tV(x) \rvert \\
    \leq & \lvert V(x) - \bg \bu V(x) \rvert + \lvert \bg \bu V (x) - \bg \bu \tV (x) \rvert \\
    \leq & \lvert V(x) - \bg \bu V (x) \rvert + \max_{k \in \calM} \lvert \bu V(k) - \bu \tV (k) \rvert 
\end{align*}
Following a standard argument similar to Equation \ref{eq:tilde_operator_gap1}, we can show that $\lvert \bu V(k) - \bu \tV (k) \rvert \leq \frac{\alpha}{1-\alpha} \lvert V - \bg \bu V \rvert_{\calS}^* $. In turn, 
\begin{align*}
    \lvert V(x) - R(k(x)) \rvert \leq & \frac{1}{1-\alpha} \lvert V - \bg \bu V \rvert_{\calS}^* \\
    \leq & \frac{1}{1-\alpha} \max_{y_1, y_2 \in \calS_k} \lvert V(y_1) - V(y_2) \rvert.
\end{align*}

To see why $\lvert \bu V(k) - \bu \tV (k) \rvert \leq \frac{\alpha}{1-\alpha} \lvert V - \bg \bu V \rvert_{\calS}^* $, let $\epsilon = \lvert V - GUV \rvert_{S}^*$, and $\bar{R}(l) = U V(l) + \frac{\alpha}{1-\alpha} \epsilon$, then
\begin{align*}
    H \bar{R}(l) =& \sum_{x \in \calS} u_{lx} [c(x) + \alpha \sum_{y \in \calS} p_{xy} \sum_{k \in \calM} g_{yk} \bar{R}(k)] \\
    =& \sum_{x \in \calS} u_{lx} [c(x) + \alpha \sum_{y \in \calS} p_{xy} \sum_{k \in \calM} g_{yk} (U V(k) + \frac{\alpha}{1-\alpha} \epsilon)] \\
    =& \frac{\alpha^2}{1-\alpha} \epsilon + UV(l) + \sum_{x \in \calS} u_{lx} \alpha \sum_{y \in \calS} p_{xy} \big( \sum_{k \in \calM} g_{yk} (U V(k) - V(y) \big) \\
    \leq & UV(l) + \frac{\alpha}{1-\alpha} \epsilon \\
    =& \bar{R}(l)
\end{align*}
Thus, we have that $H \bar{R}(l) \leq \bar{R}(l)$, from which it follows that $R(l) \leq \bar{R}(l) = U V(l) + \frac{\alpha}{1-\alpha} \epsilon$. 
The other side follows similarly. 
\eProof

 \noindent {\bf Proof of Lemma \ref{lem:gridfn}}
 Define $\calN_x = \{y: p_{xy}>0 \}$ and similarly $\wtilde \calN_x = \{y: \wtilde p_{xy}>0 \}$. With our coarse grid scheme, $\wtilde \calN_x = \cup_{y \in \calN_x} \calB(y)$, so we have 
$$ \max_{z \in \wtilde \calN_x} \lvert z_i - x_i \rvert \leq \max_{y \in \calN_x}  \lvert y_i - x_i \rvert + \max_{y \in \calN_x} \gridfn(\lvert y_i \rvert) $$ 
for $\gridfn(\cdot)$ increasing. Then $ \wtilde \Delta_x \leq \Gamma \gridfn(\lVert x \rVert + \Delta_x ) $ for some constant $\Gamma$, and taking $\gridfn(z) = z^{\fraks}$, $0 \leq \fraks < \frac{1}{2}$, is sufficient to satisfy the required condition for $\varepsilon = 1 - 2 \fraks$.
\eProof

\noindent {\bf Proof of Lemma \ref{lem:secondapprox}.}
Fix $y\in\calS$ and let $\hat{x}^1,\ldots,\hat{x}^{2^d}$ be the corners of the box $\calB$ that contains $y$. The second order Taylor expansion of the function $M_2(z)=zz^{\intercal}$ has for $k = 1, ..., 2^d$
\begin{align*} 
M_2(\hat{x}^k) & = M_2(y)+ DM_2(x)'(\hat{x}^k-y)\pm \frac{1}{2}\Gamma \|\hat{x}^k-y\|^2,
\end{align*} 
where we use the fact that $\|D^2 M_2\|\leq \Gamma(d)$. Taking the corner of the box that is coordinate-wise closest to the origin, and call it $k_0$, every point $y$ in the box satisfies $y_i\in [\hat{x}^{k_0}_i,\hat{x}^{k_0}_i+\gridfn(|\hat{x}^{k_0}_i|)]$ 
so that $\|\hat{x}^k-y\|^2\leq d \max_{i}\gridfn^2(|\hat{x}^{k_0}_i|) \leq d \gridfn^2(\|y\|) $

Since by construction, $\sum_{l}\sg_{yl}x_l=y$, we have that $\sum_{l}\sg_{kl}(x_l-y)=0$. In turn, there exists some constant $\Gamma$ independent of $x,\alpha$ such that
$$\tEx_x[X_1 X_1 ^{\intercal}] =\sum_{y,l}p_{xy} \sg_{yl} (y + x_l -y) (y + x_l -y)^{\intercal} = \Ex_x[X_1 X_1 ^{\intercal}]\pm \Gamma \gridfn^2(\|y\|).$$ Finally, since $\max_{y:p_{xy}>0}\gridfn^2(\|y\|)\leq \Gamma \gridfn^2(\|x\|+\Delta_x)$, we conclude that
$$ \lVert \wtilde{\sigma}^2(x) - \sigma^2(x) \rVert \leq \Gamma \gridfn^2(\|x\|+\Delta_x),$$ as required. 
\eProof

\noindent {\bf Proof of Lemma \ref{lem:grid_growth}.} Recall that $\calS = \times_{i=1}^d [\ell_i,u_i]$. Consider first the case that $\ell_i\geq 0$. Fix $\fraks\in [\frac{1}{3},\frac{1}{2})$ and define for $k\in\mathbb{Z}_+$ the recursion
$f(0)=\ell_i$ and \be f(k+1) = \lceil f(k)+f^\fraks(k) \rceil + 1,\label{eq:frecursion}\ee that defines, recall \eqref{eq:f_recursion}, the grid points on the $i^{th}$ axis. 
The number $k_f^* = \inf \{ k: f(k) \geq  u_i \}$ is then the number of grid points on this axis. We bound this number by considering a continuous lower bound on $f$. 

Defined for $x\geq 0$ the function $h(x) = ((1-\fraks)x)^\frac{1}{1-\fraks} + \ell_i$. Its Taylor expansion has form
\[ h(k+1) = h(k) + ((1 -\fraks) k)^{\frac{\fraks}{1-\fraks}} \pm \frac{1}{2}\fraks ((1-\fraks)k)^{\frac{2\fraks-1}{1-\fraks}}=h(k)+h^{\fraks}(k)\pm \frac{1}{2}\fraks ((1-\fraks)k)^{\frac{2\fraks-1}{1-\fraks}}.\]
Note that $2\fraks-1 < 0$ for $\fraks < \frac{1}{2}$, so we have $((1-\fraks)k)^{\frac{2\fraks-1}{1-\fraks}} \leq 1$, thus the last term is upper bounded by $\frac{1}{4}$. Then combined with \eqref{eq:frecursion}, we know that if there is a $k \in \mathbb{Z_+}$ such that $h(k) \leq f(k)$, we have $h(k') \leq f(k')$ for all $k' \in \mathbb{Z_+}$ where $k' \geq k$.

Since $h(0)=f(0)=\ell_i$, it is clear that $h(k)\leq f(k)\mbox{ for all }k\in\mathbb{Z}_+.$ We argue the slightly stronger claim
\be  h(x) \leq f(\lfloor x \rfloor ) \mbox{ for $x \geq 2$.} \label{eq:hfclaim}\ee Let us use \eqref{eq:hfclaim} to complete the proof of the lemma. 

Defining now $k_h^* := \inf \{ k: h(k) =  u_i \}$ we have $f(\lfloor k_h^* \rfloor) \geq h(k_h^*) = u_i $, thus $k_f^* \leq k_h^*$. Since $h$ is a continuous increasing function $h(k_h^*) = \ell_i+((1-\fraks) k_h^*)^{\frac{1}{1-\fraks}} = u_i$ so that 
$ k_h^* = (u_i - \ell_i)^{1- \fraks}/(1-\fraks),$ and, in turn, 
$$k_f^* \leq \frac{(u_i - \ell_i)^{1- \fraks}}{1-\fraks}.$$ 

The case that the $i^{th}$ axis has $u_i\leq 0$. Define $f_{-}(0)=u_i$ and recursively $f_{-}(-k) = \lfloor f_{-}(-(k-1))-(-f_{-}(-(k-1)))^\fraks \rfloor - 1$. It follows identically that $$k_{f_{-}}^* \leq \frac{(- \ell_i - (-u_i))^{1- \fraks}}{1-\fraks}.$$

If $\ell_i\leq 0$ and $u_i\geq 0$, we treat the negative portion $[\ell_i, 0]$ and the positive portion $[0, u_i]$ separately to obtain that the number of grid points on the $x$ axis satisfies  $k^*\leq \frac{u_i^{1- \fraks} + (-\ell_i)^{1- \fraks}}{1-\fraks}.$ 
In all cases then, $k$ is upper bounded by 
\begin{align*}
     \prod_i \{ \mathrm{1}_{\{u_i > 0\}} \frac{u_i^{1 - \fraks}}{1-\fraks}  + \mathrm{1}_{\{\ell_i < 0\}} \frac{(-\ell_i)^{1 - \fraks}}{1-\fraks}   \} 
   \leq  \prod_i  \frac{2(\frac{u_i-\ell_i+1}{2})^{1-\fraks}}{1- \fraks} 
    \leq  \prod_i \frac{\sqrt{2} r_i^{1-\fraks}}{1- \fraks} 
    =  \left(\frac{\sqrt{2}}{1- \fraks}\right)^d N^{1-\fraks},
\end{align*} where we used $N=\Pi_{i}r_i$. 

It remains only to prove \eqref{eq:hfclaim}. Because $h$ is decreasing in $\fraks$, the maximum value over $\fraks\in [\frac{1}{3},\frac{1}{2})$ is achieved at $\fraks = \frac{1}{3}$. For the basis of the induction, notice that if $\ell_i =0$, $f(1) = 1$ and $f(2) = 3$, then for $x\in [2,3)$, 
\[ h(x) = ( (1-\fraks)x )^ {\frac{1}{1-\fraks}} \leq 2^{\frac{3}{2}}  < 3 = f(\lfloor x \rfloor).  \] If 
$\ell_i \geq 1$, $f(1) \geq \lfloor \ell_i + \ell_i^\fraks \rfloor +1\geq \ell_i + 2$ so that, again,
\[ h(x) = ( (1-\fraks)x )^ {\frac{1}{1-\fraks}} \leq \left(\frac{4}{3}\right)^{\frac{3}{2}} + \ell_i < 2 + \ell_i < f(\lfloor x \rfloor).  \] 
Suppose now that \eqref{eq:hfclaim} holds up to $k_1 \in \mathbb{Z}_+$, and for $x\in [k_1, k_1+1)$. Then, for such $x$, we have by the induction assumption that $h(x)\leq f(k_1)$ and, in turn, for $y=x+1\in [k_1+1,k_1+2)$
\[ h(y) = h(x) + h^{\fraks}(x) \pm \frac{1}{4} \leq f(k_1) + f^{\fraks}(k_1) \pm \frac{1}{4} < f(k_1+1)=f(\lfloor y\rfloor), \] and this completes the induction.  
\eProof

\noindent {\bf Proof of Equation \eqref{eq:opt_delta}. } The proof follows a standard argument; see e.g. \cite[Proposition 6.2]{bertsekasneuro}. Let $J_1(x) = V^{ \ast }(x) + \frac{\alpha}{1-\alpha}\delta_{V^{\pi^*}}[P^{\pi^*}, \wtilde{P}^{\pi^*}].$
\begin{align*}
\wtilde T J_1(x) &=  \min_{a \in \mathrm{A}(x)} \left[ c(x,a) + \alpha  \sum_{y \in S} \wtilde p_{xy}^a J_1 (y) \right] \\
&\leq \min_{a \in \mathrm{A}(x)} \left[ c(x,a) + \alpha  \sum_{y \in S} \wtilde p_{xy}^a \{ V^{ \ast }(y) + \frac{\alpha}{1-\alpha}\delta_{V^{\pi^*}}[P^{\pi^*}, \wtilde{P}^{\pi^*}] \} \right]\\
&=  \min_{a \in \mathrm{A}(x)} \left[ c(x,a) + \frac{\alpha ^2}{1-\alpha}\delta_{V^{\pi^*}}[P^{\pi^*}, \wtilde{P}^{\pi^*}] +\alpha \sum_{y \in S} [p_{xy}^aV^{ \ast }(y) +  \wtilde p_{xy}^aV^{ \ast }(y) -  p_{xy}^aV^{ \ast }(y) ] \right] \\
&\leq c(x,\pi ^ \ast(x))  +\alpha  P^{\pi ^ \ast}V^{ \ast }(x) + \alpha \mid  \wtilde P^{ \pi ^ \ast}V^{ \ast }(x) - P^{\pi ^ \ast}V^{ \ast }(x) \mid  +  \frac{\alpha ^2}{1-\alpha} \delta_{V^{\pi^*}}[P^{\pi^*}, \wtilde{P}^{\pi^*}]\\
&\leq V^{ \ast }(x) + \frac{\alpha}{1-\alpha}\delta_{V^{\pi^*}}[P^{\pi^*}, \wtilde{P}^{\pi^*}] = J_1(x)
\end{align*}
The above shows that $ J_1(x) \geq \wtilde T J_1(x) $. Since $\wtilde T J_1(x)  \to \wtilde V(x)$, we have $J_1(x) \geq \wtilde V(x)$ by monotonicity.
Now repeat for the other side and let $J_2(x) = \wtilde V ^{\ast}(x) + \frac{\alpha}{1-\alpha}\delta_{\tV^{\tpi}}[P^{\tpi}, \wtilde{P}^{\tpi}].$
\begin{align*}
T J_2(x) &=  \min_{a \in \mathrm{A}(x)} \left[ c(x,a) + \alpha  \sum_{y \in S} p_{xy}^{ u}J_2 (y) \right] \\
&= \min_{a \in \mathrm{A}(x)} \left[ c(x,a) + \alpha  \sum_{y \in S} p_{xy}^a [ \wtilde V ^{\ast}(y) + \frac{\alpha  }{1-\alpha}\delta_{\tV^{\tpi}}[P^{\tpi}, \wtilde{P}^{\tpi}] ] \right] \\
&=  \min_{a \in \mathrm{A}(x)} \left[ c(x,a) + \frac{\alpha ^2  }{1-\alpha}\delta_{\tV^{\tpi}}[P^{\tpi}, \wtilde{P}^{\tpi}] +\alpha \sum_{y \in S} [ \wtilde p_{xy}^a\wtilde V ^{\ast}(y) + p_{xy}^a\wtilde V ^{\ast}(y) - \wtilde p_{xy}^a\wtilde V ^{\ast}(y) ] \right] \\
&\leq c(x, \wtilde \pi ^ \ast(x))  +\alpha \wtilde P^{\wtilde \pi ^ \ast}\wtilde V ^{\ast}(x) + \alpha \mid  P^{\wtilde \pi ^ \ast} \wtilde V ^{\ast}(x) - \wtilde P^{\wtilde \pi ^ \ast}\wtilde V ^{\ast}(x) \mid  +  \frac{\alpha ^2  }{1-\alpha}\delta_{\tV^{\tpi}}[P^{\tpi}, \wtilde{P}^{\tpi}]\\
& \leq \wtilde V ^{\ast}(x) + \frac{\alpha}{1-\alpha}\delta_{\tV^{\tpi}}[P^{\tpi}, \wtilde{P}^{\tpi}]= J_2(x)
\end{align*}
So we have $ J_2(x) \geq T J_2(x)$. Since ${T} J_2(x)  \to  V ^{\ast}(x)$, by monotonicity we have $  J_2(x) \geq V ^{\ast}(x)$. We conclude that 
$$\mid V^*(x) - \tV^*(x) \mid  \leq  \frac{ \alpha}{1-\alpha} (\delta_{V^{\pi^*}}[P^{\pi^*}, \wtilde{P}^{\pi^*}] + \delta_{\tV^{\tpi}}[P^{\tpi}, \wtilde{P}^{\tpi}]). $$
as stated. \eProof

\noindent {\bf Proof of Theorem \ref{thm:backtodelta2}. } 
Analogously to Theorem \ref{thm:backtodelta}, since $\delta$ can be expressed explicitly as $\delta_V[P,\tP]\leq \max_{x \in \mathcal{S}} \sum_{y}p_{xy}|V(y)-\sum_{l}\sg_{yl}V(x_l)|$, by \eqref{eq:opt_delta} 

\begin{align*}
    \mid V^*(x) - \tV^*(x) \mid_{\calS^*}  \leq & \max_{x \in \mathcal{S}} \frac{ 1}{1-\alpha} \sum_{y}\Big(p^{\pi^*}_{xy}|V^{\pi^*}(y)-\sum_{l}\sg_{yl}V^{\pi^*}(x_l)|+\wtilde{p}_{xy}^{\tpi}|\tV^{\tpi}(y)-\sum_{l}\sg_{yl}\tV^{\tpi}(x_l)|\Big)\\
    \leq & \frac{ 1}{1-\alpha}\left( |V^*-GW^*|_{\calS}^*  +|\tV^*-G\tW^*|_{\calS}^*\right)
\end{align*}
\eProof

\noindent {\bf Proof of Lemma \ref{lem:interior}.} The existence and uniqueness of a solution $\hV\in \mathcal{C}^{2,1-\vartheta}(\bar{\calB_r})$ to the Dirichlet problem follows from Assumption \ref{asum:primitives} and \cite[Theorem 6.14]{gilbarg2015elliptic}. Indeed, The assumed Lipschitz continuity of $\mu,\sigma^2$ implies, in particular, that they are both H\"{o}lder continuous for any $\vartheta\in(0,1)$ on any bounded set and, in particular, on $\calB_r$. 

Let $u$ be this solution and let us write 
$$u(x) =\frac{c(x)}{1-\alpha} - \frac{f}{1-\alpha},\mbox{ where } f:=c-(1-\alpha)u.$$ Since $c$ twice continuously differentiable, $f$ inherits its smoothness from that of $u$ as established above. In particular, 
$$D^2 u = \frac{1}{1-\alpha}(D^2c - D^2 f).$$ 
Then, 
$$|D^2 u |_{\brh}^* \leq \frac{1}{1-\alpha}\lpb|D^2c|_{\brh}^* + |D^2 f|_{\brh}^*\rpb \leq \frac{1}{1-\alpha}\Gamma \lpb \|x\|^{k-2}+\left(\frac{1}{1-\alpha}\right)^{\frac{k-2}{2}}+ |D^2 f|_{\brh}^*\rpb.$$ 

To complete the bounds we must bound $|D^2 f|_{\brh}^*$. The function $f$, notice, solves the equation \be \frac{1}{2}trace(\sigma^2(y)D^2 f(y))+\mu(y)'Df(y)-(1-\alpha)f(y) = \frac{1}{2}trace(\sigma^2(y)D^2 c(y)) + \mu(y)' D c(y).\label{eq:fequation}\ee Bounds on the derivative of $f$ will now follow from general derivative estimates for PDEs. To be self contained we quote here \cite[Theorem 6.2]{gilbarg2015elliptic}. 

Also by \cite[Page 61]{gilbarg2015elliptic}, for $\theta\in (0,1)$ and a set $\br\subset \mathbb{R}^d$, $f\in\calC^{2,\theta}(\br)$ $$|f|_{2,\theta,\br}^*:= \sup_{x\in\br}|u(x)|+\sup_{x\in\br} d_x\|Df(x)\| +\sup_{x\in\br} d_x^2\|D^2 f(x)\|  + \sup_{x,y\in\br}d_{x,y}^{2+\theta}\frac{\|D^2f(x)-D^2f(y)\|}{\|y-x\|^{\theta}}, $$ 
where $d_x= dist(x,\partial \br)\leq r$ is the distance from $x$ to the boundary $\br$ and $d_{x,y}=\min\{d_x,d_y\}$. 

With some abuse of notation, let $$\brh:=\{y:\|y-x\|\leq \varrho\}=x\pm \varrho,\mbox{ and } \brhh:=\{y:\|y-x\|\leq \varrho/2\}=x\pm \varrho/2.$$

Then,  $|f|_{2,\theta,\brh}^*\geq \sup_{y\in\brh} d_y^2\|D^2 f(y)\|\geq \sup_{y\in\brhh} d_y^2\|D^2 f(y)\|.$ For all $y\in \brhh$, notice, $d_y\geq \varrho/2$ (the distance from the boundary is greater than $\varrho/2$) so that 
$$|f|_{2,\theta,\brh}^*\geq \sup_{y\in \brhh}d_y^2\|D^2 f(y)\|\geq \frac{\varrho^2}{4}\|D^2 f\|_{\brhh}^*,$$ and, thus, that 
$$|D^2 f|_{\brhh}^*\leq \frac{4}{\varrho^2} |f|_{2,\theta,\brh}^*.$$ In this way a bound on  $|f|_{2,\theta,\brh}^*$ will produce the bound stated in Lemma \ref{lem:interior}.

\begin{theorem}[Theorem 6.2 in \cite{gilbarg2015elliptic}] Let $\Omega$ be an open subset of $\mathbb{R}^d$, and let $u\in\mathcal{C}^{2,\theta}(\Omega)$ be a bounded solution in $\Omega$ of the equation 
 $$\frac{1}{2} trace(\sigma^2(y)'D^2u(y)) + \mu(y)'Du(y) -\beta(y)u(y) = g(y)$$
where $f$ is in $\mathcal{C}^{\theta}(\Omega)$ and there are positive constants $\lambda,\Lambda$ such that the coefficients satisfy
$$ \lambda^{-1}\|\xi\|^2\geq \sum_{i,j}\xi_i\xi_j\sigma_{ij}(x)\geq \lambda \|\xi\|^2, \mbox{ for all } \xi,x \in\mathbb{R}^d, $$ and
$$ |\sigma^2|^{(0)}_{0,\theta,\Omega},|\mu|^{(1)}_{0,\theta,\Omega}, |\beta|^{(2)}_{0,\theta,\Omega}\leq \Lambda. $$ Then, 
$$|u|^*_{2,\theta,\Omega}\leq C\left(|u|_{\Omega}^* + |g|^{(2)}_{0,\theta,\Omega}\right), $$ where $C=C(d,\theta,\lambda,\Lambda)$. \end{theorem}

In this theorem we take $$g(y):= \frac{1}{2}trace(\sigma^2(y)D^2c(y))+  \mu(y)'D c(y)\mbox{ and } \beta = 1-\alpha.$$ Per our observations, with $\varrho(x)\equiv (1-\alpha)^{-1/2}$, $|\sigma^2|^{(0)}_{0,\theta,\brh},|\mu|^{(1)}_{0,\theta,\brh},|\beta|^{(2)}_{0,\theta,\brh}\leq \Gamma$ so we can take $\Lambda=\Gamma,$ to conclude that \footnote{Here, we use the following simple fact: for a Lipschitz continuous function $f$, $$\sup_{x,y}d_{x,y}^{2+\theta}\frac{|f(y)-f(x)|}{\|y-x\|^{\theta}}\leq \varrho^{2+\theta}[f]_{\brh}^*\varrho^{1-\theta}\leq \varrho[f]_{\brh}.$$} \begin{align} |D^2 f|_{\brhh}^* & \leq 4(1-\alpha) |f|_{2,\theta,\brh}^*\leq  \Gamma (1-\alpha) \left(|f |_{\brh}^* + |g|_{0,\theta,\brh}^{(2)}\right) \nonumber \\ & = \Gamma\left((1-\alpha)|f |_{\brh}^* + \sqrt{1-\alpha}|g|_{\brh}^{*}+[g]_{\brh}^{*}\right). \label{eq:intermediatebound}\end{align} 

By our assumption on $c$ and $\mu$
$$ |g|_{\brh}^*\leq |\mu|_{\brh}^*|D c|_{\brh}^*+|\sigma^2|_{\brh}^*|D^2 c|_{\brh}^*\leq \sqrt{1-\alpha}\|x\|^{k-1}+\|x\|^{k-2},$$ and 
\begin{align*} [g]_{\brh}^*& \leq |\sigma^2|_{\brh}^*[D^2c]_{\brh}^*+[\sigma^2]_{\brh}^*|D^2c|_{\brh}^*
+|\mu|_{\brh}^*[D c]_{\brh}^*+[\mu]_{\brh}^*|D c|_{\brh}^*\\ & \leq 
\Gamma \left((1-\alpha)\|x\|^{k-1} + \sqrt{1-\alpha} \|x\|^{k-2}+\|x\|^{k-3}\right),\end{align*}  so that  \be 
\sqrt{1-\alpha}|g|_{0,\brh}^{*}+[g]_{0,\brh}^{*} \leq \Gamma\left((1-\alpha)\|x\|^{k-1}+\sqrt{1-\alpha}\|x\|^{k-2}+\|x\|
^{k-3}\right). \label{eq:gbound} \ee
It remains to bound $|f|_{\brh}^*$ where, recall, $f=c-(1-\alpha)u$. We do so directly by studying a related diffusion process. Specifically, consider the process 
\be \hX_i(t) = x_i+ \int_0^t \alpha \mu_i(\hX_s)ds + \sum_{j=1}^d\int_0^t \alpha \sigma_{ij}(\hX_s)dB_j(s),\label{eq:diffeq}\ee where $B_j(\cdot)$ is a standard Brownian motion. Our requirement in Assumption \ref{asum:primitives} guarantee the existence of this process as a strong solution of this stochastic differential equation; see e.g. \cite[Theorem 5.4]{klebaner2005introduction}. It is then a standard argument that the function 
\be u(x) = \Ex_x\left[\int_0^{\tau} e^{-(1-\alpha)s} c(\hX_s)ds\right]\label{eq:udiff}, \ee where $\tau=\inf\{t\geq 0:\hX_t\in \partial \calB_{r_0^2}$\} is the solution to the PDE \eqref{eq:PDE} with the boundary condition $u(x)=0$ when $x\in\partial\calB_{r_0}$. By Ito's formula \cite[Chapter 4]{klebaner2005introduction} 
\begin{align*} c(\hX_s)&= c(x) + \sum_{i} \alpha\int_0^s\mu_i(\hX_u)c_i(\hX_u)du+ \alpha\frac{1}{2}\sum_{ij} \int_0^s\sigma_{ij}(\hX_u)c_{ij}(\hX_u)du + \sum_{ij}\int_0^s c_{i}(\hX_u)\sigma_{ij}^2(\hX_u)dB_j(u).\end{align*}

Using the boundedness of $\mu$ and $\sigma^2$ it is easily proves that, for all $x\in\calB_{r_0}$ (recall that \eqref{eq:udiff} is defined in the larger ball $\calB_{r_0^2}$), $c(x)\Ex_x[\int_{t=\tau}^{\infty}e^{-(1-\alpha)s}]\leq \epsilon$. Thus, we conclude that 
$$\left|u(x)-\frac{c(x)}{1-\alpha}\right| = \left|\Ex_x\lsb \int_0^{\tau} e^{-(1-\alpha)s} c(\hX_s)ds\rsb -\frac{c(x)}{1-\alpha}  \right|\leq \epsilon + \int_0^{\infty}e^{-(1-\alpha)s}\left( \mathbb{A}(s)+\mathbb{B}(s)+\mathbb{C}(s)ds\right),$$ where  \begin{align*} \mathbb{A}(s)&:=\Ex_x\left[\sum_{i}\left|\int_0^s \mu_i(\hX_u)c_i(\hX_u)du\right|\right],\\ \mathbb{B}(s) & := \Ex_x\left[ \sum_{i,j} \left|\int_0^s\sigma_{ij}^2(\hX_u)c_{ij}(\hX_u)du\right|\rsb, \\ \mathbb{C}(s)& :=\Ex_x\left[\sum_{ij}\left|\int_0^s c_{i}(\hX_u)\sigma_{ij}(\hX_u)dB_j(u)\right|\right].\end{align*} 
Since $|D^ic(x)|\leq \Gamma(1+\|x\|^{k-i})$ for $i=0,1,2$ and since $|\mu|,\sigma^2$ are globally bounded, we have  
\begin{align*}
   \mathbb{A}(s)\leq \Gamma \left(s+  \sum_{i} \Ex_x\lsb\int_0^s|\mu|_{\calB_{r_0^2}}^*\|\hX_u\|^{k-1}du\rsb\right),~ \mathbb{B}(s)\leq  \Gamma\left(s+\sum_{ij} \Ex_x\lsb \int_0^s|\sigma^2|_{\calB_{r_0^2}}^*\|\hX_u\|^{k-2}du\rsb \right),  
      \end{align*}  and 
            $$ \mathbb{C}(s) \leq \Gamma\left(s+ \sqrt{\Ex_x\left[\int_0^s(|\sigma|_{\calB_{r_0^2}}^*)^2\|\hX_u\|^{2(k-1)}du\right]}\right).$$ This last bound follows, again, from a standard result on Brownian integrals \cite[Theorem 4.3]{klebaner2005introduction}.
From \eqref{eq:diffeq} and the global boundedness of $\mu$ and $\sigma^2$ we have, for any $l\in\mathbb{Z}_+$ and $t$, that (recalling \eqref{eq:diffeq}) 
    \begin{align*} \Ex_x[|\hX(t)|^l]&\leq \Gamma\left( \|x\|^l+\sum_{i}\Ex_x[\int_0^t |\mu_i(\hX_{u})|du]    +\sum_{j}|\Ex_x\int_0^t \sigma_{ij}(\hX_u)dB_j(u)|^l\right]\\ &
    \leq \Gamma(|x_i|^{l} + \sqrt{1-\alpha}^lt^l + t^{l/2}),\end{align*} where we use our assumption that $|\mu|_{\calB_{r_0^2}}^*\leq \Gamma \sqrt{1-\alpha}$. Thus, 
        $$\mathbb{A}(s)\leq \int_0^s (1+\sqrt{1-\alpha}\Ex_x[\|\hX_u\|^{k-1})du\leq \Gamma(1+ \sqrt{1-\alpha}(s\|x\|^{k-1} + \sqrt{1-\alpha}^{k-1} s^{k} + s^{\frac{k+1}{2}})).$$ 
        
        We can repeat the same for $\mathbb{B}(s),\mathbb{C}(s)$. Multiplying by $(1-\alpha)$ we conclude 
    \begin{align} |c-(1-\alpha)u|&=(1-\alpha)|\frac{c}{1-\alpha}-u|\nonumber\\&\leq (1-\alpha)\epsilon + (1-\alpha)\int_0^{\infty}e^{-(1-\alpha)s}(\mathbb{A}(s)+\mathbb{B}(s)+\mathbb{C}(s))ds\nonumber  \\&  \leq \Gamma \left (\frac{\|x\|^{k-1}}{\sqrt{1-\alpha}} + \left(\frac{1}{1-\alpha}\right)^{\frac{k+1}{2}}\right).\label{eq:fbound}\end{align} 
    
    We then  have that \be \Gamma(1-\alpha)|f|_{\brh}^* \leq \Gamma\left(\sqrt{1-\alpha}\|x\|^{k-1}+ \left(\frac{1}{1-\alpha}\right)^{\frac{k-1}{2}}\right). \ee 
    
    Combining this with \eqref{eq:gbound} we have that 
     $$ \frac{|D^2f|_{\brh}^*}{1-\alpha} \leq \Gamma\left(\frac{\|x\|^{k+1}}{\sqrt{1-\alpha}}+\left(\frac{1}{1-\alpha}\right)^{\frac{k-1}{2}}\right).$$ Finally, recalling $\|D^2c\|\leq \Gamma(1+\|x\|^{k-2}),$  we also have 
                    \begin{align*} |D^2u|_{\brh}^*& = \frac{|D^2 c|_{\brh}^*}{1-\alpha}+ \frac{|D^2 f|_{\brh}^*}{1-\alpha} \\&   \leq \Gamma\left(\frac{\|x\|^{k-1}}{\sqrt{1-\alpha}}+\frac{\|x\|^{k-2}}{1-\alpha}+\left(\frac{1}{1-\alpha}\right)^{\frac{k+1}{2}}\right)\\ & 
                    \leq \Gamma\left(\frac{\|x\|^{k-1}}{\sqrt{1-\alpha}}+\left(\frac{1}{1-\alpha}\right)^{\frac{k+1}{2}}\right),\end{align*}  as stated.  \eProof

\noindent {\bf Proof of Lemma  \ref{lem:intergratedDer}.} The proof of the first part is a simpler version of that of \cite[Theorem 1]{braverman2018taylor} and we refer the reader there.   

We turn to second part. By the definition of $\Delta_x$, $\|X\|_{t+1}\leq \|X_t\|+\Delta_x\leq \Gamma(1+\|X_t\|+\sqrt{\|X_t\|}).$ In particular, given $\kappa$, there exists $m(\kappa)$ such that if $\|x\|\geq m(\kappa)$, $\|X_{t+1}\|\leq (1+\kappa)\|X_t\|.$ Overall, 
\begin{equation} \|X_{t+1}\|\leq \max\{(1+\kappa)\|X_t\|\1_{\{\|X_t\| \geq m(\kappa)\}},2m(\kappa)\}.\label{eq:jumps} \end{equation} By Assumption \ref{asum:primitives}, $|c(x)|\leq \Gamma (1+\|x\|^k)$ so that $|c(X_{t+1})| \lesssim 1+\max\{(1+\kappa)\|X_t\|\1_{\{\|X_t\| \geq m(\kappa)\}},2m(\kappa)\}^k$.
% $|c(X_{t+1})| \lesssim 1 + (\|X_t\|+\sqrt{\|X_{t}\|}^k)$
Thus, 
$$\Ex_x\left[\sum_{t=\tau(r_0^2)+1}^{\infty}\alpha^t c(X_{t})\right] \lesssim \frac{\Ex_x[\alpha^{\tau(r_0^2)}]}{1-\alpha}+ \Ex_x\left[\sum_{t=\tau(r_0^2)+1}^{\infty}\alpha^t ((1+\kappa)^t\|x\|)^k\right] .$$ Choosing $\kappa$ such that $\beta = \alpha (1+\kappa)^{k}<1$ we then have (notice that $\alpha<\beta$) 
$$\Ex_x\left[\sum_{t=\tau(r_0^2)+1}^{\infty}\alpha^t c(X_{t})\right]\lesssim  \frac{\Ex_x[\beta^{\tau(r_0^2)}]}{1-\alpha}.$$ Equation \eqref{eq:jumps} implies that, for $x\in \calB_{r_0}$ $\|X_t\|\leq (1+\kappa)^t \|x\|\leq (1+\kappa)^t r_0$ with probability $1$. In turn, $\tau(r_0^2)=\inf\{t\geq 0:X_t\notin \calB_{r_0}\}\geq  \frac{log(r)}{\log(1+\kappa)}$ with probability $1$,  so that $\Ex_x[\beta^{\tau(r_0^2)}]\downarrow 0$ as $r_0\uparrow \infty$. Choosing $r_0$ large enough then concludes the proof. \eProof

\section{m-step moment coupling \label{sec:mstep}} 

This section is an informal complement to the first comment in the concluding remarks \S \ref{sec:concluding}. There, we argued that choosing $\bg,\bu$ to match the ($m-1$)-step moment, i.e., $\sum_{z}(\bg \bu)_{yz} = \Ex_y[X_{m-1}]$, $\tP$ matches the $m$-step moment. 
$$\sum_{z} (P\bg \bu)_{yz} = \Ex_y[X_m].$$ 

A similar property holds for the second moment: if $\bg \bu$ matches the second moment $\Ex_y[X_{m-1}X_{m-1}^{\intercal}]$. then the $\tP$ matches the second moment at time $m$. We refer to this generalization as $m$-step coupling.

As an extension of \moma, one expects that a sister chain based on $m$-step coupling approximate the value of the focal chain ``sampled'' every $m$ step, namely that given $\beta\in (0,1)$:
$$V^m(x):=\Ex_x\lsb\sum_{t=0}^{\infty}\beta^t c(X_{mt})\rsb\approx\tEx_x\lsb\sum_{t=0}^{\infty}\beta^t c(X_{t})\rsb = (I-\beta \wtilde{P})^{-1}c.$$ Since what we want to eventually approximate is the value  $V(x)=\Ex_x[\sum_{t=0}^{\infty}\alpha^t c(X_t)]$ of the focal chain, the following simple relationship is useful.

\begin{lemma}\label{lem:secondstep} Consider a Markov reward process $(\calS,P,c,\alpha)$. Let $V^m(x):=\Ex_x\lsb \sum_{t=0}^{\infty}\alpha^{mt} c(X_{mt})\rsb$ and $V(x)=\Ex_x\lsb \sum_{t=0}^{\infty}\alpha^{t} c(X_{t})\rsb$. Then, 
$$V^m(x)=\frac{V(x)}{1+\sum_{k=1}^{m-1}\alpha^k}-\frac{1}{1+\sum_{k=1}^{m-1}\alpha^k}\left(\sum_{k=1}^{m-1}\alpha^k (\Ex_x[V^m(X_k)]-V^m(x)])\right).$$
 \end{lemma}

\bProof Notice that 
\begin{align*} V(x) & = \Ex\lsb \sum_{t=0}^{\infty}\alpha^t c(X_t)\rsb \\ & = \Ex\lsb\sum_{t=0}^{\infty}\alpha^{mt} c(X_{mt})\rsb + \sum_{k=1}^{m-1} \alpha^k\Ex\lsb \Ex_{X_k}\lsb\sum_{t=0}^{\infty} \alpha^{mt}c(X_{mt})\rsb\rsb \\& = V^m(x)(1+\sum_{k=1}^{m-1}\alpha^k) + \sum_{k=1}^{m-1}\alpha^k (\Ex_x[V^m(X_k)]-V^m(x)]).  
\end{align*} \eProof

\begin{figure}[h]
\includegraphics[scale=0.3]{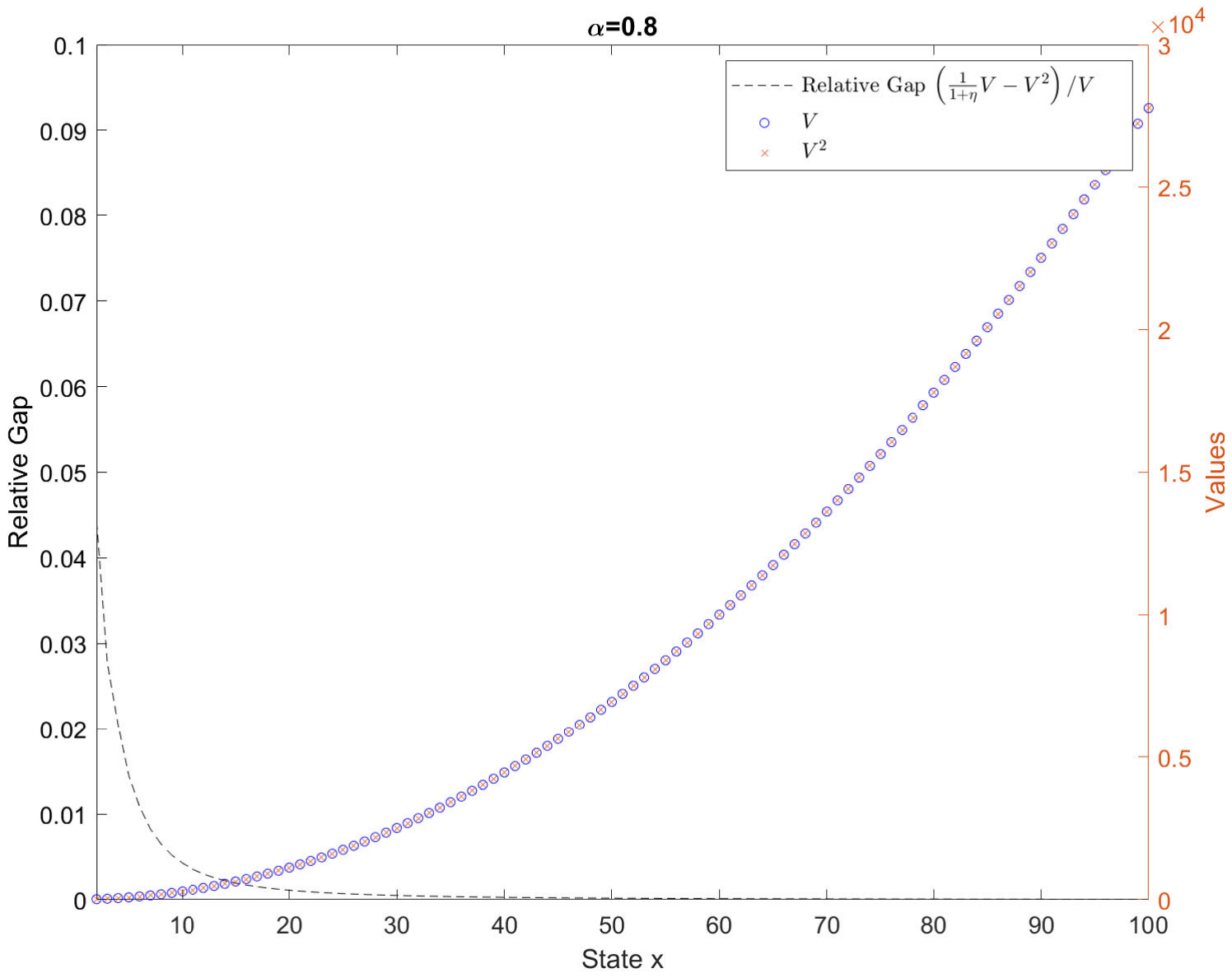}
\includegraphics[scale=0.3]{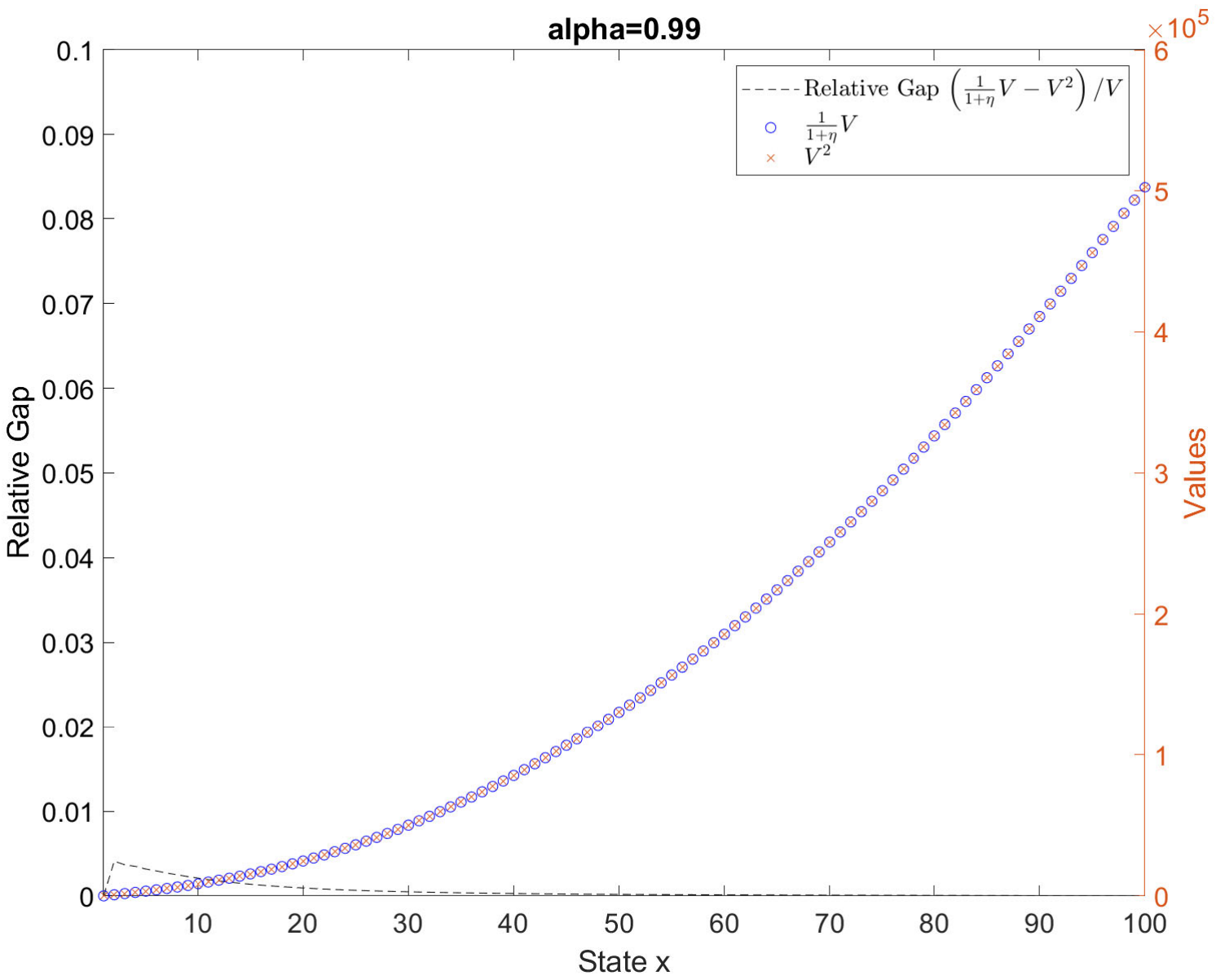}
\caption{A numerical illustration of Lemma \ref{lem:secondstep} (LEFT) $\alpha=0.8$ (RIGHT) $\alpha=0.99$. The relative error is below \%5 in the former and below \%1 in the latter. \label{fig:secondstep}}
\end{figure} 

Supposing that the chain is ergodic, we would have that the first term approaches $(1-\alpha)V(x)$ as $m\uparrow \infty$ and then $\alpha\uparrow 1$, while the second shrinks to $0$. We heuristically then take the approximation 
$$V^m(x) \approx  \frac{V(x)}{1+\sum_{k=1}^{m-1}\alpha^k}.$$

Since one expects, via moment matching, that the chain $\tP$ that matches the $m$-step moment has $\tV\approx V^m$ (notice that the discount factor for $\tV$ is $\beta=\alpha^m$) we arrive at the approximation 
$$(1+\sum_{k=1}^{m-1}\alpha^k)\tV_{\beta}(x) \approx  (1+\sum_{k=1}^{m-1}\alpha^k)V^m(x)\approx V(x).$$

A computational implementation of $m$-step coupling is not straightforward. With $1$-step coupling, each $\Ex_y[X_0]=y$ can be expressed as a convex combination of its enclosing box corners. It is no longer clear that $\Ex_y[X_1]$ can be expressed as the a convex combination of the values $\Ex_{x_l}[X_1]$ in the corner points $x_1,\ldots$. The grid has to be designed more carefully. 
 
Despite of this difficulty, the following numerical examples suggests that $m$-step coupling is a direction worth exploring.

\begin{figure} 
\includegraphics[scale=0.32]{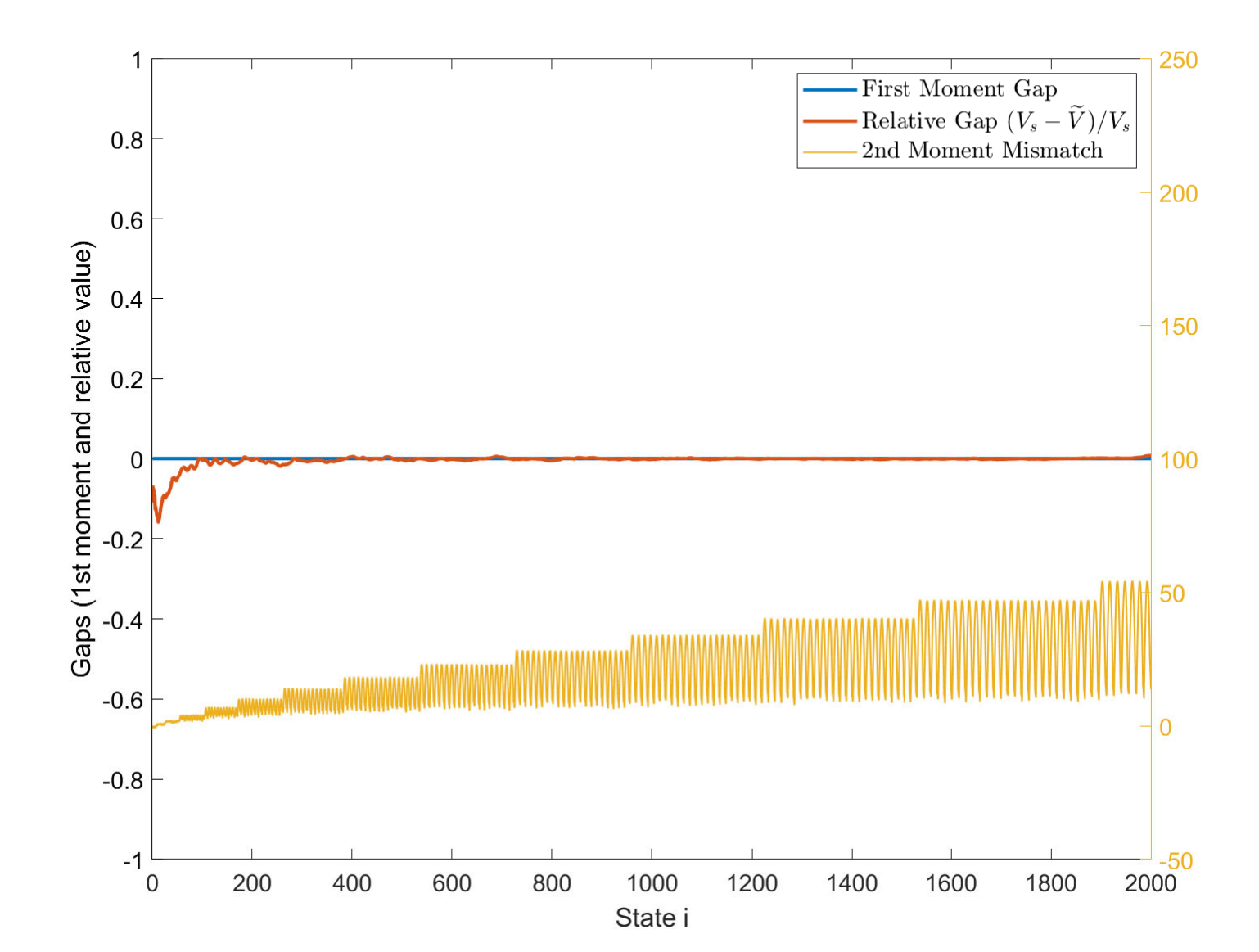}\includegraphics[scale=0.32]{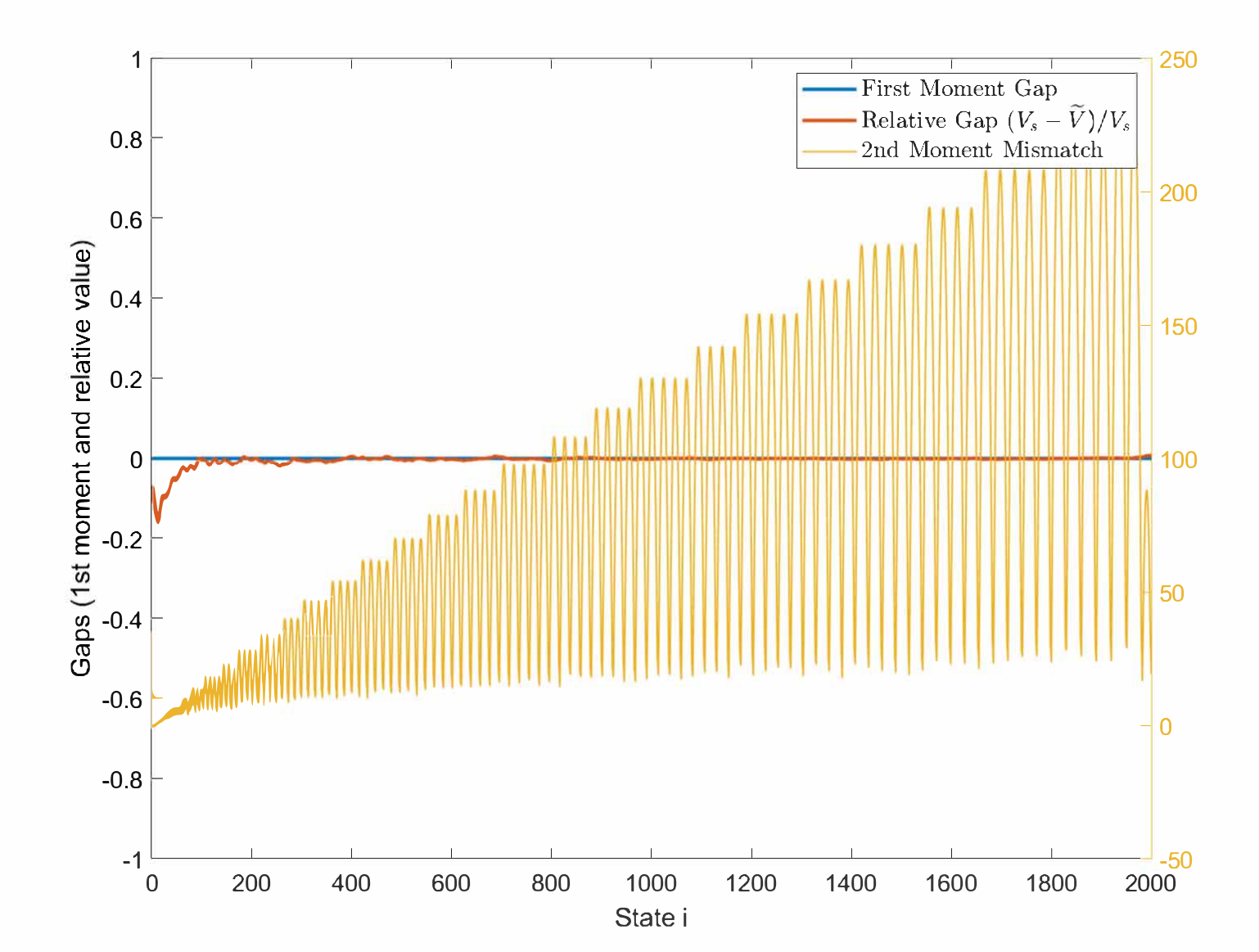}
\caption{Second moment matching, $\wtilde{P}$ vs $P^2$, with state-dependent grid (LEFT) $\fraks=0.35$ (RIGHT) $\fraks=0.45$. The growth in the latter is the larger but still sub-linear. In both case, the relative gap is smaller than $20\%$ near the origin but shrinks rapidly thereafter.  \label{fig:RW2spacing}}
\end{figure}

\begin{figure} 
\includegraphics[scale=0.3]{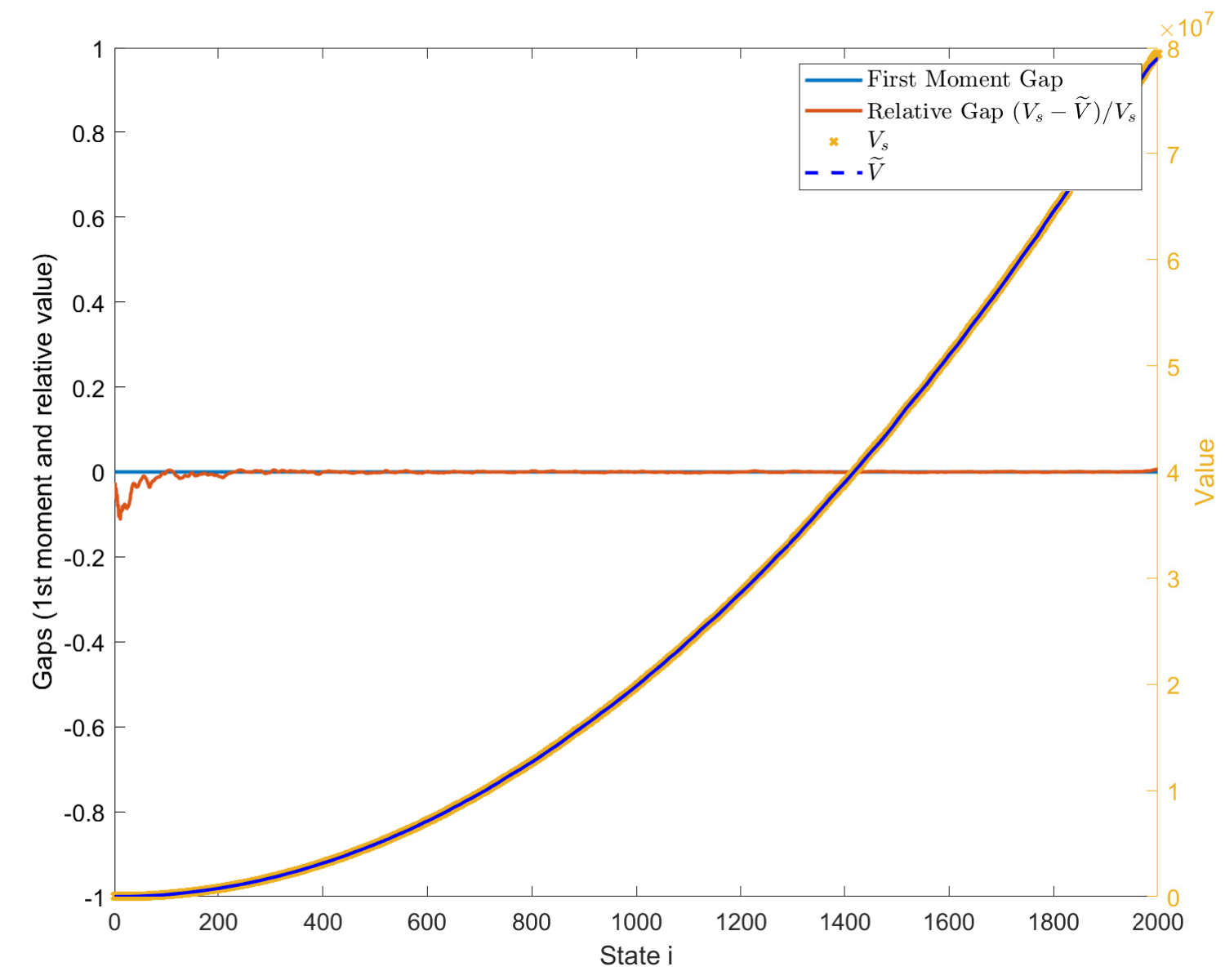}\includegraphics[scale=0.3]{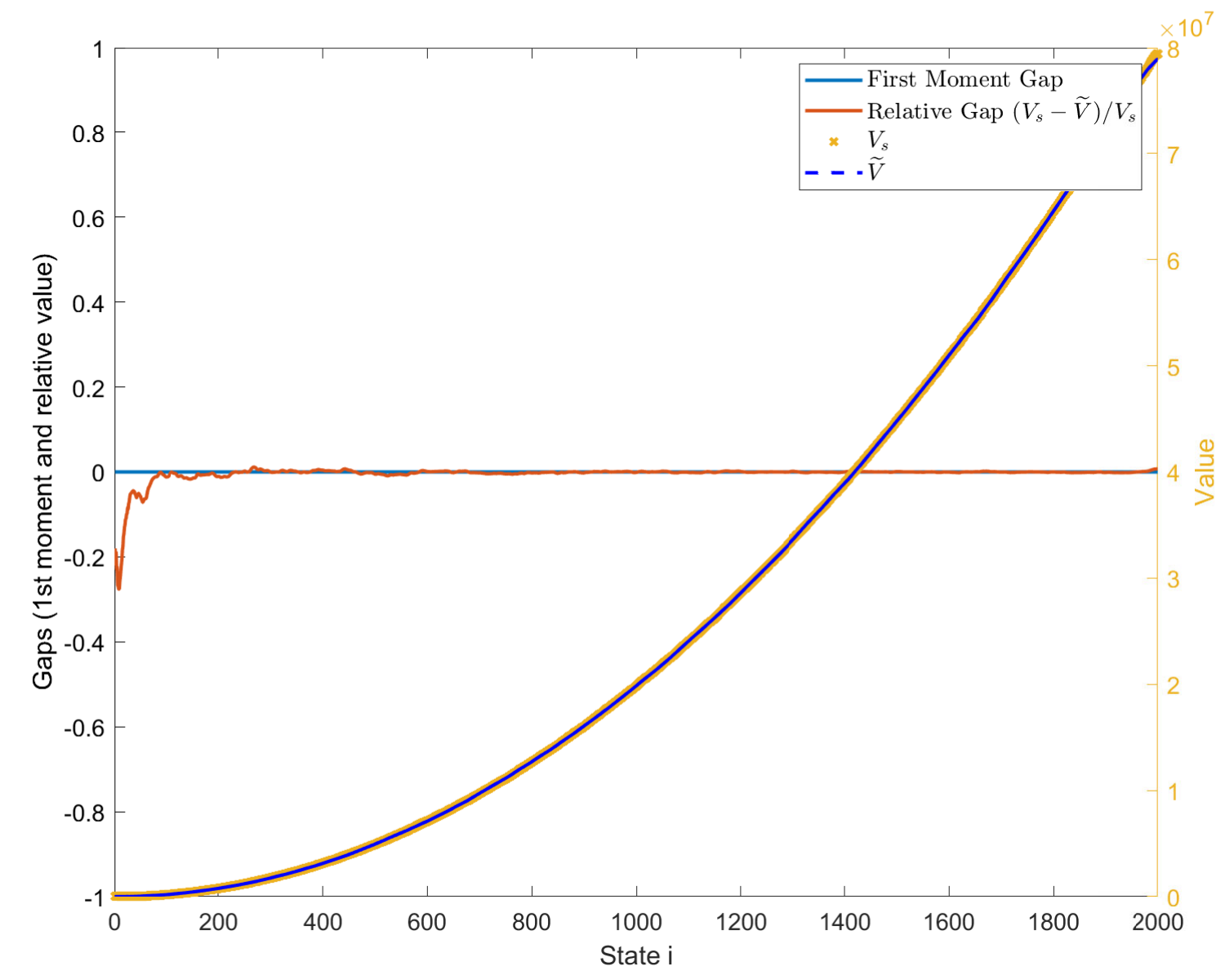}\\
\caption{Value comparison for Example \ref{fig:RW2spacing} (LEFT) $\fraks=0.35$ (RIGHT) $\fraks=0.45$, with $\calS = \{1,\ldots,2000\}$ and $\bg \bu$ constructed using the 2-step zero-mean coupling.}  \label{fig:RW2spacingvalue}
\end{figure}

\begin{figure} 
\includegraphics[scale=0.3]{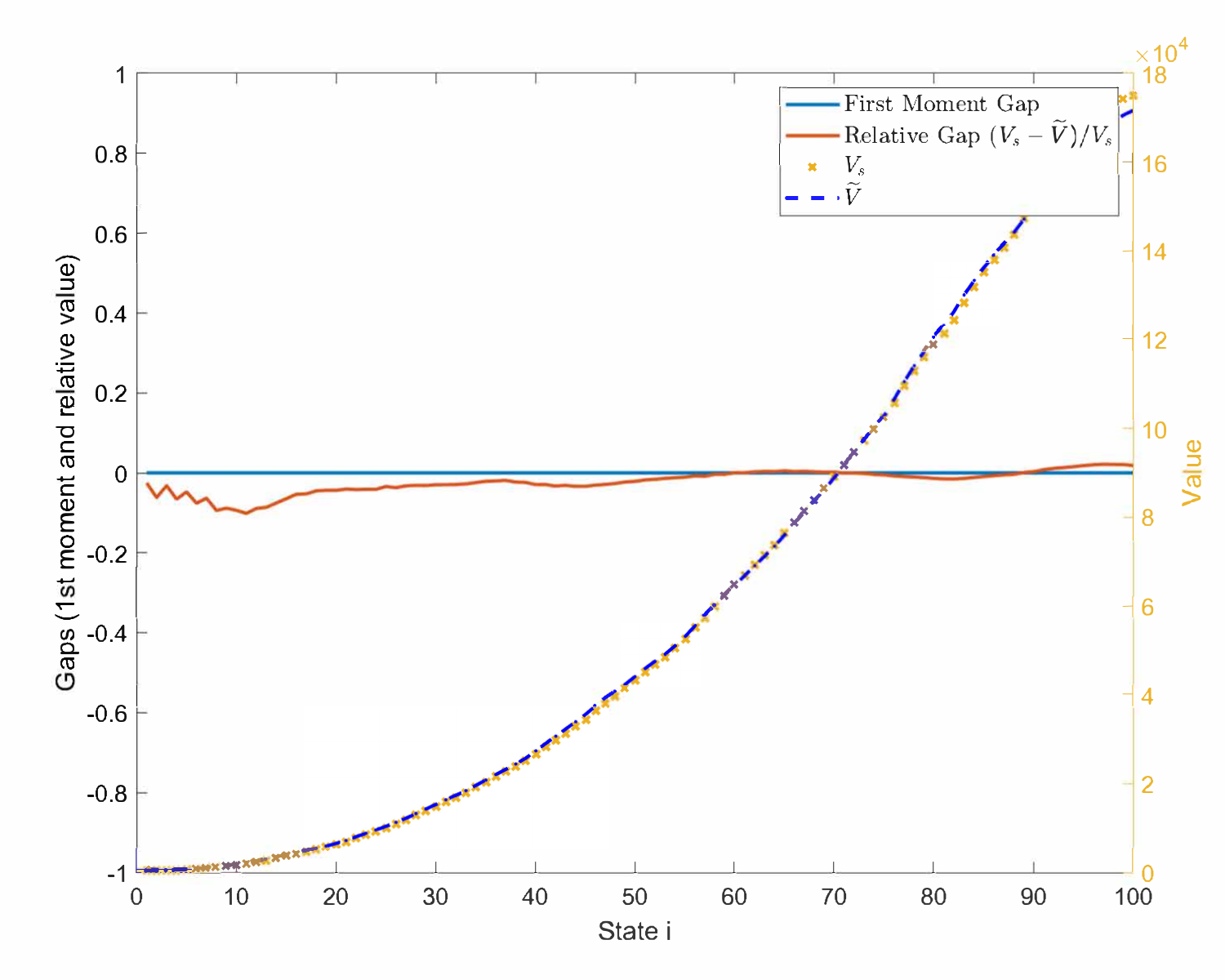}\includegraphics[scale=0.3]{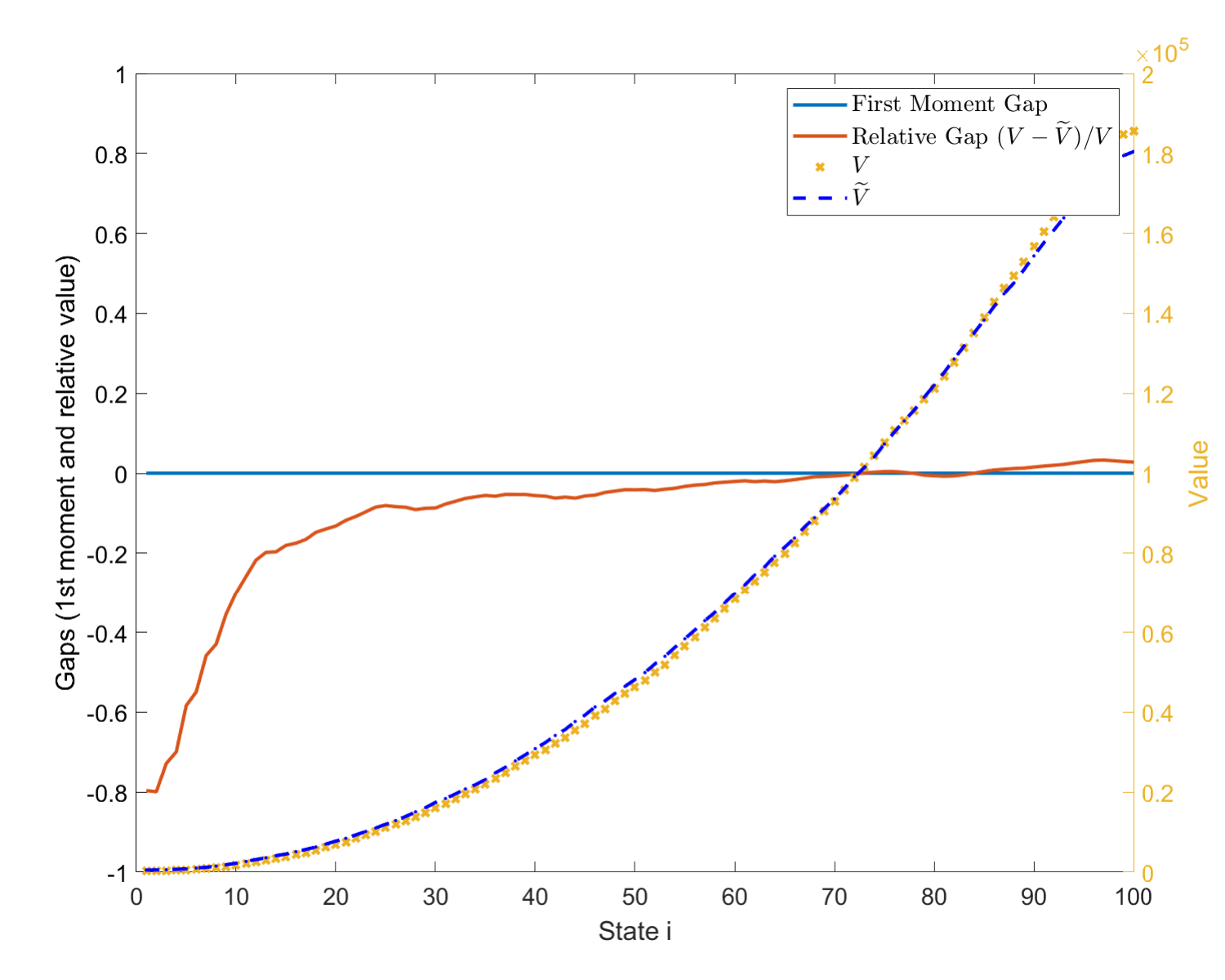}
\caption{Value comparison for Example \ref{fig:RW2spacing} for $\fraks=0.35$ and $\calS = \{1,\ldots,100\}$, with $\bg \bu$ constructed using (LEFT) the 2-step coupling vs (RIGHT) the 1-step coupling. The former shows better performance for small states.  \label{fig:RW1v2step}}
\end{figure}

\begin{figure}
\includegraphics[scale=0.3]{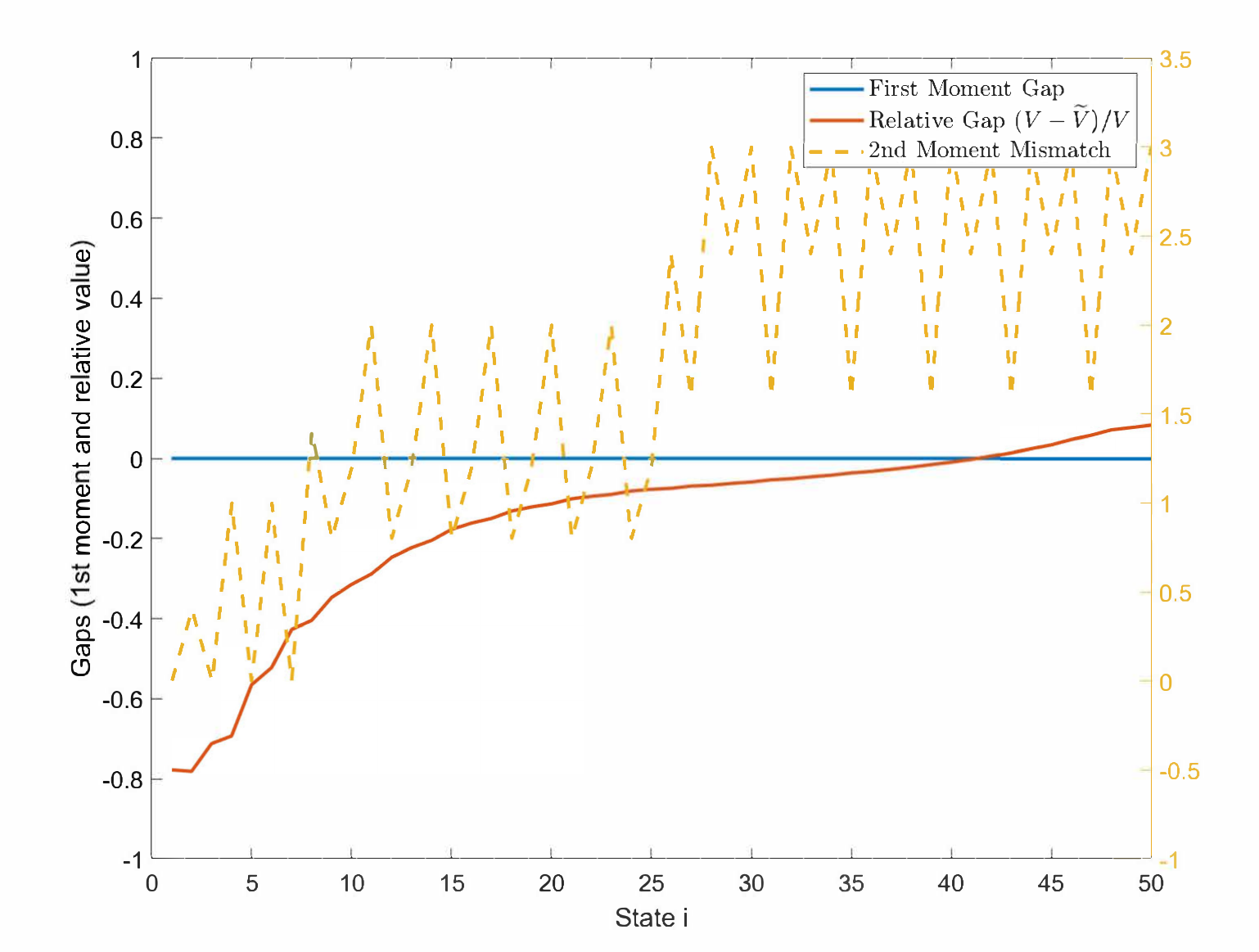}\includegraphics[scale=0.3]{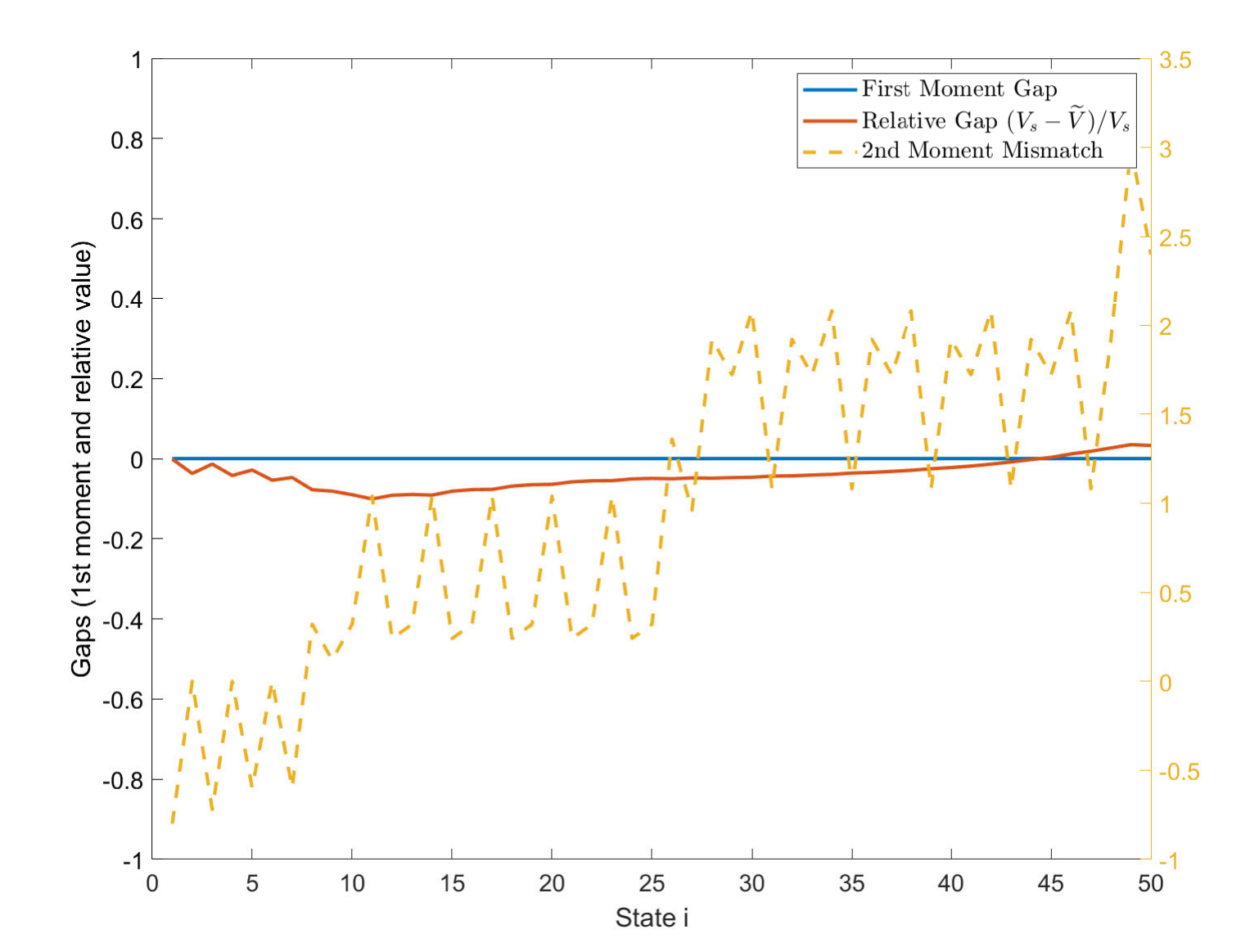}\caption{The difference in variance matching. The dashed line depicts the difference (mismatch) between $\tilde{P}M_2$ and, respectively, $PM_2$ (on the left) and $P^2M_2$ (on the right): (LEFT) 1-step coupling (RIGHT) 2-step coupling. \label{fig:variance}}\end{figure}

\begin{example} Consider a (non-absorbing) random walk on $\{1,\ldots,N\}$ with two-step coupling. We generate $\bg$ and $\tilde{P}$ for 2-step coupling, i.e., so that $\tEx_x[X_1]=\Ex_x[X_2].$ The random walk is a simple one (i.e., jumps up by one or down by 1) with ``reflecting boundaries'' $P_{12}=P_{N,N-1}=1$. Otherwise $P_{i,i+1}\in [0.5,0.6]$; the actual value was chosen as a random number  $P_{i,i+1}=0.5-0.1*rand()$. This chain then has a downward drift. 

We construct the grid based on spacing exponent $\fraks\in \{0.35,0.45\}$. Figure \ref{fig:RW2spacing} shows the growth of the second-moment (mis-) match between $\wtilde{P}$ and $P^2$. Figure \ref{fig:RW2spacingvalue} displays the value comparison $(\tilde{V}-V_s)/V_s$
where $\tilde{V}=(I-\alpha \wtilde{P})^{-1}c$ and the scaled value $V_s=\frac{1}{1+\sqrt{\alpha}}(I-\sqrt{\alpha}P)^{-1}c$ for $c(i)=i^2$ and discount $\alpha=0.95$. It also displays the first moment matching (to confirm it is 0,  as by design).  

A similar value comparison between one-step and two-step construction is illustrated in Figure \ref{fig:RW1v2step}. The one-step construction is showing inferior performance for small states $x$. This may be because the variance gap between $\tilde{P}$ and $P$ (under the one-step construction) is larger than that between $\tilde{P}$ and $P^2$ under the two-step construction; see Figure \ref{fig:variance}. Since $P^2$ has larger jumps than $P$, the two-step approximation ``suffers less'' from the coarse grid. Generally, we conjecture that the two-step approximation will be beneficial where the local variance is small, so that the gains in the better matching of the second moment overwhelm the approximation error in Lemma \ref{lem:secondstep}. \label{example:RWspacing}\hfill \bsq 
\end{example} 

\newpage 

\end{document}